\documentclass[preprint, authoryear]{elsarticle}

\usepackage{amsthm,amsmath,amssymb,natbib, subcaption}
\usepackage{graphicx, placeins}
\RequirePackage[colorlinks,citecolor=blue,urlcolor=blue]{hyperref}

\def\tilde{\widetilde}

\def\esssup{\mathop{\mathrm{ess\,sup}}}

\newcommand{\RR}{\mathbb{R}}
\newcommand{\FF}{\mathbb{F}}

\theoremstyle{localthm}
\theoremstyle{definition}
\newtheorem{Definition}{Definition}
\newtheorem{Theorem}{Theorem}
\newtheorem{Proposition}[Theorem]{Proposition}
\newtheorem{Corollary}[Theorem]{Corollary}
\newtheorem{Lemma}[Theorem]{Lemma}
\newtheorem{Example}{Example}

\newtheoremstyle{localrem}
	{5pt} 
	{5pt} 
	{\rm} 
	{} 
	{\bf} 
	{{\rm.}} 
	{.7em} 
	{} 

\begin{document}

\begin{frontmatter}

\title{Bi-$s^*$-Concave Distributions} 

\author[1]{Nilanjana Laha}
    \ead{nlaha@hsph.harvard.edu}
    
\author[2]{Zhen Miao} 
     \ead{zhenm@uw.edu} 
     
\author[3]{Jon A. Wellner\fnref{fn1}\corref{cor1}}
     \ead{jaw@stat.washington.edu}
     
 \cortext[cor1]{Corresponding author} 
 \fntext[fn1]{The research of J. A. Wellner was partially supported by NSF grant DMS-1566514, NI-AID grant 2R01 AI291968-04, 
a Simons Fellowship via the Newton Institute (INI-program STS 2018), Cambridge University, and the Saw Swee Hock Visiting 
Professorship of Statistics at the National University of Singapore (in 2019).}

\address[1]{Department of Biostatistics, Harvard University, 677 Huntington Ave, Boston, MA 02115}
\address[2]{Statistics, Box 354322, University of Washington, Seattle, WA 98195-4322}
\address[3]{Statistics, Box 354322, University of Washington, Seattle, WA 98195-4322}

\begin{abstract}
We introduce new shape-constrained classes of distribution functions on $\RR$, the bi-$s^*$-concave classes. 
In parallel to results of \cite{DUMBGEN20171} for what they called the class of bi-log-concave distribution functions, 
we show that every $s$-concave density $f$ has a bi-$s^*$-concave distribution function 
$F$ for $s^*\leq s/(s+1)$.  
 
Confidence bands building on existing nonparametric confidence bands, but accounting for the shape constraint 
of bi-$s^*$-concavity, are also considered.  
The new bands extend those developed by \cite{DUMBGEN20171} for the constraint of bi-log-concavity.  
We also make connections between bi-$s^*$-concavity and finiteness of the Cs\"org{\H o} - R\'ev\'esz constant of $F$ 
which plays an important role in the theory of quantile processes.

\end{abstract}


\end{frontmatter}
\newpage


\section{Introduction}  
Statistical methods based on shape constraints have been developing rapidly during 
the past 15 - 20 years.  From the classical univariate methods based on monotonicity 
going back to the work of \cite{MR93415} 
and \cite{MR0074746} 
in the 1950's and 1960's, research 
has progressed to consideration of convexity type constraints in a variety of problems including 
estimation of density functions, regression functions, and other ``nonparametric'' functions 
such as hazard (rate) functions.   See \cite{MR3881203} 
for a summary and overview of some of this recent activity.   

One very appealing shape constraint is log-concavity:  a (density) function $f : \RR^d \rightarrow [0,\infty]$ 
is {\sl log-concave} if $\log f$ is concave (with $\log 0 = - \infty$).  
See \cite{MR3881205}  for a recent review of the properties of  log-concave densities and their relevance
for statistical applications.  While much of the current literature has focused on point estimation,
our main focus here will be on inference for one-dimensional distribution functions and especially 
on (honest, exact) confidence bands for distribution functions which take advantage of shape constraints. 

To this end, \cite{DUMBGEN20171}   introduced the class of {\sl bi-log-concave} distribution functions defined as follows:
a distribution function $F$ on $\RR$ is bi-log-concave if both $F$ and $1-F$ are log-concave.  
They provided several different equivalent characterizations of this property, and noted (the previously known fact) 
 that if $f$ is a log-concave density, then the corresponding distribution function $F$ and survival function $1-F$ 
 are both log-concave.  But the converse is false:  there are many bi-log-concave distribution functions $F$
 with density $f$ which fail to be log-concave;  see Section 2 below for an explicit example.  
 \cite{DUMBGEN20171}  also showed how to construct confidence bands which exploit the bi-log-concave shape 
 constraint and thereby obtain narrower bands, especially in the tails, 
 with correct coverage when the bi-log-concave assumption holds.  
 
 However, a difficulty with the assumption of bi-log-concavity is that the corresponding density 
 functions inherit the requirement of exponentially decaying tails of the class of log-concave densities,
 and this rules out distribution functions $F$ with tails decaying more slowly than exponentially.   
 Here we introduce new shape-constrained families of distribution functions $F$, which we call the 
 {\sl bi-$s^*$-concave distributions},  with tails possibly decaying more slowly (or more rapidly) than
 exponentially.  As the name indicates, these families involve a parameter $s^* \in (-\infty, 1]$ 
 which allows heavier than exponential tails when $s^* < 0$, lighter than exponential tails when 
 $s^* >0$, and which correspond to exactly the bi-log-concave class introduced by \cite{DUMBGEN20171} 
 when $s^* = 0$.  
 
 Here is an outline of the rest of the paper.
 In Section~\ref{Section:Examples-Properties}  we give careful definitions of the new classes of 
 {\sl bi-$s^*$-concave distributions}.  
 We also present several helpful examples and discuss some basic properties of the new classes
 and their connections to the classes of  $s$-concave densities studied by 
 \cite{Borell1975}, 
 \cite{BRASCAMP1976366},  
 and \cite{rinott1976}.  
 (See also 
 \cite{trove.nla.gov.au/work/12919064}, and \cite{MR1898210}.)  
 Section~\ref{Section:MainTheoryResults}  contains the main theoretical results of the paper.  
 The connection between the bi-$s^*$-concave class and a key condition in the theory of 
 quantile processes, the Cs\"org{\H o} - R\'ev\'esz 
 condition,  is discussed in Corollary \ref{Corollary:CRCondition}.
Finally, we give two tail bounds for distribution functions $F\in\mathcal{P}_{s^*}$, 
see Corollary \ref{Corollary:BoundsForF}.

 In Section~\ref{Section:ConfidenceBands}   we first introduce the new confidence bands for a 
 distribution function $F \in {\cal P}_{s^*}$ assuming $s^*$ is known.  
 We also discuss some of their theoretical properties:  
 the consistency of confidence bands is discussed in Theorem~\ref{Theorem:ConsistencyOfConfidenceBands}, 
 and 
Theorem~\ref{Theorem:RateOfConfidenceBands} provides a rate of convergence for linear functionals of 
bi-$s^*$-distribution functions contained in the bands. 
This extends Theorem 5 of \cite{DUMBGEN20171}.
We then briefly discuss the algorithms used to compute 
 the new bands, and illustrate the new bands  with real and artificial data.
 Section~\ref{Section:Question}  gives a brief summary and statements of further problems. 
 An especially important remaining problem concerns construction of confidence bands when 
 $s^*$ is unknown.
The proofs for all the results in Sections~\ref{Section:Examples-Properties},   
~\ref{Section:MainTheoryResults}, and
~\ref{Section:ConfidenceBands} 
are given in Sections~\ref{Section:Proofs} and~\ref{Section:Appendix1}.
 
 We conclude this section with some notation which will be used throughout the 
 rest of the paper.
The supremum norm of a function $h:\mathbb{R}\rightarrow\mathbb{R}$ 
is denoted by $\|h\|_{\infty}\equiv\sup_{x\in\mathbb{R}}|h(x)|$, and for $K\subset\mathbb{R}$ 
we write $\|h\|_{K,\infty}\equiv \sup_{x\in K}|h(x)|$.
For a function $x\mapsto f(x)$, 
\begin{eqnarray*}
f^\prime_+(x) \equiv  \lim_{\lambda\downarrow0}\frac{f(x+\lambda)-f(x)}{\lambda},\ \ && \mbox{and} \ \ 
f^\prime_-(x) \equiv \lim_{\lambda\uparrow0}\frac{f(x+\lambda)-f(x)}{\lambda}, \\
f(x+) \equiv  \lim_{y\downarrow x}f(y),\ \  && \mbox{and} \ \ 
f(x-)\equiv\lim_{y\uparrow x}f(y),
\end{eqnarray*}
assuming that the indicated limits exist.
In general, we use $F$ and $f$ to denote a distribution function and the 
corresponding density function with respect to Lebesgue measure, and we set 
$J(F)\equiv \{x\in\RR: 0<F(x)<1\}$.

\section{Definitions, Examples, and First Properties} 
\label{Section:Examples-Properties} 
As we discussed above,  for distribution functions $F$ on $\RR$,
\cite{DUMBGEN20171} introduced a shape constraint they called {\sl bi-log-concavity} 
by requiring that both $F$ and $1-F$ be log-concave.

\indent In this paper, we generalize the bi-log-concave distribution functions by introducing 
and studying bi-$s^*$-concave distributions defined as follows.  

\par\noindent
\begin{Definition}
\label{Definition1}
\hfill \\
For $-\infty<s^*<0$, a distribution function $F$ is bi-$s^*$-concave if both 
$x\mapsto F^{s^*}(x)$ and $x\mapsto \left(1-F(x)\right)^{s^*}$ are convex functions from $\RR$ to $[0,\infty]$.\\
For $s^*=0$, a distribution function $F$ is bi-$s^*$-concave (or bi-log-concave)  if both 
$x\mapsto \log(F(x))$ and $x\mapsto \log\left(1-F(x)\right)$ are concave functions from $\RR$ to $[-\infty,0]$.\\
For $0<s^*\leq 1$, a distribution function $F$ is bi-$s^*$-concave if 
$x\mapsto F^{s^*}(x)$ is concave from $(\inf J(F),\infty)$ to $[0,1]$ and $x\mapsto \left(1-F(x)\right)^{s^*}$ 
is concave from $(-\infty,\sup J(F))$ to $[0,1]$.
\end{Definition}
\smallskip

\par\noindent
The class of bi-$s^*$-concave distribution functions is denoted by $\mathcal{P}_{s^*}$, i.e.
\[
\mathcal{P}_{s^*}\equiv\{F: F \text{ is bi-$s^*$-concave}\}.
\]

\par\noindent
\begin{Definition}  
\label{Definition2} 
(Alternative to Definition \ref{Definition1}.)\\
A distribution function $F$ is bi-$s^*$-concave if it is continuous on $\RR$ 
and satisfies the following properties on $J(F)$:\\
$\bullet$ For $-\infty<s^*<0$, 
both $x\mapsto F^{s^*}(x)$ and 
$x\mapsto \left(1-F(x)\right)^{s^*}$ are convex functions on $J(F)$.\\
$\bullet$ For $s^*=0$, 
both $x\mapsto \log(F(x))$ 
and $x\mapsto \log\left(1-F(x)\right)$ are concave functions on $J(F)$.\\
$\bullet$  For $0<s^*\leq 1$, 
both $x\mapsto F^{s^*}(x)$  
and $x\mapsto \left(1-F(x)\right)^{s^*}$ are concave functions on $J(F)$. 
\label{Definition2}
\end{Definition}

See the Appendix, Section~\ref{Section:Appendix1}, for a proof of the equivalence of Definitions 1 and 2.  
The main benefit of the second definition is that it is immediately clear that any bi-$s^*$-concave 
distribution function $F$ is continuous since continuity of $F$ is explicitly required in Definition 2.
Moreover, to verify $F\in\mathcal{P}_{s^*}$ we only need to verify the convexity 
or concavity of $F^{s^*}$ or $\left(1-F\right)^{s^*}$ on the same interval $J(F)$.

\par
Recall that a density function $f$ is $s$-concave if $f^s$ is convex for $s<0$, $f^s$ is concave 
for $s>0$, and $\log f$ is concave for $s=0$.  
Two basic properties linking $s$-concave densities 
and bi-$s^*$-concave  distribution functions are given in the following two 
propositions.   Proposition~\ref{Prop:s-concaveAndBisStar} generalizes 
the case $s=0$ as noted above, while Proposition~\ref{Prop:NestedBisStar} 
generalizes the corresponding nestedness property of the classes of $s$-concave densities; see e.g. 
\cite{trove.nla.gov.au/work/12919064}, page 86, 
and \cite{Borell1975}, page 111.   

\begin{Proposition}
\label{Prop:s-concaveAndBisStar}
Suppose a density function $f$ is $s$-concave with $s\in(-1,\infty)$.  Then 
the corresponding distribution function $F$ is bi-$s^*$-concave for all $s^*\leq s/(1+s)$.
\end{Proposition}

\begin{Proposition}
\label{Prop:NestedBisStar}
The bi-$s^*$-concave classes are nested in the following sense:
\begin{eqnarray}
\mathcal{P}_{s^*}\subset \mathcal{P}_{t^*}, \text{\ \ \ \ whenever $t^*\leq s^*\leq 1$.}
\label{NestedProperty}
\end{eqnarray}
\noindent
Moreover, the bi-$s^*$-concave classes are continuous at $s^*=0$ in the following sense:
\begin{eqnarray}
\bigcup_{s^*>0}\mathcal{P}_{s^*}=\mathcal{P}_0=\bigcap_{s^*<0}\mathcal{P}_{s^*}.
\label{ContinuityAtZero}
\end{eqnarray}
\end{Proposition}
In view of the nesting property (\ref{NestedProperty}), for each $F \in {\cal P}_{s^*}$ for some $s^*$ we define
\[
s_0^* (F) \equiv \sup \{ s^* : \ F  \ \mbox{is}\ \mbox{bi-}s^*\mbox{-concave} \} .
\]
Similarly if $f$ is $s$-concave for some $s$ we define 
\[
s_0 (f) \equiv \sup \{ s : \ f \ \mbox{is} \ s \mbox{-concave} \} .
\]
We often drop the subscript $0$ if the meaning is clear.  
For other basic properties of $s$-concave densities and bi-$s^*$-concave distribution functions, 
including results concerning closure under convolution, see 
\cite{Borell1975},
\cite{trove.nla.gov.au/work/12919064},  
 and 
\cite{MR4017135}.
\medskip

Now we introduce two important parameters, one of which will appear in connection with our characterization
of the class of bi-$s^*$-concave distribution functions in the next section and in our  
examples below.
The Cs\"org{\H o} - R\'ev\'esz constant of a bi-log-concave distribution function $F$, 
denoted by $\tilde{\gamma} (F)$, is given by 
\begin{eqnarray}
\tilde{\gamma}(F)\equiv \esssup_{x\in J(F)}F(x)(1-F(x))\frac{|f^\prime(x)|}{f^2(x)},
\label{Formula:tildegammaF}
\end{eqnarray}
 provided that $F$ is differentiable on $J(F)\equiv \{x\in\RR: 0<F(x)<1\}$ with derivative 
 $f\equiv F^\prime$ and $f$ is differentiable almost everywhere on $J(F)$ 
 with derivative $f^\prime=F^{\prime\prime}$.
Here the essential supremum is with respect to Lebesgue measure.
Alternatively (and suited for our characterization Theorem~\ref{Thm:CharacterizingThm}),
\begin{eqnarray}
\gamma (F) \equiv \esssup_{x\in J(F)}}{ \{F(x) \wedge (1-F(x)) \}\frac{|f^\prime(x)|}{f^2(x)} .
\label{Formula:gammaF}
\end{eqnarray}
Note that since $u \wedge (1-u) \le 2u(1-u) \le  2\{ u\wedge (1-u) \}$ it follows that 
$2^{-1} \gamma (F) \le \tilde{\gamma} (F) \le  \gamma (F)$, and hence finiteness of $\gamma (F)$ is equivalent 
to finiteness of $\tilde{\gamma} (F)$.  
The Cs\"org{\H o} - R\'ev\'esz constant $\tilde{\gamma}(F)$ arises in the study of quantile 
processes and transportation distances between empirical distributions and true distributions on $\RR$: see
\cite{csorgo1978}, 
\cite{MR3396731},
\cite{10.2307/3318912},   
and \cite{MR4028181}.  
It follows from the characterization Theorem  1(iv)  of DKW (2017) that $F$ is bi-log-concave  if and only 
if $\overline{\gamma} (F) \le 1$. 
We will define $\overline{\gamma} (F) \ge \gamma (F)$ and 
generalize this to the classes of bi-$s^*$-concave distribution functions 
in Section~\ref{Section:MainTheoryResults}.

Now we consider several examples of $s$-concave densities and bi-$s^*$-concave distribution functions.

\begin{Example} (Student-$t$)
\label{Ex:StudentT}
Suppose $x\mapsto f_r(x)$ is the density function of the Student-$t$ distribution with $r$ degrees of freedom defined as follows:\\
\[
f_r(x)=\frac{\Gamma((r+1)/2)}{\sqrt{\pi}\Gamma(r/2)}\left(1+\frac{x^2}{r}\right)^{-(r+1)/2}\text{\ \ \ \ for $x\in\RR$.}
\]
It is well-known (see e.g. \cite{Borell1975}) that $f_r$ is $s$-concave for any 
$s\leq -1/(1+r) = s_0 (f_r)$. Note that $s$ takes values in $(-1,0)$ since $r\in(0,\infty)$.
It follows from Proposition~\ref{Prop:s-concaveAndBisStar} that $F_r^{s^*}$ and $(1-F)^{s^*}$ are convex 
for $s^* = s/(1+s) = -1/r = s_0^* (F_r) < 0$, and hence $F_r$ is bi-$s^*$-concave for all $0 < r < \infty$.  Direct calculation 
shows that the Cs\"org{\H o} - R\'ev\'esz constant $\gamma (F_r) = 1-s^* = 1+ (1/r) \in (1,\infty)$ for $0 < r < \infty$.

In particular, this yields $\gamma(F_1)=\gamma(Cauchy)=2$.  And it suggests 
that $\gamma(F)\leq 1/(1+s)=1-s^*$ for all bi-$s^*$-concave distribution functions 
$F$ where $1/(1+s)$ varies from $1$ to $\infty$ as $s$ varies from $0$ to $-1$.
This is one of the characterizations of the  bi-$s^*$-concave class that we will prove 
in Section~\ref{Section:MainTheoryResults}.
\end{Example} 

\par\noindent
\begin{Example}  ($F_{a,b}$)
\label{Ex:Fdistrb}
Suppose that $f_{a,b}$ is the family of $F-$distributions with ``degrees of freedom" $a>0$ and $b>0$. 
(In statistical practice, if $T$ has the density $f_{a,b}$,  this would usually be denoted by 
$T\sim F_{a,b}$, where $a$ is the ``numerator degrees of freedom" and $b$ is the 
``denominator degrees of freedom".) The density is given by
\[
f_{a,b}(x)=C_{a,b}\frac{x^{b/2-1}}{(a+bx)^{(a+b)/2}}\text{ for $x\geq 0$.}
\]
(In fact, $C(a,b)=a^{a/2}b^{b/2}\text{Beta}(a/2,b/2)$, and $f_{a,b}(x)\rightarrow g_a(x)$ 
as $b\rightarrow\infty$ where $g_a$ is the Gamma density with parameters $a/2$ and $a/2$.)
It is well-known (see e.g. \cite{Borell1975}) that $f_{a,b}$ belongs to the class of $s$-concave densities, if 
$s\leq -1/(1+b/2) = s_0 (f_{a,b})$ when $a\geq 2$ and $b\geq 2$. This implies that $s\in[-1/2,0)$, and the resulting 
$s_0^*=s/(1+s)=-2/b $ is in $[-1,0)$. By Proposition \ref{Prop:s-concaveAndBisStar}, it follows 
that $F^{s^*}$ and $(1-F)^{s^*}$ are convex; i.e. $F$ is bi-$s^*$-concave.
\end{Example}

\begin{Example}  (Pareto)
\label{Ex:Pareto}
Suppose that $f_{a,b}=(a/b)(x/b)^{-(a+1)}1_{[b,\infty)}(x)$, the Pareto distribution with parameters 
$a >0$ and $b>0$. In this case, $f_{a,b}$ is $s$-concave for each $s\leq -1/(1+a)$ by noting the 
convexity of $f_{a,b}^{-1/(1+a)}=(x/b)\cdot (b/a)^{1/(1+a)}$.\\
Thus we take $s=-1/(1+a)\in(-1,0)$ for $a\in(0,\infty)$ and hence $s^*=s/(1+s)$ equals $-1/a$. 
Furthermore, it is easily seen that 
\[
CR_R(x)\equiv (1-F(x))\frac{-f^\prime(x)}{f^2(x)}=1-s^*=1+1/a \text{ for all $x> b$ }.
\]
($CR_R(\cdot)$ represents the Cs\"org{\H o} - R\'ev\'esz function in the right tail.)\\
Thus the Pareto distribution is analogous to the exponential distribution in the log-concave 
case in the sense that $x \mapsto f^s (x) = c x $  (with $c = b^{-1} (b/a)^{1/(1+a)}$) is linear.
\end{Example}

\par\noindent
\begin{Example} (Symmetrized Beta)
\label{Ex:SymmetricBeta}
Suppose that 
\[
f_r(x)=C_r (1-x^2/r)^{r/2}1_{[-\sqrt{r},\sqrt{r}]}(x),
\]
where 
\[
C_r=\Gamma((3+r)/2)/(\sqrt{\pi r}\Gamma(1+r/2))
\]
and 
$r\in(0,\infty)$.
Note that $f_r$ is an $s$-concave density with $s=2/r\in(0,\infty)$ since 
\[
f^{2/r}_r(x)=C_r^{2/r}(1-x^2/r)1_{[-\sqrt{r},\sqrt{r}]}
\]
 is concave and hence the corresponding distribution function $F_r$ is bi-$s^*$-concave with $s^*=s/(1+s)=2/(2+r)$. 
As $r\rightarrow\infty$ it is easily seen that 
\[
f_r(x)\rightarrow (2\pi)^{-1/2}\exp(-x^2/2),
\] the standard normal density. Thus $r=\infty$ corresponds to $s=0$ and $s^*=0$. 
On the other hand,
\[
g_r(x)\equiv \sqrt{r}f_r(\sqrt{r}x)=\sqrt{r}C_r(1-x^2)^{r/2}1_{[-1,1]}(x)\rightarrow 2^{-1}1_{[-1,1]}(x)
\]
as $r\rightarrow 0$. Thus $r=0$ corresponds to $s=\infty$ and $s^*=1$. 
\end{Example}
\medskip

\par\noindent
Note that just as the class of bi-log-concave distributions is considerably larger than the 
class of log-concave distributions (as shown by \cite{DUMBGEN20171}), the class of 
bi-$s^*$-concave distributions is considerably larger than the class of $s$-concave distributions. 
In particular, multimodal distributions are allowed in both the bi-log-concave and the bi-$s^*$-concave classes.

\par\noindent
\begin{Example}  (Exponential family; exponential tilt of $U(0,1)$)
\label{Ex:ExponentialFam} 
Suppose that
\[ 
f_t (x)=\exp (t x  - K(t))1_{[0,1]} (x) 
\] 
 where
\begin{equation}
\label{K(t)}
K(t) \equiv \left\{
\begin{array}{rcl}
\log (e^t-1) - \log t, & & t>0,\\
0, & & t=0,\\
\log(1-e^t)-\log(-t), & & t<0,
\end{array} \right. 
\end{equation}
for $-\infty < t < \infty$ with $K(0)\equiv 0$, and further define
$F_t (x) \equiv \int_0^x f_t (y)dy$. \\
One can verify that $f_t$ is $s$-concave only for $s\leq0$ and hence $F_t$ is bi-$s^*$-concave 
for $s^*\leq s/(1+s)\leq 0$ by Proposition~\ref{Prop:s-concaveAndBisStar}.   
However, this might not be optimal; i.e. there remains the possibility that 
$F \in {\cal P}_{s^*}$ for some $s^*>0$.  
In fact, by Theorem~\ref{Thm:CharacterizingThm}(iv) it follows that $F_t \in {\cal P}_{s^*} $ with $s^* = e^{- |t|}$.
(For an example involving a power-tilt of $U(0,1)$, see \cite{trove.nla.gov.au/work/12919064} 
(iv), page 95.) 
This also implies that the converse of Proposition \ref{Prop:s-concaveAndBisStar} does not hold here or 
in general.  
The following two examples also illustrate this point.
\end{Example}

\begin{Example} (Mixture of Gaussians shifted)
\label{Ex:GaussMixture}
(\cite{DUMBGEN20171}, page 2-3) Suppose that 
$f_{\delta} $ is the mixture $(1/2)N(-\delta,1)+(1/2)N(\delta,1)$ with $\delta>0$. 
It is well-known that $f_{\delta}$ is bimodal if $\delta>1$.  Since all $s$-concave densities
are unimodal (see e.g. \cite{trove.nla.gov.au/work/12919064} page 86), it follows that $f_{\delta}$ is not $s$-concave 
for any $\delta>1$.
\cite{DUMBGEN20171} showed (numerically) that the corresponding 
distribution $F_{\delta} $ is bi-log-concave for $\delta\leq1.34$ but not for $\delta\geq 1.35$.  
With $\delta = 1.8$ this 
 example also shows that strict inequality can occur in the second inequality in Corollary~\ref{Corollary:CRCondition} below. 
\end{Example}

\par\noindent
\begin{Example} (Mixture of shifted Student-$t$) 
\label{Ex:MixOfShiftedT}
Now suppose that $f$ is the mixture $(1/2)t_1(\cdot-\delta)+(1/2)t_1(\cdot+\delta)$ 
with $\delta>0$ where $t_r$ is the standard Student-$t$ density with $r$ degrees of freedom 
as in Example 1. 
Since $f_\delta$ is bimodal if $\delta>\delta_0 \approx 0.6$ 
and all $s$-concave densities are unimodal, 
it follows that $f_\delta$ is not $s$-concave for any $\delta> \delta_0 $. 
For values of $\delta < \delta_0$, $f_{\delta}$ is $s$-concave with $s=-1/2$, so Proposition 1 applies and shows that 
$F_{\delta} $ is bi-$s^*$-concave with $s^* = -1$.
By numerical calculation, for $\delta> \delta_0$ 
the distribution functions $F_{\delta} $
are bi-$s^*$-concave for some $s^* = s^* (\delta) \in (-\infty, 1]$ which decreases (approximately 
linearly) for large $\delta$. 
\end{Example}

\par\noindent
\begin{Example} (L\'evy with $\alpha = 1/2$) 
\label{Ex:LevyOneHalf}
This example is the completely asymmetric $\alpha-$stable (or L\'evy) law with $\alpha = 1/2$.
It gives the first passage time to the level $a>0$ for a standard Brownian motion $B$ (started at $0$ and with no drift).
See e.g. \cite{MR3930614}, pages 372 - 374.  The density is given by
\begin{eqnarray*}
f_a (t) =  \frac{a}{\sqrt{2\pi t^3}} \exp ( - a^2/ 2t) 1_{(0,\infty)} (t),
\end{eqnarray*}
and the distribution function $F_a (t) = 2P(B_t \ge a ) = 2(1-\Phi (a / \sqrt{t} )) $.  It is easily seen  
that $f_a$ is $s$-concave with $s = -2/3$, and hence $F_a$ is bi-$s^*$-concave with $s^* = -2$.  Thus 
$\gamma (F) = 3$.
\end{Example}

The following table summarizes the examples:  
\begin{table}[h]
\begin{center}
\caption{Summary of Examples 1-8}
\medskip
\begin{tabular}{| r || r | r | r | c | c | c |}  \hline \hline
Name& Example       & density                     &  d.f.            & $s$     & $s^*$                      &   $\overline{\gamma}(F)$   \\
\                       & \     &  $f$                           &  $F$           & \          &  $=s/(1+s)$             & \  $=1-s^*$         \\ \hline  \hline 
student-$t$     & 1    &  $f_r, \ r>0$              & $F_r$          & $-1/(1+r)$      &   $-1/r$         & $1+(1/r)$         \\  \hline
$F_{a,b}$         & 2    & $f_{a,b}$, $a,b>0$    & $F_{a,b}$  & $-1/(1+b/2)$   &   $-2/b$       & $1+2/b$           \\   \hline
Pareto$(a,b)$  & 3    & $f_{a,b}$, $a,b>0$    & $F_{a,b}$   & $-1/(1+a)$     &   $-1/a$       &  $1 + 1/a$        \\   \hline
Symmetric       & 4    & $f_r$, $r>0$              &  $F_r$        & $ 2/r$             &   $2/(r+2)$  &  $1/(1+2/r)$      \\ 
\ \  Beta            &  \    & \                                &  \                &  \                     &  \                 &  $= r/(r +2)$          \\  \hline
Expo family                            & 5   &  $f_t$, $t \in \RR$     & $F_t$         &  $0$                 &  $e^{-|t|}$    &  $1 - e^{- |t|}$    \\ 
Tilted $U(0,1)$                       & \    &   \                              & \                  &  \                     &  \                 &   \                       \\ \hline
Mixture,                                  & 6   & $f_{\delta} $              & $F_{\delta}$  & not $s$-         &  $0$  for       &  $1$                  \\
$N(\delta,1)$, $N(-\delta,1)$   & \    &  \                              & \                     & concave       &   $0<\delta<1.34$    &  $0<\delta<1.34$                        \\   
\                                              & \    &  \                               & \                     & for $\delta>1$  &  \      &  \               \\ \hline  
Mixture,                                  & 7   & $f_{\delta} $              & $F_{\delta}$  & not $s$-        &  bi-$s^*$-concave, \    &  $2$          \\
$T(\delta,1)$, $T(-\delta,1)$   & \    &  \                               & \                     & concave       &  some $s^*$     &  $\delta$ small                         \\      
\                                              & \    &  \                               & \                     & $\delta> .6$  &  $0<\delta<\infty$     &  \               \\ \hline  
 L\'evy $\alpha = 1/2$             & 8   &  $f_a$                       & $F_a$           & $-2/3$           & $-2$                         &  $3$          \\ \hline \hline
\end{tabular}
\end{center}
\end{table}

\par\noindent 
Example 5 shows that strict inequality can hold in the inequality $\gamma (F) \le \overline{\gamma}(F)$

\section{Main Theoretical Results} 
\label{Section:MainTheoryResults}

Here is our theorem characterizing bi-$s^*$-concave distribution functions.

\par\noindent
\begin{Theorem} 
\label{Thm:CharacterizingThm}
Let $s^*\leq 1$.
For a non-degenerate distribution function $F$, the following statements are equivalent:\\
(i) $F$ is bi-$s^*$-concave. \\
(ii) $F$ is continuous on $\RR$ and differentiable on $J(F)$ with derivative $f=F^\prime$. \\
Moreover, for $s^* \neq 0$,
\begin{eqnarray}
\label{Formula:TailBounds1}
F(y)\left\{
\begin{array}{l}
\leq F(x)\cdot \left(1+s^*\frac{f(x)}{F(x)}(y-x)\right)_+^{1/s^*}\\
\geq 1-(1-F(x))\cdot \left(1-s^*\frac{f(x)}{1-F(x)}(y-x)\right)^{1/s^*}_+
\end{array}
\right.
\end{eqnarray}
while for $s^*=0$
\begin{eqnarray}
\label{Formula:TailBounds2}
F(y)\left\{
\begin{array}{l}
\leq F(x)\cdot  \exp \left( \frac{f(x)}{F(x)}(y-x)\right) \\
\geq 1-(1-F(x))\cdot \exp \left( - \frac{f(x)}{1-F(x)}(y-x)\right)
\end{array}
\right.
\end{eqnarray}
for all $x,y\in J(F)$.\\  
(iii) $F$ is continuous on $\RR$ and differentiable on $J(F)$ with derivative $f=F^\prime$ 
such that the $s^*$-hazard function $f/(1-F)^{1-s^*}$ is non-decreasing on $J(F)$, 
and the reverse $s^*$-hazard function $f/F^{1-s^*}$ is non-increasing on $J(F)$.\\
(iv) $F$ is continuous on $\RR$ and differentiable on $J(F)$ with bounded and strictly 
positive derivative $f=F^\prime$. Furthermore, $f$ is differentiable almost everywhere 
on $J(F)$ with derivative $f^\prime=F^{\prime\prime}$ satisfying
\begin{eqnarray}
\label{Formula:fPrimeBounds}
-(1-s^*)\frac{f^2}{1-F}\leq f^\prime \leq (1-s^*) \frac{f^2}{F}\text{ almost everywhere on $J(F)$}.
\end{eqnarray}
\end{Theorem}
\medskip

\par\noindent
The following two remarks are immediately consequences of Theorem \ref{Thm:CharacterizingThm}.
See Section \ref{Section:Proofs} for a proof of Remark 1.\\
\textbf{Remark 1.} \\
(i)  The proof of Theorem 3(iv) implies that if $s^* >1$, then not both $F^{s^*}$ and $(1-F)^{s^*}$ 
can be concave. \\
(ii) If $F$ is a bi-$s^*$-concave distribution function for $0<s^*\leq 1$,
then $\inf J(F)>-\infty$ and $\sup J(F)<\infty$.\\
(iii) If $F$ is a bi-$s^*$-concave distribution function for $s^*<0$,
then it follows that
\begin{eqnarray}
\label{Formula:Remark1}
\left(0,T(F)\right)\subset\left\{t\in\mathbb{R}^+:\int |x|^t dF(x)<\infty\right\},
\end{eqnarray}
with
\begin{eqnarray}
T(F)&\equiv&
\left\{
\begin{matrix}
\infty & \text{ if }\inf J(F)>-\infty \text{ and } \sup J(F)<\infty,\\
-\frac{1}{s^*} & \text{ otherwise.}
\end{matrix}
\right.
\end{eqnarray}
\smallskip

\par\noindent
\textbf{Remark 2.}
Suppose that $F$ is a bi-$s^*$-concave distribution function,  and define
\[
T_1(F)\equiv \sup_{x\in J(F)}\frac{f}{F^{1-s^*}}(x),\ \text{and}\ 
T_2(F)\equiv \sup_{x\in J(F)}\frac{f}{(1-F)^{1-s^*}}(x).
\]
Since $f/F^{1-s^*}$ is monotonically non-increasing on $J(F)$, 
it follows that for any $x,x_0\in J(F)$ with $x<x_0$,
\[
\frac{f}{F^{1-s^*}}(x)\geq \frac{\frac{1}{s^*}F^{s^*}(x)-\frac{1}{s^*}F^{s^*}(x_0)}{x-x_0}
\]
and hence
\begin{eqnarray*}
T_1(F)
=\sup_{x\in J(F)}\frac{f}{F^{1-s^*}}(x)
=\lim_{x\rightarrow \inf J(F)}\frac{f}{F^{1-s^*}}(x)
\left\{
\begin{matrix}
>0, &\\
=\infty &\text{ if }\inf J(F)>-\infty.
\end{matrix}
\right.
\end{eqnarray*}
Analogously one can show that 
\begin{eqnarray*}
T_2(F)
\left\{
\begin{matrix}
>0, &\\
=\infty &\text{ if }\sup J(F)<\infty.
\end{matrix}
\right.
\end{eqnarray*}

\begin{Corollary}
\label{Corollary:CRCondition}
(Connection with the Cs\"org{\H o} - R\'ev\'esz constant.)\\
Suppose $F$ is a bi-$s^*$-concave distribution function for $s^*\leq1$. 
Then with $\tilde{\gamma}(F)$ and $\gamma (F)$ as defined in (\ref{Formula:tildegammaF}) 
and (\ref{Formula:gammaF}), we have
\begin{eqnarray}
\label{Formula:InequalityOfCRCondition}
\frac{1}{2} \gamma(F) \le  \tilde{\gamma}(F) \le \gamma (F) \le \overline{\gamma} (F)\le 1-s^*,
\end{eqnarray}
where
\[
\overline{\gamma}(F)\equiv\max\{\tilde{CR}(F),\tilde{CR}(\bar{F})\},\ \bar{F}\equiv 1-F, 
\]
\[
\tilde{CR}(F)\equiv\esssup_{x\in J(F)}\frac{F(x)F^{\prime\prime}(x)}{\left(F^\prime(x)\right)^2},
\]
and
\[ 
\gamma (F) \equiv \esssup_{x\in J(F)}\frac{\{ F(x) \wedge (1-F(x))\} | F^{\prime\prime}(x)|}{\left(F^\prime(x)\right)^2 } 
= \esssup_{x\in J(F)}\frac{\{ F(x) \wedge (1-F(x))\} | f^{\prime}(x)|}{\left(f(x)\right)^2 } .
\]
\end{Corollary}

\par\noindent
\par\noindent
\textbf{Remark 3.}
By Theorem \ref{Thm:CharacterizingThm}, one can verify that $\tilde{CR}(F)$ is well-defined for any $F\in\mathcal{P}_{s^*}$.
Note that
\[ 
\tilde{CR}(\overline{F} )\equiv\esssup_{x\in J(F)}\frac{\overline{F}(x)(-F^{\prime\prime}(x))}{\left(F^\prime(x)\right)^2}.
\]
The first two inequalities in Corollary~\ref{Corollary:CRCondition} follow (as we noted before) from 
$2^{-1} \{u \wedge (1-u)\} \le u(1-u) \le  u \wedge (1-u) $ for $0 \le u \le 1$.
Thus finiteness of $\tilde{\gamma} (F)$ implies finiteness of $\gamma (F)$ and vice-versa.   
Examples show that strict inequality may hold in the inner inequalities in (\ref{Formula:InequalityOfCRCondition}).
On the other hand, if $f$ is non-decreasing on $(a,F^{-1} (1/2))$ and $f$ is non-increasing on $(F^{-1} (1/2), b)$ 
where $J(F) = (a,b)$, then $\gamma = \overline{\gamma} $ 
by inspection of the proof of $\gamma (F) \le \overline{\gamma} (F)$.
\begin{Corollary}
\label{Corollary:BoundsForF}
(Bounds for $F\in\mathcal{P}_{s^*}$, where $s^*\neq0$.)\\
For any $s^*\in(-\infty,0)\cup(0,1]$ and $F\in\mathcal{P}_{s^*}$, 
\begin{eqnarray}
F_L(x)\leq F(x)\leq F_U(x),
\label{Corollary4:Formula1}
\end{eqnarray}
where $F_L(x)\equiv \frac{1}{s^*}\left(F^{s^*}(x)-(1-s^*)\right)$ and $F_U(x)\equiv \frac{1}{s^*}\left(1-\left(1-F(x)\right)^{s^*}\right)$.\\
Moreover, $F_U(x)$ is a convex function on $J(F)$,  and $F_L(x)$ is a concave function on $J(F)$.
For $s^* = 0$ and $F \in {\cal P}_0$, (\ref{Corollary4:Formula1}) holds with 
$F_L (x) = 1 + \log F(x)$ and $F_U (x) = - \log (1-F(x))$.
\end{Corollary}

\section{Confidence bands for bi-$s^*$-concave distribution functions}
\label{Section:ConfidenceBands}

Our goal in this section is to define confidence bands for $F$ which exploit the shape constraint $F \in {\cal P}_{s_0^*}$.
We start with some known unconstrained nonparametric bands and define new bands under the assumption
that the true distribution function $F$ satisfies the shape constraint $F \in {\cal P}_{s_0^*}  $ 
where $s_0^*$ is known.    

\subsection{Definitions and Basic Properties}
\label{subsec:ConfBand-Defns-BasicProp} 
Let $X_1,\ldots,X_n$ be i.i.d. random variables with continuous distribution function $F$. 
A $(1-\alpha)$-confidence band, denoted by $(L_n,U_n)$, for $F$ means that both $L_n$ 
and $U_n$ are monotonically non-decreasing functions on $\RR$ depending on $\alpha$ and $X_1,\ldots,X_n$ only,
moreover, $L_n$ and $U_n$ have to satisfy $L_n<1$, $U_n>0$ and 
\[
P\left(L_n(x)\leq F(x)\leq U_n(x)\text{ for all }x\in\RR\right)=1-\alpha.
\]
The following two bands are discussed in \cite{DUMBGEN20171} and we briefly restate them here.
\smallskip

\par\noindent
\textbf{Example} (Komogorov-Smirnov band).
A Komogorov-Smirnov band $(L_n,U_n)$ is given by
\[
[L_n(x),U_n(x)]\equiv 
\left[\FF_n(x)-\frac{\kappa_{\alpha,n}^{KS}}{\sqrt{n}},
       \FF_n(x)+\frac{\kappa_{\alpha,n}^{KS}}{\sqrt{n}}\right]\cap [0,1],
\]
where $\FF_n$ is the empirical distribution function and $\kappa_{\alpha,n}^{KS}$ denotes the 
$(1-\alpha)$-quantile of $\sup_{x\in\RR}n^{1/2}|\FF_n(x)-F(x)|$, see \cite{MR3396731}
Note that $\kappa_{\alpha,n}^{KS}\leq \sqrt{\log(2/\alpha)/2}$ by Massart's (1990) inequality, see \cite{massart1990}.
\medskip

\par\noindent
\textbf{Example} (Weighted Komogorov-Smirnov band).
A Weighted Komogorov-Smirnov band $(L_n,U_n)$ is as follows: for any $\gamma\in[0,1/2)$,
\[
[L_n(x),U_n(x)]\equiv \left[t_i-\frac{\kappa_{\alpha,n}^{WKS}}{\sqrt{n}}(t_i(1-t_i))^\gamma,t_{i+1}
+\frac{\kappa_{\alpha,n}^{WKS}}{\sqrt{n}}(t_{i+1}(1-t_{i+1}))^\gamma\right]\cap [0,1],
\]
for $i\in\{0,1,\ldots,n\}$ and $x\in[X_{(i)},X_{(i+1)})$, 
where $\{X_{(i)}\}_{i=1}^n$ denotes the order statistics of $\{X_i\}_{i=1}^n$, 
$X_{(0)}\equiv -\infty$, $X_{(n+1)}\equiv\infty$, $t_i\equiv i/(n+1)$ for $i=1,\ldots,n$, 
and $\kappa_{\alpha,n}^{WKS}$ denotes the $(1-\alpha)$-quantile of the following test statistics 
\[
\sqrt{n}\max_{i=1,\ldots,n}\frac{|F(X_{(i)})-t_i|}{(t_i(1-t_i))^\gamma}.
\]
Note that $\kappa_{\alpha,n}^{WKS}=O(1)$.
\medskip

\par\noindent
A further example of a nonparametric confidence band due to 
\cite{doi:10.1080/01621459.1995.10476543}  and refined 
by \cite{DuembgenAndWellner2014} was considered by \cite{DUMBGEN20171}.   
We will not consider this third possibility further here due to space constraints.
\medskip

\par\noindent
Now we turn to confidence bands for bi-$s^*$-distribution functions.  
Our approach will be to refine the three unconstrained bands given in the three examples.\\
Suppose $F$ is a bi-$s^*$-concave distribution function.
A nonparametric ($1-\alpha$) confidence band $(L_n,U_n)$ for $F$ may be refined as follows:
\begin{eqnarray*}
L_n^o(x)&\equiv& \inf\{G(x):G\in\mathcal{P}_{s^*}, L_n\leq G\leq U_n\},\\
U_n^o(x)&\equiv& \sup\{G(x):G\in\mathcal{P}_{s^*}, L_n\leq G\leq U_n\}.
\end{eqnarray*}
If there is no bi-$s^*$-concave distribution function $F$ fitting into the band $(L_n,U_n)$, 
we set $L_n^o\equiv 1$ and $U_n^o\equiv 0$ and we conclude with confidence 
$1-\alpha$ that $F$ is not bi-$s^*$-concave.
But in the case of $F\in\mathcal{P}_{s^*}$, this happens with probability at most $\alpha$.
\medskip

\par\noindent
The following lemma implies two properties of our shape-constrained band $(L_n^o,U_n^o)$.
The first one is that both $L_n^o$ and $U_n^o$ are Lipschitz continuous on 
$\RR$, unless $\inf\{x\in\RR:L_n(x)>0\}\geq\sup\{x\in\RR:U_n(x)<1\}$.
The second one is that $L_n^o(x)$ converges polynomially fast to $0$ as 
$x\rightarrow-\infty$ and $U_n^o(x)$ converges polynomially fast to $1$ as $x\rightarrow\infty$
as long as $\lim_{x\rightarrow\infty}L_n(x)>\lim_{x\rightarrow-\infty}U_n(x)$.
\begin{Lemma}
\label{Lemma:LemmaInConfidenceBands}
For real numbers $a<b$, $0<u<v<1$ and $s^*\in(-\infty,0)\cup(0,1]$, define
\[
\gamma_1\equiv \frac{\frac{1}{s^*}\left(v^{s^*}-u^{s^*}\right)}{b-a}
\text{ and }
\gamma_2\equiv \frac{\frac{-1}{s^*}\left((1-v)^{s^*}-(1-u)^{s^*}\right)}{b-a}.
\]
(i) If $L_n(a)\geq u$ and $U_n(b)\leq v$, then $L_n^o$ and $U_n^o$ 
are Lipschitz-continuous on $\mathbb{R}$ with Lipschitz constant $\max\{\gamma_1,\gamma_2\}$.\\
(ii) If $U_n(a)\leq u$ and $L_n(b)\geq v$, then 
\[
U_n^o(x)\leq\left(u^{s^*}+s^*\gamma_1(x-a)\right)^{1/s^*}_+ \text{ for }x\leq a,
\]
\[
1-L_n^o(x)\leq\left((1-v)^{s^*}-s^*\gamma_2(x-b)\right)_+^{1/s^*} \text{ for }x\geq b.
\]
\end{Lemma}
\medskip

\par\noindent
The following theorem implies the consistency of our proposed confidence band $(L_n^o,U_n^o)$.
\begin{Theorem}
\label{Theorem:ConsistencyOfConfidenceBands}
Suppose that the original confidence band $(L_n,U_n)$ is consistent in the sense 
that for any fixed $x\in\mathbb{R}$, both $L_n(x)$ and $U_n(x)$ tend to $F(x)$ in probability.\\
(i) Suppose that $F\notin\mathcal{P}_{s^*}$. Then $P(L_n^o\leq U_n^o)\rightarrow0$.\\
(ii) Suppose that $F\in\mathcal{P}_{s^*}$ with $s^*\neq0$. Then $P(L_n^o\leq U_n^o)\geq 1-\alpha$, and 
\begin{eqnarray}
\label{Formula:ConfidenceBand1}
\sup_{G\in\mathcal{P}_{s^*}:L_n\leq G\leq U_n}\|G-F\|_{\infty}\rightarrow_p0,
\end{eqnarray}
where $\sup(\emptyset)\equiv0$. Moreover, for any compact interval $K\subset J(F)$,
\begin{eqnarray}
\label{Formula:ConfidenceBand2}
\sup_{G\in\mathcal{P}_{s^*}:L_n\leq G\leq U_n}\|h_G-h_F\|_{K,\infty}\rightarrow_p0,
\end{eqnarray}
where $h_G$ stands for any of the three functions $G^\prime,\left(G^{s^*}\right)^\prime$, 
and $\left((1-G)^{s^*}\right)^\prime$. Finally, for any fixed $x_1\in J(F)$ and $0<b_1<f(x_1)/F^{1-s^*}(x_1)$,
\begin{eqnarray}
\label{Formula:ConfidenceBand3}
P\left(U_n^o(x)\leq \left(U_n^{s^*}(x^\prime)+s^*b_1(x-x^\prime)\right)^{1/s^*}_+\text{ for }x\leq x^\prime\leq x_1\right)\rightarrow1,
\end{eqnarray}
while for any fixed $x_2\in J(F)$ and $0<b_2<f(x_2)/(1-F(x_2))^{1-s^*}$,
\begin{eqnarray}
\label{Formula:ConfidenceBand4}
P \left (1-L_n^o(x)\leq \left((1-L_n(x^\prime))^{s^*}-s^*b_2(x-x^\prime)\right)^{1/s^*}_+
\text{ for }x\geq x^\prime\geq x_2 \right ) \rightarrow1.
\end{eqnarray}
\end{Theorem}
\medskip

\par\noindent
The following result provides the consistency of confidence bands for 
functionals $\int \phi dF$ of $F$ with well-behaved integrands $\phi:\RR\rightarrow\RR$.
\begin{Corollary}
\label{Corollary:CorollaryInConfidenceBands}
Suppose that the original confidence band $(L_n,U_n)$ is consistent, and let $F\in\mathcal{P}_{s^*}$ with $s^*<0$.
Let $\phi:\mathbb{R}\mapsto\mathbb{R}$ be absolutely continuous with a 
continuous derivative $\phi^\prime$ satisfying the following constraint: there exist constants $a>0$ and $k<-1/s^*$ such that\\
\[
|\phi^\prime(x)|\leq a|x|^{k-1}.
\]
Then
\[
\sup_{G:L_n^o\leq G\leq U_n^o}\left|\int\phi dG-\int \phi dF\right|\rightarrow_p0.
\]
\end{Corollary}
\medskip

\par\noindent
The following theorem provides rates of convergence, with the following condition on the original confidence band $(L_n,U_n)$:
\smallskip

\par\noindent
Condition (*): For certain constants $\gamma\in[0,1/2)$ and $\kappa,\lambda>0$, 
\[
\max\{\FF_n-L_n,U_n-\FF_n\}\leq \kappa n^{-1/2}\left(\FF_n(1-\FF_n)\right)^\gamma
\]
on the interval $\{\lambda n^{-1/(2-2\gamma)}\leq \FF_n\leq 1-\lambda n^{-1/(2-2\gamma)}\}$.
\smallskip

\par\noindent
As stated in \cite{DUMBGEN20171}, this condition is satisfied 
with $\gamma=0$ in the case of the Kolmogorov-Smirnov band.
In the case of the weighted Kolmogorov-Smirnov band, it is satisfied for the given value of $\gamma\in[0,1/2)$.
For the refined version of Owen's band, it is satisfied for any fixed number $\gamma\in(0,1/2)$.
\begin{Theorem}
\label{Theorem:RateOfConfidenceBands}
Suppose that $F\in\mathcal{P}_{s^*}$ with $s^*<0$ and let $(L_n,U_n)$ satisfy Condition (*). 
Let $\phi:\mathbb{R}\mapsto\mathbb{R}$ be absolutely continuous with a continuous derivative $\phi^\prime$.\\
Suppose that $|\phi^\prime(x)|=O(|x|^{k-1})$ as $|x|\rightarrow\infty$ for some numbers $k<-1/s^*$. 
Then 
\begin{eqnarray}
\label{Formula:Theorem9}
\sup_{G:L_n^o\leq G\leq U_n^o}\left|\int \phi dG-\int \phi dF\right|
=
O_p \left(n^{-\frac{1}{2}\left(1\wedge \frac{ks^*+1}{1-\gamma}\right)}\right).
\end{eqnarray}
\end{Theorem}
\noindent
\textbf{Remark: } 
(i) From (\ref{Formula:Theorem9}), one can verify that the convergence 
rate is $n^{-1/2}$ as long as $k<\gamma/(-s^*)$.\\
(ii) From (\ref{Formula:Theorem9}), one can verify that when $\gamma=0$, 
the convergence rate is $n^{-1/2+k/(-s^*)}$ and we have a 
``power deficit" (or ``polynomial rate deficit") relative to $n^{-1/2}$.

\subsection{Implementation and illustration of the confidence bands}
\label{subsec:Implementation}

In this section, we discuss the implementation of confidence bands for bi-$s^*$-concave distribution functions. 
This extends the treatment of \cite{DUMBGEN20171} from $s^*=0$ to general values $s^*\in(-\infty,1]$.

Recall the procedure ConcInt$(\cdot, \cdot)$ developed in \cite{DUMBGEN20171}.
Given any finite set $\mathcal{T} = \{t_0,...,t_m\}$ of real numbers 
$t_0 < t_1 < \cdots < t_m$ and any pair $(l,u)$ of functions $l,u: \mathcal{T} 
\rightarrow [-\infty,\infty)$ with $l<u$ pointwise and $l(t) > -\infty$ for at least two 
different points $t\in \mathcal{T}$ , this procedure computes the pair $(l^o, u^o)$ where
\begin{eqnarray*}
l^o(x)&\equiv&\inf\left\{g(x):\text{ $g$ is concave on $\mathbb{R}$}, \  l\leq g\leq u\text{ on }\mathcal{T}\right\},\\
u^o(x)&\equiv&\sup\left\{g(x):\text{ $g$ is concave on $\mathbb{R}$},\  l\leq g\leq u\text{ on }\mathcal{T}\right\}.
\end{eqnarray*}
First note that $l^o$ is the smallest concave majorant of $l$ on $\mathcal{T}$; thus it may 
be computed by a version of the pool-adjacent-violators algorithm; see for example \cite{Robertson1988}.
Then we obtain indices $0 \leq j(0) < j(1) < \cdots < j(b) \leq m$ such that
\begin{eqnarray*}
l^o
\left\{
\begin{aligned}
&\equiv-\infty \text{ on }\mathbb{R}\backslash [t_{j(0)},t_{j(b)}],\\
&\text{is linear on } [t_{j(a-1)},t_{j(a)}] \text{ for }1\leq a\leq b,\\
&\text{change slope at }t_{j(a)} \text{ if }1\leq a\leq b.
\end{aligned}
\right.
\end{eqnarray*}
With $l^o$ in hand, we then check to see if $l^o \leq u$ on $\mathcal{T}$. 
If this fails, then there is no concave function lying between $l$ and $u$, and the procedure returns an error message. 
If this test succeeds, then we compute $u^o(x)$ as
$$
\min\left\{u(s)+\frac{u(s)-l^o(r)}{s-r}(x-s): r\in\mathcal{T}_o, r<s\leq x\text{ or }x\leq s<r\right\},
$$
where $\mathcal{T}_o=\{t_{j(0)},t_{j(1)},\ldots,t_{j(b)}\}$.
(The rest of the description of the procedure ConcInt$(\cdot,\cdot)$ is just as in \cite{DUMBGEN20171}.)

When $s^* < 0$, let $g(v;s^*)\equiv g(v)\equiv -v^{s^*}$ and 
$h(v;s^*)\equiv h(v)\equiv(-v)^{1/s^*}$. (This is the most important new case. 
When $s=s^*=0$, $g(v)\equiv \log(v), h(v)\equiv \exp(v)$. When $s^*> 0$, $g(v) \equiv v^{s^*}$ 
and $h(v)\equiv v^{1/s^*}$.) Here is pseudocode for the computation of $(L_n^o,U_n^o)$.
\begin{align*}
&(L_n^o, U_n^o)\leftarrow (L_n,U_n)\\
&(l^o,u^o)\leftarrow \text{ConcInt}(g(L_n^o),g(U_n^o))\\
&(\tilde{L}_n^o,\tilde{U}_n^o)\leftarrow (h(l^o),h(u^o))\\
&(l^o,u^o)\leftarrow \text{ConcInt}(g(1-\tilde{U}_n^o),g(1-\tilde{L}_n^o))\\
&(\tilde{L}_n^o,\tilde{U}_n^o)\leftarrow (1-h(u^o),1-h(l^o))\\
&\text{while }(\tilde{L}_n^o,\tilde{U}_n^o)\neq (L_n^o,U_n^o)\text{ do}\\
&\ \ \ \ (L_n^o, U_n^o)\leftarrow (\tilde{L}^o_n,\tilde{U}^o_n)\\
&\ \ \ \ (l^o, u^o)\leftarrow \text{ConcInt}(g(L^o_n),g(U^o_n))\\
&\ \ \ \ (\tilde{L}^o_n,\tilde{U}^o_n)\leftarrow (h(l^o),h(u^o))\\
&\ \ \ \ (l^o,u^o)\leftarrow \text{ConcInt}(g(1-\tilde{U}_n^o),g(1-\tilde{L}^o_n))\\
&\ \ \ \ (\tilde{L}_n^o,\tilde{U}_n^o)\leftarrow (1-h(u^o),1-h(l^o))\\
&\text{end while} .
\end{align*}
\medskip

\par\noindent
\textbf{Illustration of the confidence bands}
\smallskip

\par
\begin{figure}[!htb]
\center{\includegraphics[width=\textwidth]{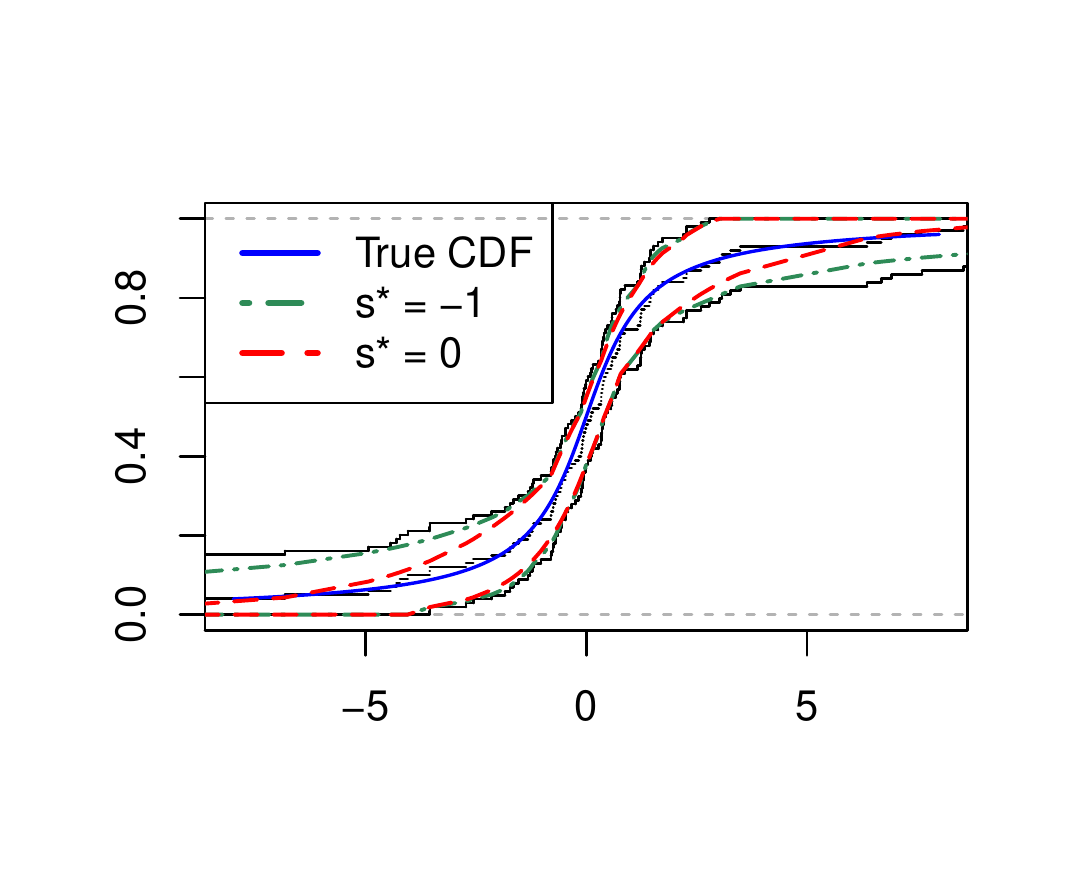}}
\caption{
Confidence bands for bi-$s^*$-concave distribution functions based on KS bands.
The black curve is the distribution function of the Student-$t$ distribution with $1$ degree of freedom.
The two gray-black lines give the KS band and lines in other colors are 
refined confidence bands under the bi-$s^*$-concave assumption.
The step function in the middle is the empirical  distribution function.
}
\label{plot:KS1c} 
\end{figure}

\begin{figure}[!htb]
\center{\includegraphics[width=\textwidth]{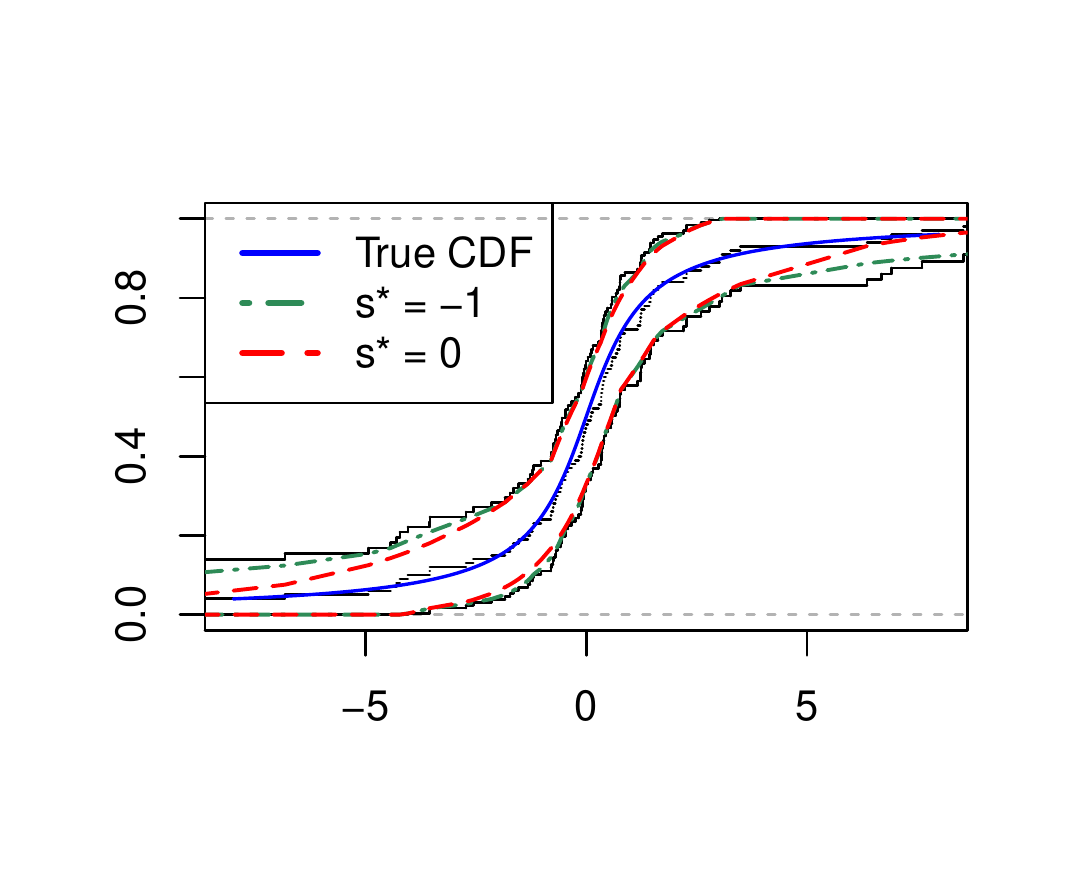}}
\caption{Confidence bands for bi-$s^*$-concave distribution functions based on WKS bands.
The black curve is the distribution function of the Student-$t$ distribution with $1$ degree of freedom.
The two gray-black lines give the WKS band and lines in other colors are 
refined confidence bands under the bi-$s^*$-concave assumption.
The step function in the middle is the empirical distribution function.
}
\label{plot:WKS2c} 
\end{figure}

\noindent
To get some feeling for the new confidence bands in a setting in which $s_0^*$ is known,
we generated a sample of size $n = 100$ from the Student-$t$ distribution with $r=1$ degrees of freedom. 
This distribution belongs to $\mathcal{P}_{s^*}$ for every $s^*\leq-1 \equiv s_0^*$. 
We constructed Kolmogorov-Smirnov (KS) and weighted Kolmogorov-Smirnov (WKS) bands 
with $\gamma=0.4$ as the initial starting bands $(L_n,U_n)$.
We 
then computed and plotted our shape constrained confidence bands 
$(L_n^0 , U_n^0)$ under the (correct) assumption that $s^*=-1$ 
and  the (incorrect) assumption that $s^* = 0$  for both the  KS and WKS bands as initial nonparametric bands with 
for $\alpha=0.05$; see Figure~\ref{plot:KS1c} and Figure~\ref{plot:WKS2c}. 
To see the components of Figures 1 and 2 separately, see 
Section 8, Appendix 2,  Figures 1-2 and 3-4 respectively.  

Note that when $s^*=0$,
$s^* $ is miss-specified and the resulting bands are not guaranteed to have 
coverage probabiltiy $.95$.  An indication of this is that the shape constrained 
bands computed under the assumption $s^* = 0$ do not contain 
the empirical distribution.  

From these two plots, an immediate observation is that the confidence bands for smaller $s^*$ 
are wider than those with larger $s^*$.  This is a direct consequence of the nested property 
of the bi-$s^*$-concave classes; see Proposition~\ref{Prop:NestedBisStar}.  
Also note that the shape constrained band with $s^* = -1$ does improve on the KS band, especially 
in the tail.  

\medskip

\begin{figure}[htb]
\center{\includegraphics[width=\textwidth]{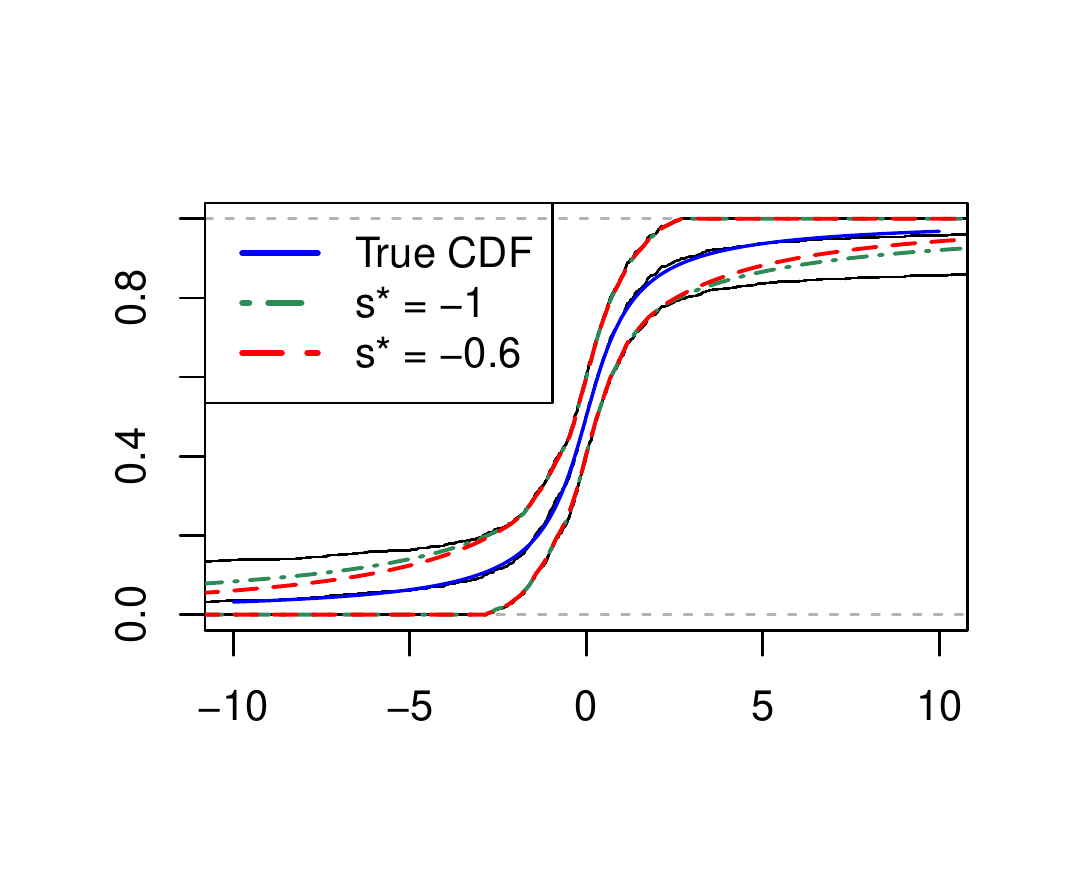}}
\caption{
Confidence Bands for bi-$s^*$-concave distribution functions from KS bands based  
on a sample of size $1000$ from the Student-$t$ distribution with $1$ degree of freedom.
The two gray-black lines give the initial bands,  lines in other colors are 
refined confidence bands under the bi-$s^*$-concave assumption. 
The step function (black) in the middle is the empirical distribution function.
}
\label{plot:KS3c}
\end{figure}

\begin{figure}
\includegraphics[width=\textwidth]{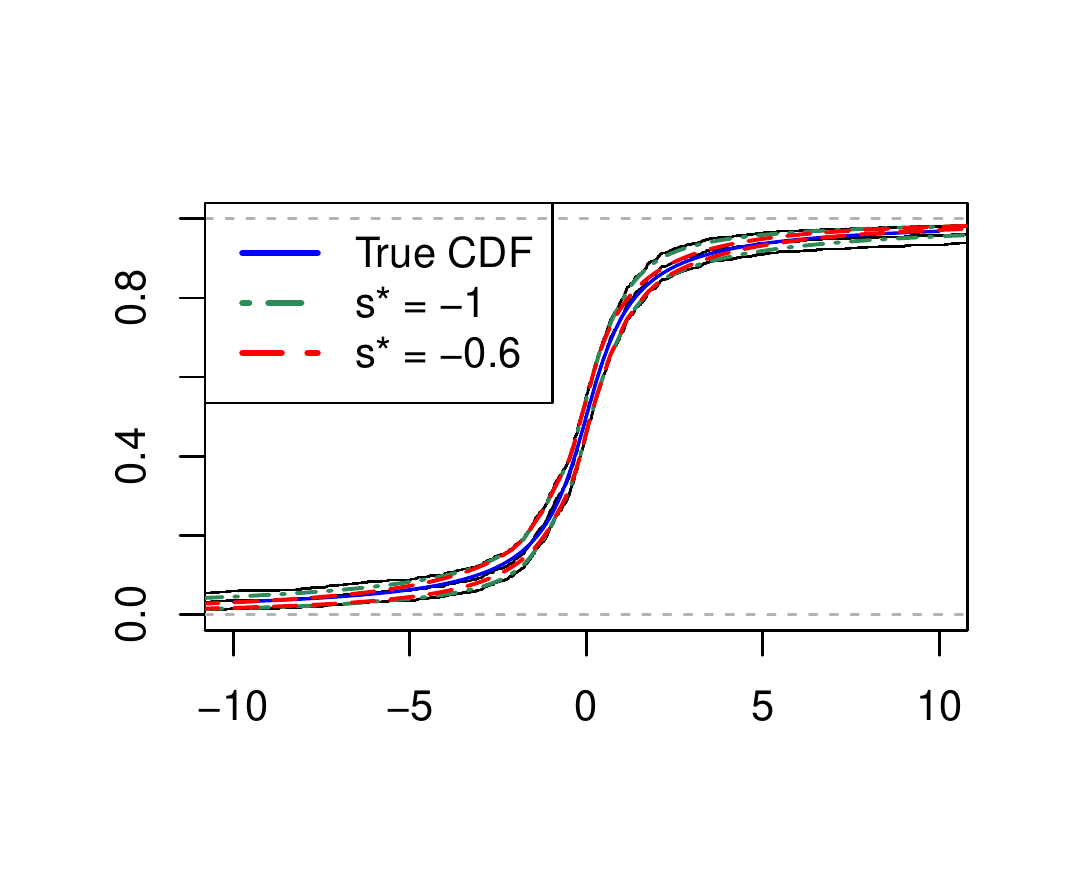}
\caption{
Confidence Bands for bi-$s^*$-concave distribution functions from  WKS bands based 
on a sample of size $1000$ from the Student-$t$ distribution with one degree of freedom.
The two gray-black lines give the initial bands,  lines in other colors are 
refined confidence bands under the bi-$s^*$-concave assumption. 
The step function (black) in the middle is the empirical distribution function.
}
\label{plot:WKS4c} 
\end{figure}

\par\noindent
\noindent\textbf{An Application}
\smallskip

\par\noindent
\cite{DUMBGEN20171} gave an application of bi-log-concave confidence bands to a dataset from \cite{Woolridge2000}. 
It contains approximate annual salaries of the CEOs of 177 randomly chosen companies in the U.S. 
The salary is rounded to multiples of 1000 USD. 
We denote the $i$-th observed approximate salary by $Y_{i,raw}$.
\cite{DUMBGEN20171} assume that the unobserved true salary $Y_{i,true}$ lies within $(Y_{i,raw} - 1, Y_{i,raw} + 1)$. 
Let us assume that $G_{true}$ is the unknown distribution of $Y_{true}$.  
For income data it is sometimes assumed that $\log_{10}Y_{true}$ is Gaussian (see \cite{Kleiber2003}). 
Since Gaussian densities are all log-concave and hence have bi-log-concave distribution functions 
(by Proposition~\ref{Prop:s-concaveAndBisStar}),
it is natural to consider replacing the Gaussian assumption by the assumption of bi-log-concavity.
\cite{DUMBGEN20171} therefore assumed that $X = \log_{10}Y_{true}$ is bi-log-concave and constructed 95\% 
confidence bands $(L_n,U_n)$ (see Figure 4 of \cite{DUMBGEN20171}) where $L_n$ is computed 
with the empirical distribution of $\log_{10}(Y_{i,raw} - 1)_{i=1}^n$ and $U_n$ is 
computed with that of $\log_{10}(Y_{i,raw} +1)_{i=1}^n$. 

Here we assume that the distribution of $X$ is bi-$s^*$-concave for some $s^*$ and compute confidence bands for different values of $s^*$.  
Now we are confronted with the issue of choosing $s^*$:  if we want narrower confidence 
bands we would assume some value of $s^* \in (0,1]$, 
while if we are not willing to assume $s^*=0$ (the choice made by \cite{DUMBGEN20171}, then we would assume some value of $s^* < 0$ (leading to 
the larger classes ${\cal P}_{s^*}$ with $s^*<0$. 
It is of some interest to know if the CEO data could be modeled by use of the bi-$s^*$ classes with $s^* \in (0,1]$ since this would 
result in still narrower confidence bands.  
But it is also of interest to try to use the data to choose $s^*$.  
\smallskip

\par\noindent
\textbf{Choosing $s^*$} 

 Since $F$ can be a member of $ \mathcal{P}_{s^*}$ for various values of 
 $s^*$, each $s^*$  leads to a  different set of bands.
  However, due to the nesting property of $ \mathcal{P}_{s^*}$, a larger $s^*$ always 
  yields a narrower  confidence band. Thus, it is of  interest to estimate
  \[
  s^*_0 (F) :=\sup\{s^*\in (-\infty,1] : F\in \mathcal {P}_{s^*}\} 
  \]
 since $s^*=s^*_0$  generates the narrowest  bands at a given confidence level.   
  If $F$ is not bi-$s^*$-concave for any $s^*\leq 1$, then we set $s^*_0(F) =-\infty$. 
 Now  $s^*_0$ is connected to  the Cs\"org{\H o} - R\'ev\'esz constant since 
  $s^*=s^*_0$ when $\overline{\gamma}(F)=1-s^*$ and $F\in\mathcal{P}^*_s$. 
  For example,  
  the Student-$t$ distribution with $r$ ``degree of freedom''  has $s^*_0=-1/r$. 
  However, this connection cannot 
  be easily exploited for practical estimation purposes due to difficulties in estimating $\gamma (F)$ or $\overline{\gamma}(F)$.
  So we take an alternative route to making inference about $s^*_0$.
  
Starting from  an initial $1-\alpha$  band $(L_n,U_n)$, a bound on  $s^*_0$ is given by 
$$
\overline{s}^*_n = \sup\{ s^* \in (-\infty, 1 ]: \ (L_n, U_n) \ \mbox{contains some d.f.}  \ F \in {\cal P}_{s^*} \}.
$$
Clearly, for $s^* >\overline{s}^*_n$,
there is no bi-$s^*$-concave distribution function fitting into the band $(L_n, U_n)$.  
Since this happens with probability at most $\alpha \in (0,1)$ when the true distribution 
function $F \in {\cal P}_{s^*}$, it follows that $(-\infty, \overline{s}^*_n]$ is   
a confidence set for $s^*_0$ with coverage probability at least $1-\alpha$.    
Our simulations suggest that  $\overline{s}^*_n$ is generally considerably larger than $s^*_0$, 
and hence not suitable as an estimator, especially for $\alpha=0.05$.   

Instead,  we propose an estimator of $s^*_0$ based on the $\FF_n$ measure of the set where 
the empirical measure remains between the shape-constrained band for $s^*$.  More formally,
let $L_n^o (s^*)$ and $U_n^o(s^*)$ denote the $1-\alpha$ level bi-$s^*$-concave confidence 
 bands based on the initial bands $L_n$ and $U_n$ and the assumption $F \in {\cal P}_{s^*}$.
 Define
 \begin{eqnarray*}
 \omega(s^*)
 & := & n^{-1} \sum_{i=1}^n 1\{ L_n^o (s^*)(X_i) \le \FF_n(X_i)  \le  U_n^o (s^*)(X_i) \} \\
 && \qquad \qquad \cdot 1\{ L_n^o (s^*)(X_i)  \le  U_n^o (s^*)(X_i) \} \\
 & = & \FF_n  \left ( \{ L_n^o (s^*) \le \FF_n  \le  U_n^o (s^*) \} \cap \{ L_n^o (s^*)  \le  U_n^o (s^*) \} \right ) .
 \end{eqnarray*}
A higher  value of $\omega(s^*)$ indicates that $(L_n^o (s^*),U_n^o(s^*))$ contains a greater portion of $\FF_n$. 
 Since the bands $(L_n(s^*),U_n(s^*))$ become narrower as $s^*$ increases, $\omega(s^*)$ 
 decreases in $s^*$, and eventually becomes zero when $s^*>\overline{s}^*_n$.  
A plausible estimator of $s^*_0$ is therefore given by
 \begin{equation}\label{def: widehat s n}
 \widehat{s}_n^*=\min \{ s^*\in(-\infty, \overline{s}^*_n]:\omega(s^*)>\rho\},
\end{equation}
where $\rho$ is a threshold  taking values in $(0,1)$. The calculation of $\widehat{s}_n^*$ thus 
depends on $\alpha$ and $\rho$. 

\begin{figure}[!htb]
\center{\includegraphics[width=\textwidth]{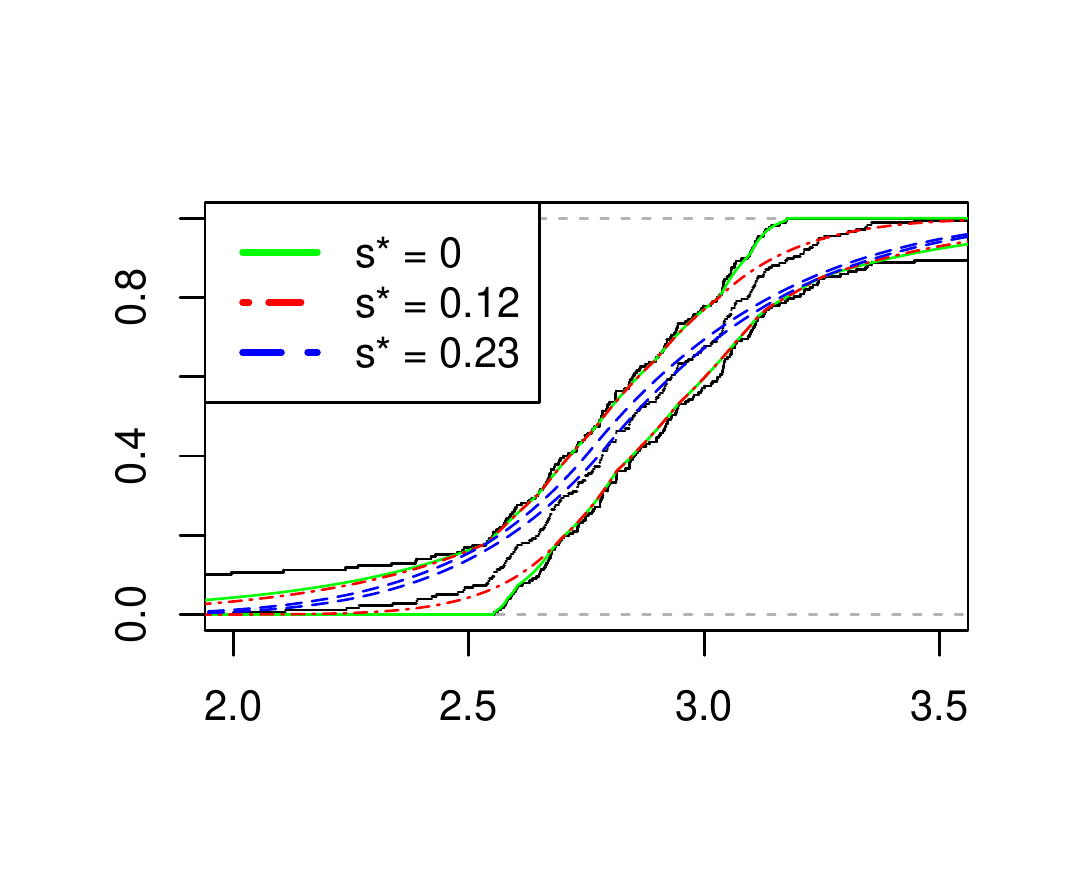}}
\caption{
Confidence Bands from an initial KS band for the CEO salary data.
The step function in the middle is the empirical distribution function.
The two gray-black lines give the KS band and lines in other colors are refined confidence bands under the bi-$s^*$-concave assumption.
}
\label{plot:KS5c} 
\end{figure}

\begin{figure}[!htb]
\center{\includegraphics[width=\textwidth]{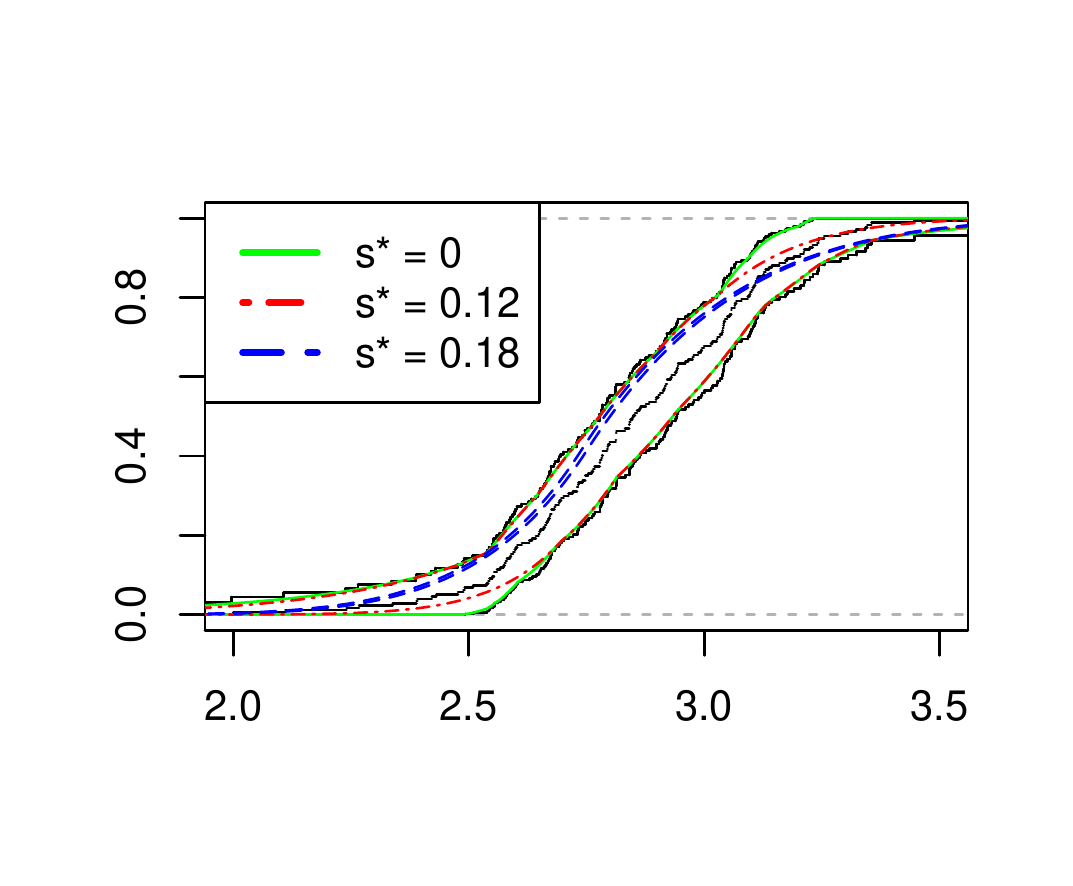}}
\caption{
Confidence Bands from an initial WKS band for the CEO salary data.
The step function in the middle is the empirical distribution function.
The two gray-black lines give the WKS band and lines in other colors are refined 
confidence bands under the bi-$s^*$-concave assumption.
}
\label{plot:WKS6c} 
\end{figure}

In the case of the CEO data,  $\overline{s}_n^* \approx 0.23$ for the KS initial band, 
and  $\overline{s}_n^* \approx 0.18$ for the WKS band.
Taking $\alpha=0.05$ and  $\rho = .95$,  leads to $\widehat{s}_n^* = 0.12$ , 
 while 
 taking    $\alpha=0.05$ and  $\rho = 0.95$,  leads to $\widehat{s}_n^* = .12$.  
The resulting bands are given in Figures~\ref{plot:KS5c} and~\ref{plot:WKS6c}.
Also see Section 8, Appendix 2, Figures 9 - 10 
and Figures 11-12 
for the steps in constructing 
 Figures~\ref{plot:KS5c} and~\ref{plot:WKS6c}.

We should emphasize that our current theory says little about the coverage probabilities 
of the bands $(L_n^o (s^*), U_n^o (s^*) )$.   
Discussion of the consistency of $\widehat{s}_n^*$ is beyond the scope 
of the present paper, but this and further issues concerning inference for both $s^*$ and $F \in {\cal P}_s$ 
seem to be interesting directions for future research.

\section{Summary and further problems} 
\label{Section:Question}
\par\noindent
In this paper we have:\\ 
$\bullet$  Defined new classes of shape-constrained distribution functions, the bi-$s^*$-classes extending the 
bi-log-concave class of distribution functions defined by  \cite{DUMBGEN20171}. \\
$\bullet$  Characterized the new classes and connected our characterization to an important parameter, the 
Cs\"org{\H o} - R\'ev\'esz constant associated with a distribution function $F$.\\
$\bullet$ Used the new bi-$s^*$-concave classes to define refined confidence bands for distribution functions which exploit 
the shape constraint, thereby producing more accurate (narrower) bands with honest coverage when the shape constraint holds.  

Thus we have shown that if we know the parameter $s^* \in (-\infty, 1]$ determining the class, we can construct refined
confidence bands which improve on any given nonparametric confidence bands if the given value of $s^*$ is correct.
It follows from the construction of our bands 
that they have conservative coverage probabilities under the (null) hypothesis that the  true distribution function 
is in ${\cal P}_{s^*}$ and that $s^*$ is correctly specified. \\
$\bullet$
What if we do not know $s^*$?  Can we estimate it from the data?   
As becomes clear from the discussion 
of the CEO data via Figures 5 and 6, our methods provide one-sided confidence bounds for the true $s^{*}$ of the form 
$(-\infty, \overline{s}_n^*]$ under the assumption that $F \in {\cal P}_{s^*}$ for some $s^*$.   
It remains to develop inference methods for $s^*$ and $(s^*, F)$ jointly.  It will also be of interest to 
have a more complete understanding of the power behavior of tests related to  $\overline{s}_n^*$ and 
$\widehat{s}_n^*$. \\ 
$\bullet$ The stable laws are known to be unimodal;
see e.g. \cite{MR771431} 
for some history.
In connection with Example~\ref{Ex:LevyOneHalf} we have the following:  \\
{\bf Conjecture:}  the $\alpha-$stable laws are 
$s-$concave with $s = - 1/(1+\alpha)$ for $0 < \alpha < 2$.  

\section{Proofs}
\label{Section:Proofs}
\textbf{Proof of Theorem \ref{Thm:CharacterizingThm}:}\\
Throughout our proof we will denote $\inf J(F)$ and $\sup J(F)$ by $a$ and $b$ respectively. 
Moreover, we assume $s^*<0$ in the following proof and leave the case of $s^*>0$ for the Appendix. 
Note that the case $s^*=0$ is proved by \cite{DUMBGEN20171}.\\
(i) implies (ii): \\
Suppose $F\in\mathcal{P}_{s^*}$. To prove that $F$ is continuous on $\RR$, 
we first note that $x\mapsto F^{s^*}(x)$ and $x\mapsto \left(1-F(x)\right)^{s^*}(x)$ are convex functions on $\RR$. 
By Theorem 10.1 (page 82) of \cite{10.2307/j.ctt14bs1ff}, $F^{s^*}$ and 
$\left(1-F(x)\right)^{s^*}$ are continuous on any open convex sets in their effective domains.
In particular, $F^{s^*}$ and $\left(1-F\right)^{s^*}$ are continuous on $(a,\infty)$ and $(-\infty,b)$ respectively. 
This implies that $F$ is continuous on $(a,\infty)$ and $(-\infty,b)$, or equivalently, 
on $(a,\infty)\cup(-\infty,b)=(-\infty,\infty)$ since $F$ is non-degenerate.\\
To prove that $F$ is differentiable on $J(F)$, note that $J(F)=(a,b)$ since $F$ is continuous on $\RR$.
By Theorem 23.1 (page 213) of \cite{10.2307/j.ctt14bs1ff}, for any $x\in J(F)$, the convexity of 
$F^{s^*}$ on $J(F)$ implies the existence of $\left(F^{s^*}\right)^{\prime}_+(x)$ and $\left(F^{s^*}\right)^{\prime}_-(x)$.
Moreover, $\left(F^{s^*}\right)^{\prime}_-(x)\leq \left(F^{s^*}\right)^{\prime}_+(x)$ by 
Theorem 24.1 (page 227) in \cite{10.2307/j.ctt14bs1ff}. Since $F=\left(F^{s^*}\right)^{1/s^*}$ 
on $J(F)$, the chain rule guarantees the existence of $F^\prime_\pm(x)$ and 
\[
F^\prime_\pm(x)=\frac{1}{s^*}\left(F^{s^*}\right)^{1/s^*-1}(x\pm)\left(F^{s^*}\right)^{\prime}_\pm(x).
\]
Since $F$ is continuous on $J(F)$, then
\[
F_\pm^\prime(x)=\frac{1}{s^*}\left(F^{s^*}\right)^{1/s^*-1}(x)\left(F^{s^*}\right)^{\prime}_\pm(x).
\]
Hence $F^\prime_-(x)\geq F^\prime_+(x)$ by noting that 
$\left(F^{s^*}\right)^{\prime}_-(x)\leq \left(F^{s^*}\right)^{\prime}_+(x)$ and $s^*<0$. \\
Similarly, one can prove $F^\prime_-(x)\leq F^\prime_+(x)$ by the convexity of $\left(1-F\right)^{s^*}$ on $J(F)$.\\
Thus $F^\prime_-(x)=F^\prime_+(x)=F^\prime(x)$ for any $x\in J(F)$, or equivalently, $F$ is differentiable on $J(F)$. 
The derivative of $F$ is denoted by $f$, i.e. $f\equiv F^\prime$.\\
To prove (\ref{Formula:TailBounds1}), note that the convexity of $x\mapsto F^{s^*}(x)$ on $J(F)$ implies that, 
for any $x,y\in J(F)$,
\[
F^{s^*}(y)-F^{s^*}(x)\geq (y-x)\left(F^{s^*}\right)^{\prime}(x)=(y-x)s^*F^{s^*-1}(x)f(x),
\]
or, with $x_+=\max\{x,0\}$,
\[
\frac{F^{s^*}(y)}{F^{s^*}(x)}\geq \left(1+s^*\frac{f(x)}{F(x)}(y-x)\right)_+.
\]
Hence,
\[
\frac{F(y)}{F(x)}\leq \left(1+s^*\frac{f(x)}{F(x)}(y-x)\right)_+^{1/s^*},
\]
or, equivalently,
\[
F(y)\leq F(x)\left(1+s^*\frac{f(x)}{F(x)}(y-x)\right)_+^{1/s^*}.
\]
Analogously, the convexity of $\left(1-F(x)\right)^{s^*}$ on $J(F)$ implies that 
\[
\left(1-F(y)\right)^{s^*}-\left(1-F(x)\right)^{s^*}\geq -(y-x)s^*\left(1-F(x)\right)^{s^*-1}f(x),
\]
or, equivalently,
\[
\left(\frac{1-F(y)}{1-F(x)}\right)^{s^*}\geq \left(1-s^*\frac{f(x)}{1-F(x)}(y-x)\right)_+,
\]
which yields
\[
F(y)\geq 1- \left(1-F(x)\right)\left(1-s^*\frac{f(x)}{1-F(x)}(y-x)\right)^{1/s^*}_+.
\]
The proof of (\ref{Formula:TailBounds1}) is complete.\\
(ii) implies (iii):\\
Applying (\ref{Formula:TailBounds1}) yields that for any $x,y\in J(F)$ with $x<y$, 
\[
\frac{F^{s^*}(x)}{F^{s^*}(y)}\geq 1+s^*\frac{f(y)}{F(y)}(x-y),
\]
and
\[
\frac{F^{s^*}(y)}{F^{s^*}(x)}\geq 1+s^*\frac{f(x)}{F(x)}(y-x),
\]
or, equivalently, 
\[
F^{s^*}(x)\geq F^{s^*}(y)+s^*\frac{f(y)}{F^{1-s^*}(y)}(x-y),
\]
and
\[
F^{s^*}(y)\geq F^{s^*}(x)+s^*\frac{f(x)}{F^{1-s^*}(x)}(y-x).
\]
By defining $h\equiv f/F^{1-s^*}$ on $J(F)$, it follows that 
\[
F^{s^*}(x)\geq F^{s^*}(y)+s^*h(y)(x-y),
\]
and
\[
F^{s^*}(y)\geq F^{s^*}(x)+s^*h(x)(y-x).
\]
After summing up the last two inequalities, it follows that 
\[
F^{s^*}(x)+F^{s^*}(y)\geq F^{s^*}(y)+s^*h(y)(x-y)+F^{s^*}(x)+s^*h(x)(y-x),
\]
or, equivalently,
\[
0\geq s^*\left(h(x)-h(y)\right)(y-x).
\]
Hence $h(x)\geq h(y)$, or equivalently, $h(\cdot)$ is a monotonically non-increasing function on $J(F)$.\\
The proof of the monotonicity of
$\tilde{h}\equiv f/(1-F)^{1-s^*}$ is similar and hence is omitted.\\
(iii) implies (iv):\\
If (iii) holds, it immediately follows that $f>0$ on $J(F)=(a,b)$. 
If not, suppose that $f(x_0)=0$ for some $x_0\in J(F)$.
It follows that $h(x_0)=f(x_0)/F^{1-s^*}(x_0)=0$.
Since $h$ is monotonically non-increasing on $J(F)$,
$h(x)=0$ for all $x\in [x_0,b)$,
or, equivalently,
$f=0$ on $[x_0,b)$.
Similarly, the non-decreasing monotonicity of $x\mapsto\tilde{h}(x)$ on $J(F)$ implies that $f=0$ on $(a,x_0]$.
Then $f=0$ on $J(F)$, which violates the continuity assumption in (iii) and hence $f>0$ on $J(F)$.\\
To prove $f$ is bounded on $J(F)$, note that the monotonicities of $h$ and $\tilde{h}$ imply that for any $x,x_0\in J(F)$,
\begin{eqnarray*}
f(x)=\left\{
\begin{tabular}{ll }
$F^{1-s^*}(x)h(x)\leq h(x)\leq h(x_0)$, &if $x\geq x_0$,\\
$(1-F(x))^{1-s^*}\tilde{h}(x)\leq \tilde{h}(x)\leq \tilde{h}(x_0)$, & if $x\leq x_0$.
\end{tabular}\right.
\end{eqnarray*}
Hence $f(x)\leq \max\{h(x_0),\tilde{h}(x_0)\}$ for any $x,x_0\in J(F)$.\\
To prove that $f$ is differentiable on $J(F)$ almost everywhere, we first 
prove that $f$ is Lipschitz continuous on $(c,d)$ for any $c,d\in J(F)$ with $c<d$.\\
By the non-increasing monotonicity of $h$ on $J(F)$, the following arguments 
yield an upper bound of $\left(f(y)-f(x)\right)/(y-x)$ for any $x,y\in(c,d)$:
\begin{eqnarray*}
\frac{f(y)-f(x)}{y-x}
&=&\frac{F^{1-s^*}(y)h(y)-F^{1-s^*}(x)h(x)}{y-x}\\
&=&h(y)\frac{F^{1-s^*}(y)-F^{1-s^*}(x)}{y-x}+F^{1-s^*}(x)\frac{h(y)-h(x)}{y-x}\\
&\leq&h(y)\frac{F^{1-s^*}(y)-F^{1-s^*}(x)}{y-x}\\
&=&h(y)(1-s^*)f(z)F^{-s^*}(z),
\end{eqnarray*}
where the last equality follows from the mean value theorem and $z$ is between $x$ and $y$.\\
Since $-s^*>0$, it follows that $F^{-s^*}<1$ and hence 
\[
\frac{f(y)-f(x)}{y-x}<(1-s^*)f(z)h(y)\leq (1-s^*)\max\{h(x_0),\tilde{h}(x_0)\}h(c)
\]
for $x,y\in (c,d)$.\\
Similar arguments imply that 
\begin{eqnarray*}
\frac{f(y)-f(x)}{y-x}
&=&\frac{\bar{F}^{1-s^*}(y)\tilde{h}(y)-\bar{F}^{1-s^*}(x)\tilde{h}(x)}{y-x}\\
&=&\tilde{h}(y)\frac{\bar{F}^{1-s^*}(y)-\bar{F}^{1-s^*}(x)}{y-x}+\bar{F}^{1-s^*}(x)\frac{\tilde{h}(y)-\tilde{h}(x)}{y-x}\\
&\geq&\tilde{h}(y)\frac{\bar{F}^{1-s^*}(y)-\bar{F}^{1-s^*}(x)}{y-x}\\
&=&-\tilde{h}(y)(1-s^*)\bar{F}^{-s^*}(z)f(z)\\
&\geq&-(1-s^*)\max\{h(x_0),\tilde{h}(x_0)\}\tilde{h}(d).
\end{eqnarray*}
Hence 
\[
\left|\frac{f(y)-f(x)}{y-x}\right|\leq (1-s^*)\max\{h(x_0),\tilde{h}(x_0)\}\max\{h(c),\tilde{h}(d)\}.
\]
The last display shows that $f$ is Lipschitz continuous on $(c,d)$.\\
By Proposition 4.1(iii) of \cite{Shorack2017}, page 82, $f$ is absolutely continuous on $(c,d)$, 
and hence $f$ is differentiable on $(c,d)$ almost everywhere.\\
Since $(c,d)$ is an arbitrary interval in $(a,b)$, the differentiability of $f$ on $(c,d)$ 
implies the differentiability of $f$ on $(a,b)$ and hence $f$ is differentiable on 
$(a,b)$ with $f^\prime=F^{\prime\prime}$ almost everywhere.\\
Since $f$ is differentiable almost everywhere, the non-increasing monotonicity of $h$ on $J(F)$ implies that
\[
h^\prime(x)\leq 0 \text{ almost everywhere on $J(F)$,}
\]
or, equivalently,
\[
\log(h)^\prime(x)\leq 0 \text{ almost everywhere on $J(F)$.}
\]
Straight-forward calculation yields that the last display is equivalent to 
\[
\frac{f^\prime}{f}-(1-s^*)\frac{f}{F}\leq 0 \text{ almost everywhere on $J(F)$,}
\]
or,
\[
f^{\prime} \leq(1-s^*)\frac{f^2}{F}\text{ almost everywhere on $J(F)$,}
\]
which is the right hand side of (\ref{Formula:fPrimeBounds}).\\
Similarly, the non-decreasing monotonicity of $\tilde{h}$ implies the left hand side of (\ref{Formula:fPrimeBounds}).\\
(iv) implies (i):\\
Since $F$ is continuous on $\RR$, it suffices to prove that $F^{s^*}$ is convex on $J(F)$ by Definition \ref{Definition2}. 
Since we assume that $F$ is differentiable on $J(F)$ with derivative $f=F^\prime$, 
the convexity of $F^{s^*}$ on $J(F)$ can be proved by the increasing monotonicity of the first derivative of $F^{s^*}$ on $J(F)$. 
Since $f$ is differentiable almost everywhere on $J(F)$, the increasing monotonicity 
of $\left(F^{s^*}\right)^{\prime}$ on $J(F)$ can be proved by the non-negativity of 
$\left(F^{s^*}\right)^{\prime\prime}$ on $J(F)$ almost everywhere, which follows from
\begin{eqnarray*}
\left(F^{s^*}\right)^{\prime\prime}(x)=s^*F^{s^*-1}(x)\left(-(1-s^*)\frac{f^2(x)}{F(x)}+f^\prime(x)\right)\geq 0,
\end{eqnarray*}
where $f=F^\prime, f^\prime=F^{\prime\prime}$. 
The last inequality follows from the right hand side of (\ref{Formula:fPrimeBounds}).\\
Similarly, the convexity of $\left(1-F(x)\right)^{s^*}$, or $\bar{F}^{s^*}$, on $J(F)$ can be proved by the following arguments:
\[
\left(\bar{F}^{s^*}\right)^{\prime\prime}(x)=s^*\bar{F}^{s^*-1}(x)\left(-(1-s^*)\frac{f^2(x)}{\bar{F}(x)}-f^\prime(x)\right)\geq 0,
\]
where the last inequality follows from the left part of (\ref{Formula:fPrimeBounds}).
\hfill $\Box$
\medskip

\par\noindent
\textbf{Proof of Proposition \ref{Prop:s-concaveAndBisStar}:}\\
First some background and definitions:\\
\begin{itemize}
  \item Let $a,b\geq 0$ and $\theta\in(0,1)$. The generalized mean of order $s\in\RR$ is defined by 
\begin{eqnarray*}
M_s(a,b;\theta)=\left\{
\begin{tabular}{ll }
$\left((1-\theta)a^s+\theta b^s\right)^{1/s}$, &if $\pm s\in(0,\infty)$,\\
$a^{1-\theta}b^\theta$, & if $s=0$,\\
$\max\{a,b\}$, &  if $s=\infty$,\\
$\min\{a,b\}$, &   if $s=-\infty$,
\end{tabular}\right.
\end{eqnarray*}
  \item Let $(M,d)$ be a metric space with Borel $\sigma-$field $\mathcal{M}$. A measure $\mu$ 
  on $\mathcal{M}$ is called $t$-concave if for nonempty sets $A,B\in\mathcal{M}$ and $0<\theta<1$ we have
\[
\mu_*\left((1-\theta)A+\theta B\right)\geq M_t\left(\mu_*(A),\mu_*(B);\theta\right)
\]
where $\mu_*$ is the inner measure corresponding to $\mu$   
(which is needed in general in view of examples noted by \cite{9e37eb77de0544d1837f114247289da9}).	
\item A non-negative real-valued function $h$ on $(M,d)$ is called $s$-concave if for $x,y\in M$ and $0<\theta<1$ we have 
\[
h\left((1-\theta)x+\theta y\right)\geq M_s\left(h(x),h(y);\theta\right).
\]
See Chapter 3.3 in \cite{trove.nla.gov.au/work/12919064} for more details of the definitions of $M_s(a,b;\theta)$, $t$-concave and $s$-concave.
 \item Suppose $(M,d)=(\mathbb{R}^k,|\cdot|)$, $k-$dimensional Euclidean space with the 
 usual Euclidean metric and suppose that $f$ is an $s$-concave density function with respect to 
 Lebesgue measure $\lambda$ on $\mathcal{B}_k$, and consider the probability measure $\mu$ on $\mathcal{B}_k$ defined by 
\[
\mu(B)=\int_B f d\lambda \text{ for all $B\in\mathcal{B}_k$}.
\]
Then by a theorem of \cite{Borell1975}, \cite{BRASCAMP1976366} and \cite{rinott1976}, the measure 
$\mu$ is $s^*$-concave where $s^*=1/(1+ks)$ if $s\in(-1/k,\infty)$ and $s^*=0$ if $s=0$.
\item Here we are in the case $k=1$. Thus for $s\in(-1,\infty)$ the measure $\mu$ is $s^*$ concave: 
for $s\in(-1,\infty)$, $A,B\in\mathcal{B}_1,$ and $0<\theta<1$,
\begin{eqnarray}
\mu_*\left((1-\theta)A+\theta B\right)\geq M_{s^*}\left(\mu_*(A),\mu_*(B);\theta\right);
\label{Proof:Prop2}
\end{eqnarray}
here $\mu_*$ denotes inner measure corresponding to $\mu$. 
\end{itemize}
With this preparation we can give our proof of Proposition \ref{Prop:s-concaveAndBisStar}: if 
$A=(-\infty,x]$ and $B=(-\infty,y]$ for $x,y\in J(F)$, it is easily seen that 
\begin{eqnarray*}
 (1-\theta)A+\theta B 
&=&\{(1-\theta)x^\prime +\theta y^\prime: x^\prime\leq x, y^\prime\leq y\}\\
&\subset& \{(1-\theta)x^\prime +\theta y^\prime: (1-\theta)x^\prime+\theta y^\prime \leq (1-\theta)x+\theta y\}\\
&=& (-\infty,(1-\theta)x+\theta y].
\end{eqnarray*}
Therefore, with the second inequality following from (\ref{Proof:Prop2}), 
\begin{eqnarray*}
F\left((1-\theta)x+\theta y\right)
&=&\mu\left((-\infty,(1-\theta)x+\theta y]\right)\\
&\geq&\mu\left((1-\theta)(-\infty,x]+\theta(-\infty, y]\right)\\
&\geq&M_{s^*}\left(\mu((-\infty,x]),\mu((-\infty, y]);\theta\right)=M_{s^*}(F(x),F(y);\theta);
\end{eqnarray*}
i.e. $F$ is $s^*$-concave. Similarly, taking $A=(x,\infty)$ and $B=(y,\infty)$ it follows that $1-F$ is $s^*$-concave.\\
Note that this argument contains the case $s^*=0$.
\hfill $\Box$
\medskip

\par\noindent
\textbf{Proof of Proposition \ref{Prop:NestedBisStar}:}\\
By Theorem \ref{Thm:CharacterizingThm}, for any $F\in\mathcal{P}_{s^*}$, $F$ is continuous 
on $\RR$ and differentiable on $J(F)$ with derivative $f=F^\prime$. Furthermore, $f$ is 
differentiable almost everywhere on $J(F)$ with derivative $f^\prime=F^{\prime\prime}$ 
satisfying (\ref{Formula:fPrimeBounds}).\\
For any $t^*\leq s^*$,  by noting that $1-s^*\leq 1-t^*$ and $-(1-s^*)\geq -(1-t^*)$, 
it follows that 
\[
-(1-t^*)\frac{f^2}{1-F}\leq -(1-s^*)\frac{f^2}{1-F}\leq f^\prime \leq (1-s^*) \frac{f^2}{F} \leq (1-t^*) \frac{f^2}{F},
\]
almost everywhere on $J(F)$. Hence $F\in\mathcal{P}_t^*$ by Theorem \ref{Thm:CharacterizingThm}. 
This proves (\ref{NestedProperty}).\\
To prove (\ref{ContinuityAtZero}), 
note that for any $F\in\cup_{s^*>0}\mathcal{P}_{s^*}$,
$F$ is continuous on $\RR$ and differentiable on $J(F)$ with derivative $f=F^\prime$. 
Furthermore, $f$ is differentiable almost everywhere on $J(F)$ with derivative 
$f^\prime=F^{\prime\prime}$ satisfying (\ref{Formula:fPrimeBounds}), i.e.
\[
-(1-s^*)\frac{f^2}{1-F}\leq f^\prime \leq (1-s^*) \frac{f^2}{F} \text{ almost everywhere on $J(F)$},
\]
for all $s^*>0$.
By taking $s^*\rightarrow 0$, it follows that
\[
-(1-0)\frac{f^2}{1-F}\leq f^\prime \leq (1-0) \frac{f^2}{F} \text{ almost everywhere on $J(F)$}.
\]
The last display is equivalent to $F\in\mathcal{P}_0$ by Theorem \ref{Thm:CharacterizingThm}. 
This proves that the left hand side of (\ref{ContinuityAtZero}) holds.  Similarly, one can prove the 
right hand side of (\ref{ContinuityAtZero}); the details are omitted.
\hfill $\Box$
\medskip

\par\noindent
\textbf{Proof of Corollary \ref{Corollary:CRCondition}:}\\
To prove the right part of (\ref{Formula:InequalityOfCRCondition}), note that (\ref{Formula:fPrimeBounds}) implies that 
\[
1-s^*\geq \frac{Ff^\prime}{f^2}\text{ and }1-s^*\geq -\frac{(1-F)f^\prime}{f^2}
\]
almost everywhere on $J(F)$, or equivalently,
\[
1-s^*\geq \max \left \{\esssup_{x\in J(F)}\frac{Ff^\prime}{f^2},\esssup_{x\in J(F)}-\frac{(1-F)f^\prime}{f^2}\right \}.
\]
Replacing $\esssup_{x\in J(F)}Ff^\prime/f^2$ and 
$\esssup_{x\in J(F)}-(1-F)f^\prime/f^2$ by $\tilde{CR}(F)$ and $\tilde{CR}(\overline{F})$, 
it follows that
\[
1-s^*\geq\max\{\tilde{CR}(F),\tilde{CR}(\overline{F})\}=\overline{\gamma}(F).
\]
One can prove the left two inequalities of (\ref{Formula:InequalityOfCRCondition}) by the following arguments:
\begin{align*}
&\overline{\gamma}(F) =
\max\{\tilde{CR}(F),\tilde{CR}(\overline{F})\}\\
&=
\max\left\{\esssup_{x\in J(F)}\frac{F(x) f^{\prime}(x)}{f(x)^2}, \esssup_{x\in J(F)}-\frac{(1-F(x)) f^\prime(x)}{f(x)^2}\right\}\\
&=
\max\left\{\esssup_{x\in J(F)}\frac{F(x) f^{\prime}(x)}{f(x)^2}1_{[f^{\prime}(x)\geq 0]}, 
                \esssup_{x\in J(F)}-\frac{(1-F(x)) f^\prime(x)}{f(x)^2}1_{[f^{\prime}(x)\leq 0]}\right\}\\
&=
\max\left\{\esssup_{x\in J(F)}\frac{F(x) |f^{\prime}(x)|}{f(x)^2}1_{[f^{\prime}(x)\geq 0]}, 
\esssup_{x\in J(F)}\frac{(1-F(x)) |f^\prime(x)|}{f(x)^2}1_{[f^{\prime}(x)\leq 0]}\right\}\\
\begin{split}
&\geq
\max\left\{\esssup_{x\in J(F)}\frac{F(x)\wedge (1-F(x)) |f^{\prime}(x)|}{f(x)^2}1_{[f^{\prime}(x)\geq 0]}, \right.\\
& \left.\ \ \ \ \ \ \ \ \ \ \ \ \esssup_{x\in J(F)}\frac{F(x)\wedge (1-F(x))|f^\prime(x)|}{f(x)^2}1_{[f^{\prime}(x)\leq 0]}\right\}
\end{split}\\
&=
\esssup_{x\in J(F)}\frac{F(x)\wedge(1-F(x))|f^{\prime}(x)|}{f(x)^2}\notag\\
&=\gamma(F) \ge \tilde{\gamma} (F) 
\end{align*}
where the last inequality holds since $ u \wedge (1-u) \ge u(1-u)$ for $0 \le u \le 1$.
\hfill $\Box$
\medskip

\par\noindent
\textbf{Proof of Corollary \ref{Corollary:BoundsForF}:}\\
Note that for $s^*<0$ and $y>-1$, we have $(1+y)^{s^*}\geq 1+s^* y$. 
Replacing $y$ by $-F(x)$, where $x\in J(F)$, it follows that 
\[
(1-F(x))^{s^*}\geq 1-s^*F(x),
\]
or, by rearranging,
\[
F(x)\leq \frac{1}{s^*}\left(1-(1-F(x)\right)^{s^*})=F_U(x),
\]
where $F_U$ is a convex function on $J(F)$ if $F\in\mathcal{P}_{s^*}$. 
This proves the right hand side of (\ref{Corollary4:Formula1}) for $s^*<0$.
Similarly, replacing $y$ by $-(1-F(x))$, where $x\in J(F)$, by rearranging terms,  it follows that
\[
F(x)\geq \frac{1}{s^*}\left(F^{s^*}(x)-(1-s^*)\right)=F_L(x),
\]
which proves the left hand side of (\ref{Corollary4:Formula1}) for $s^*<0$.\\
Similarly, for $1\geq s^*>0$ and $y>-1$, we have $(1+y)^{s^*}\leq 1+s^* y$. 
Replacing $y$ by $-F(x)$, where $x\in J(F)$, it follows that 
\[
(1-F(x))^{s^*}\leq 1-s^*F(x),
\]
or, by rearranging,
\[
F(x)\leq \frac{1}{s^*}\left(1-(1-F(x)\right)^{s^*})=F_U(x),
\]
where $F_U$ is a convex function on $J(F)$ if $F\in\mathcal{P}_{s^*}$. 
This proves the right hand side of (\ref{Corollary4:Formula1}) for $s^*>0$.\\
Similarly, replacing $y$ by $-(1-F(x))$, where $x\in J(F)$, by rearranging terms,  it follows that
\[
F(x)\geq \frac{1}{s^*}\left(F^{s^*}(x)-(1-s^*)\right)=F_L(x),
\]
which proves the left hand side of (\ref{Corollary4:Formula1}) for $s^*>0$.\\
\medskip

\par\noindent
\textbf{Proof of Lemma \ref{Lemma:LemmaInConfidenceBands}}:\\
If there is no $G\in\mathcal{P}_{s^*}$ fitting in between $L_n$ and $U_n$, 
it follows that $L_n^o\equiv 1$ and $U_n^o\equiv 0$ and assertions in both (i) 
and (ii) are trivial. In the following proof, we let $G\in\mathcal{P}_{s^*}$ such that $L_n\leq G\leq U_n$.\\
(i)
It suffices to prove that for any $x\in J(G)$ the density function $g=G^\prime$ satisfies
$g(x)\leq\max\{\gamma_1,\gamma_2\}$,
because this is equivalent to Lipschitz-continuity of $G$ with the latter constant, 
and this property carries over to the pointwise infimum $L_n^o$ and supremum $U_n^o$.\\
To prove $g(x)\leq \max\{\gamma_1,\gamma_2\}$, note
that $g/G^{1-s^*}$ is monotonically non-increasing on $J(G)$ (see Theorem \ref{Thm:CharacterizingThm}(iii)), 
it follows that for $x\geq b$
\begin{eqnarray*}
\frac{g(x)}{G^{1-s^*}(x)}
&\leq&\frac{g(b)}{G^{1-s^*}(b)}
=\left(\frac{1}{s^*}G^{s^*}\right)^\prime(b)\\
&\leq&\frac{\frac{1}{s^*}G^{s^*}(b)-\frac{1}{s^*}G^{s^*}(a)}{b-a}\\
&\leq&\frac{\frac{1}{s^*}\left(v^{s^*}-u^{s^*}\right)}{b-a}
=\gamma_1.
\end{eqnarray*}
The last inequality follows from noting that $x\mapsto (1/s^*)x^{s^*}$ is a monotonically 
non-decreasing function for all $s^*\neq0$, $G(b)\leq U_n(b)\leq v$ and $G(a)\geq L_n(a)\geq u$.
Hence 
\[
g(x)\leq G^{1-s^*}(x)\gamma_1\leq\gamma_1 \text{ for }x\geq b.
\]
Similarly, by noting that $g/(1-G)^{1-s^*}$ is monotonically 
non-decreasing on $J(G)$ (see Theorem \ref{Thm:CharacterizingThm}(iii)), 
it follows that for $x\leq a$
\begin{eqnarray*}
\frac{g(x)}{(1-G(x))^{1-s^*}}
&\leq&\frac{g(a)}{(1-G(a))^{1-s^*}}\\
&=&\left(\frac{-1}{s^*}(1-G)^{s^*}\right)^\prime(a)\\
&\leq&\frac{\frac{-1}{s^*}(1-G(b))^{s^*}-\frac{-1}{s^*}(1-G(a))^{s^*}}{b-a}\\
&\leq&\frac{\frac{-1}{s^*}\left((1-v)^{s^*}-(1-u)^{s^*}\right)}{b-a}
=\gamma_2.
\end{eqnarray*}
The last inequality follows from noting that $x\mapsto - (1/s^*)(1-x)^{s^*}$ is a monotonically 
non-decreasing function for all $s^*\neq0$, $G(b)\leq v$ and $G(a)\geq u$. Hence 
\[
g(x)\leq (1-G(x))^{1-s^*}\gamma_2\leq\gamma_2 \text{ for }x\leq a.
\]
For $a<x<b$, analogously, we get two following inequalities 
\begin{eqnarray*}
g(x)
&=&G^{1-s^*}(x)\frac{g(x)}{G^{1-s^*}(x)}\\
&\leq&G^{1-s^*}(x)\frac{\frac{1}{s^*}G^{s^*}(x)-\frac{1}{s^*}G^{s^*}(a)}{x-a}\\
&=&\frac{1}{s^*}\frac{1}{x-a}\left(G(x)-G^{s^*}(a)G^{1-s^*}(x)\right)
\end{eqnarray*}
and
\begin{eqnarray*}
g(x)
&=&(1-G(x))^{1-s^*}\frac{g(x)}{(1-G(x))^{1-s^*}}\\
&\leq&(1-G(x))^{1-s^*}\frac{\frac{-1}{s^*}(1-G(b))^{s^*}-\frac{-1}{s^*}(1-G(x))^{s^*}}{b-x}\\
&=&\frac{1}{s^*}\frac{1}{b-x}\left(1-G(x)-(1-G(b))^{s^*}(1-G(x))^{1-s^*}\right).
\end{eqnarray*}
The former inequality times $(x-a)$ plus the latter inequality times $(b-x)$ yields 
\begin{eqnarray*}
g(x)&\leq&
\frac{1}{s^*}\frac{1-G^{s^*}(a)G^{1-s^*}(x)-(1-G(b))^{s^*}(1-G(x))^{1-s^*}}{b-a}
=\frac{h(G(x))}{b-a},
\end{eqnarray*}
where 
\[
h(y)
\equiv
\frac{1}{s^*}\left(1-G^{s^*}(a)y^{1-s^*}-(1-G(b))^{s^*}(1-y)^{1-s^*}\right)
\text{ for } y\in(0,1).
\]
Since 
\[
h^{\prime\prime}(y)
=(1-s^*) \left(G^{s^*}(a)y^{-s^*-1}+(1-G(b))^{s^*}(1-y)^{-s^*-1}\right)
\geq0,
\]
it follows that $h(y)$ is convex on $(0,1)$ and hence 
\[
g(x)\leq \max_{y\in\{G(a),G(b)\}}\frac{h(y)}{b-a}=\max\left\{\frac{h(G(a))}{b-a},\frac{h(G(b))}{b-a}\right\}.
\]
Note that 
\[
\frac{h(G(a))}{b-a}=(1-G(a))^{1-s^*}\frac{\frac{-1}{s^*}(1-G(b))^{s^*}-\frac{-1}{s^*}(1-G(a))^{s^*}}{b-a}\leq \gamma_2
\]
and
\[
\frac{h(G(b))}{b-a}=G(b)^{1-s^*}\frac{\frac{1}{s^*}G^{s^*}(b)-\frac{1}{s^*}G^{s^*}(a)}{b-a}\leq \gamma_1.
\]
Hence $g(x)\leq\max\{\gamma_1,\gamma_2\}$ for $a<x<b$.\\
(ii)
By Theorem \ref{Thm:CharacterizingThm}(ii), it follows that for $x\leq a$
\[
G(x)
\leq G(a)\left(1+s^*\frac{g(a)}{G(a)}(x-a)\right)^{1/s^*}_+
=\left(G^{s^*}(a)+s^*\frac{g(a)}{G^{1-s^*}(a)}(x-a)\right)^{1/s^*}_+.
\]
By Theorem \ref{Thm:CharacterizingThm}(iii), the non-increasing monotonicity of $g/G^{1-s^*}$ implies that 
\[
\frac{g(a)}{G^{1-s^*}(a)}
=\left(\frac{1}{s^*}G^{s^*}\right)^\prime(a)
\geq \frac{\frac{1}{s^*}G^{s^*}(b)-\frac{1}{s^*}G^{s^*}(a)}{b-a}
\geq \frac{\frac{1}{s^*}v^{s^*}-\frac{1}{s^*}u^{s^*}}{b-a}=\gamma_1.
\]
The last inequality follows from noting that $G(a)\leq U_n(a)\leq u$ and $G(b)\geq L_n(b)\geq v$. 
Since $x-a\leq0$, it follows that
\begin{eqnarray*}
G(x)
&\leq&\left(G^{s^*}(a)+s^*\frac{g(a)}{G^{1-s^*}(a)}(x-a)\right)^{1/s^*}_+\\
&\leq&\left(G^{s^*}(a)+s^*\gamma_1(x-a)\right)^{1/s^*}_+\\
&\leq&\left(u^{s^*}+s^*\gamma_1(x-a)\right)^{1/s^*}_+.
\end{eqnarray*}
The last inequality follows from noting that $G(a)\leq u$.\\
On the other hand, by Theorem \ref{Thm:CharacterizingThm}(ii), it follows that for $x\geq b$
\begin{eqnarray*}
1-G(x)
&\leq&(1-G(b))\left(1-s^*\frac{g(b)}{1-G(b)}(x-b)\right)_+^{1/s^*}\\
&=&\left((1-G(b))^{s^*}-s^*\frac{g(b)}{(1-G(b))^{1-s^*}}(x-b)\right)_+^{1/s^*}\\
&\leq&\left((1-v)^{s^*}-s^*\frac{g(b)}{(1-G(b))^{1-s^*}}(x-b)\right)_+^{1/s^*}.
\end{eqnarray*}
The last inequality follows from noting that $1-G(b)\leq 1-v$.
By Theorem \ref{Thm:CharacterizingThm}(iii), the non-decreasing monotonicity of $g/(1-G)^{1-s^*}$ implies that 
\begin{eqnarray*}
\frac{g(b)}{(1-G)^{1-s^*}(b)}
&=&(\frac{-1}{s^*}(1-G)^{s^*})^\prime(b)\\
&\geq&\frac{\frac{-1}{s^*}(1-G(b))^{s^*}-\frac{-1}{s^*}(1-G(a))^{s^*}}{b-a}\\
&=&\frac{\frac{1}{s^*}\left((1-G(a))^{s^*}-(1-G(b))^{s^*}\right)}{b-a}\\
&\geq&\frac{\frac{1}{s^*}\left((1-u)^{s^*}-(1-v)^{s^*}\right)}{b-a}\\
&=&\gamma_2 \ .
\end{eqnarray*}
The last inequality follows from noting that $G(a)\leq U_n(a)\leq u$ and $G(b)\geq L_n(b)\geq v$. 
Since $x-b\geq0$, it follows that 
\[
1-G(x)
\leq
\left((1-v)^{s^*}-s^*\gamma_2(x-b)\right)_+^{1/s^*}. 
\]
\hfill $\Box$
\medskip

\par\noindent
\textbf{Proof of Theorem \ref{Theorem:ConsistencyOfConfidenceBands}:}\\
The following proof is analogous to the proof of Theorem 3 in \cite{DUMBGEN20171}, 
in which they proved the result in the case $s^*=0$.
In the following proof we assume that $s^*\neq 0$.\\
(i) Suppose $s^*>0$. 
Since $F$ is not bi-$s^*$-concave, it follows that $F^{s^*}$ or $(1-F)^{s^*}$ is not concave. 
Without loss of generality, we assume that $F^{s^*}$ is not concave 
and hence there exist real numbers $x_0<x_1<x_2$ such that 
$F^{s^*}(x_1)<(1-\lambda)F^{s^*}(x_0)+\lambda F^{s^*}(x_2)$, where $\lambda\equiv (x_1-x_0)/(x_2-x_0)\in(0,1)$.
By the consistency of $L_n$ and $U_n$, it follows that, with 
probability tending to one, $U_n^{s^*}(x_1)<(1-\lambda)L_n^{s^*}(x_0)+\lambda L_n^{s^*}(x_2)$ and hence
\[
G^{s^*}(x_1)
<
(1-\lambda)G^{s^*}(x_0)+\lambda G^{s^*}(x_2),
\]
for any $G$ such that $L_n\leq G\leq U_n$.
Therefore, there are no bi-$s^*$-concave distribution functions 
fitting between $L_n$ and $U_n$ and hence $L_n^o=1$ and $U_n^o=0$ with probability tending to one.\\
The proof of the case $s^*<0$ is similar and hence is omitted.\\
(ii) Suppose $F\in\mathcal{P}_{s^*}$. 
Note that since $(L_n, U_n)$ is a $(1-\alpha)$ confidence band for $F$, 
it follows that $P(L_n^o\leq U_n^o)\geq P(L_n\leq F\leq U_n)\geq 1-\alpha$.\\
If $\{G\in\mathcal{P}_{s^*}:L_n\leq G\leq U_n\}$ is empty, it follows 
that $L_n^o=1$ and $U_n^o=0$ and hence the assertions are trivial.
In the following proof, we assume that $\{G\in\mathcal{P}_{s^*}:L_n\leq G\leq U_n\}$ 
is not empty.\\
To prove (\ref{Formula:ConfidenceBand1}), we first prove that 
$\|L_n-F\|_{\infty}\rightarrow_p0$ and $\|U_n-F\|_{\infty}\rightarrow_p0$.
By the continuity of $F$, for any $m\in\mathbb{N}^+$ with $m\geq 2$, 
there exist real numbers $\{x_i\}_{i=1}^{m-1}$ such that $F(x_i)=i/m$, $i=1,\ldots,m-1$.
Furthermore, define $x_0=-\infty$ and $x_m=\infty$.\\
By the non-decreasing monotonicity of $L_n$ and $F$, it follows that for $x\in[x_{i-1},x_i]$
\[
L_n(x)-F(x)\leq L_n(x_i)-F(x_{i-1})=L_n(x_i)-(F(x_i)-\frac{1}{m})=L_n(x_i)-F(x_i)+\frac{1}{m},
\]
and
\begin{flalign*}
L_n(x)-F(x)&\geq L_n(x_{i-1})-F(x_{i})&&\\
&=L_n(x_{i-1})-(F(x_{i-1})+\frac{1}{m})=L_n(x_{i-1})-F(x_{i-1})-\frac{1}{m}.&&
\end{flalign*}
Hence 
\[
|L_n(x)-F(x)|\leq \max_{i=1,\ldots,m-1}|L_n(x_{i})-F(x_{i})|+\frac{1}{m}
\]
for $x\in[x_{i-1},x_i]$.
Note that 
\begin{eqnarray*}
\|L_n-F\|_{\infty}
&=&\sup_{x\in\mathbb{R}}|L_n(x)-F(x)|=\max_{i=1,\ldots,m}\sup_{x\in[x_{i-1},x_i]}|L_n(x)-F(x)|,
\end{eqnarray*}
it follows that 
\[
\|L_n-F\|_\infty\leq\max_{i=1,\ldots,m-1}|L_n(x_{i})-F(x_{i})|+\frac{1}{m},
\]
and hence pointwise convergence implies uniform convergence.
An analogous proof shows that $\|U_n-F\|_\infty\rightarrow_p0$ and is omitted.\\
Combining $\|L_n-F\|_{\infty}\rightarrow_p0$ and $\|U_n-F\|_{\infty}\rightarrow_p0$ implies that
\[
\sup_{G\in\mathcal{P}_{s^*}:L_n\leq G\leq U_n}\|G-F\|_{\infty}
\leq
\|L_n-F\|_\infty+\|U_n-F\|_\infty
\rightarrow_p0.
\]
To prove (\ref{Formula:ConfidenceBand2}) in the case that $h_G=\left(G^{s^*}\right)^\prime$,
it suffices to prove that 
\begin{eqnarray}
\label{Formula:ProofOfTheorem2_1}
\sup_{G\in\mathcal{P}_{s^*}:L_n\leq G\leq U_n} 
\left\|\left(G^{s^*}/s^*\right)^\prime-\left(F^{s^*}/s^*\right)^\prime\right\|_{K,\infty}\rightarrow_p0.
\end{eqnarray}
Note that $h_G/s^*=G^\prime/G^{1-s^*}$.
Since $K$ is a compact interval in $J(F)$ and $h_F/s^*=f/F^{1-s^*}$ is continuous and non-increasing on $J(F)$,
for any fixed $\epsilon>0$ there exist points $a_0<a_1<\dots<a_m<a_{m+1}$ in $J(F)$ such that $K\subset [a_1,a_m]$ and
\[
0\leq \frac{1}{s^*}h_F(a_{i-1})-\frac{1}{s^*}h_F(a_i)\leq\epsilon\text{ for }1\leq i\leq m+1.
\]
For $G\in\mathcal{P}_{s^*}$ with $L_n\leq G\leq U_n$, for any $x\in K$ 
it follows from the monotonicity of $h_F/s^*$ and $h_G/s^*$ that 
\begin{eqnarray*}
\sup_{x\in K}\left(\frac{1}{s^*}h_G(x)-\frac{1}{s^*}h_F(x)\right)
&\leq&
\max_{i=1,\ldots,m-1}\left(\frac{1}{s^*}h_G(a_i)-\frac{1}{s^*}h_F(a_{i+1})\right)\\
&\leq&
\max_{i=1,\ldots,m-1}\left(\frac{\frac{1}{s^*}G^{s^*}(a_i)-\frac{1}{s^*}G^{s^*}(a_{i-1})}{a_i-a_{i-1}}\right.\\
& &\left.\ \ \ \ \ \ \ \ \ \ \ \ \ \ \ \ -\frac{1}{s^*}h_F(a_{i+1})\right)\\
&\leq&
\max_{i=1,\ldots,m-1}\left(\frac{\frac{1}{s^*}U_n^{s^*}(a_i)-\frac{1}{s^*}L_n^{s^*}(a_{i-1})}{a_i-a_{i-1}}\right.\\
& &\left.\ \ \ \ \ \ \ \ \ \ \ \ \ \ \ \ -\frac{1}{s^*}h_F(a_{i+1})\right)\\
&=&
\max_{i=1,\ldots,m-1}\left(\frac{\frac{1}{s^*}F^{s^*}(a_i)-\frac{1}{s^*}F^{s^*}(a_{i-1})}{a_i-a_{i-1}}\right.\\
& &\left.\ \ \ \ \ \ \ \ \ \ \ \ \ \ \ \ -\frac{1}{s^*}h_F(a_{i+1})\right)+o_p(1)\\
&\leq&
\max_{i=1,\ldots,m-1}\left(\frac{1}{s^*}h_F(a_{i-1})-\frac{1}{s^*}h_F(a_{i+1})\right)\\
& &+o_p(1)\\
&\leq&
2\epsilon+o_p(1).
\end{eqnarray*}
Analogously,
\begin{eqnarray*}
\sup_{x\in K}\left(\frac{1}{s^*}h_F(x)-\frac{1}{s^*}h_G(x)\right)
&\leq&
\max_{i=1,\ldots,m-1}\left(\frac{1}{s^*}h_F(a_i)-\frac{1}{s^*}h_F(a_{i+2})\right)+o_p(1)\\
&\leq& 2\epsilon+o_p(1).
\end{eqnarray*}
Since $\epsilon>0$ is arbitrarily small, this shows that (\ref{Formula:ProofOfTheorem2_1}) holds.\\
The proof of (\ref{Formula:ConfidenceBand2}) in the case that 
$h_G=\left((1-G)^{s^*}\right)^\prime$ is similar and hence is omitted.\\
Since $G^\prime=G^{1-s^*}\left(G^{s^*}/s^*\right)^\prime$, it follows from 
(\ref{Formula:ProofOfTheorem2_1}) that (\ref{Formula:ConfidenceBand2}) holds in the case that $h_G=G^\prime$.\\
Finally, let $x_1<\sup J(F)$ and $b_1<f(x_1)/F^{1-s^*}(x_1)$. 
As in the proof of Lemma \ref{Lemma:LemmaInConfidenceBands}(ii) an 
analogous argument  implies that for any $x_1^\prime>x_1$, $x_1^\prime\in J(F)$,
\begin{eqnarray*}
U_n^o(x)
&\leq&\left(U_n^{s^*}(x^\prime)+s^*\frac{\frac{1}{s^*}L_n^{s^*}(x_1^\prime)
      -\frac{1}{s^*}U_n^{s^*}(x_1)}{x_1^\prime-x_1}(x-x^\prime)\right)^{1/s^*}_+
\end{eqnarray*}
for all $x\leq x^\prime\leq x_1$. \\
Note that by the consistency of $L_n$ and $U_n$ and letting $x_1^\prime\downarrow x_1$, it follows that.
\[
\frac{\frac{1}{s^*}L_n^{s^*}(x_1^\prime)-\frac{1}{s^*}U_n^{s^*}(x_1)}{x_1^\prime-x_1}
\rightarrow_p
\frac{\frac{1}{s^*}F^{s^*}(x_1^\prime)-\frac{1}{s^*}F^{s^*}(x_1)}{x_1^\prime-x_1}
>b_1.
\]
Hence with probability tending to one, 
\[
U_n^o(x)\leq\left(U_n^{s^*}(x^\prime)+s^*b_1(x-x^\prime)\right)^{1/s^*}_+,
\]
for all $x\leq x^\prime\leq x_1$. The proof of (\ref{Formula:ConfidenceBand4}) is similar and hence is omitted.
\hfill $\Box$
\medskip

\par\noindent
\textbf{Proof of Remark 1:}\\
(i) By Theorem \ref{Thm:CharacterizingThm}(ii), 
if $s^*>0$ and $\inf J(F)=-\infty$, it follows that for arbitrary $x\in J(F)$,
\[
F(y)\leq F(x)\cdot\left(1+s^*\frac{f(x)}{F(x)}(y-x)\right)_+^{1/s^*}=0
\]
for small enough $y$ such that 
\[
1+s^*\frac{f(x)}{F(x)}(y-x)<0.
\]
This violates the assumption that $\inf J(F)=-\infty$ and hence $\inf J(F)>-\infty$.\\
The finiteness of $\sup J(F)$ can be proved similarly and hence is omitted.\\
(ii) We first note that (\ref{Formula:Remark1}) holds automatically if $\inf J(F)>-\infty$ and $\sup J(F)<\infty$.\\
In the following proof, we focus on the case that $\inf J(F)=-\infty$ and $\sup J(F)<\infty$. 
To prove (\ref{Formula:Remark1}), it suffices to show that 
$\int |x|^t dF(x)$ is finite for $t\in (0,(-1)/s^*)$.\\
Note that 
\begin{eqnarray*}
\int |x|^t dF(x)
&=&E|X|^t=\int_0^\infty P(|X|^t>a)da=\int_0^\infty P(|X|>a^{1/t})da\\
&=&\int_0^\infty ta^{t-1}P(|X|>a)da\\
&=&\int_0^\infty ta^{t-1}P(X>a)da+\int_0^\infty ta^{t-1}P(X<-a)da.
\end{eqnarray*}
Since $\sup J(F)$ is finite, the first term of the last display is finite and hence it 
suffices to prove that $ ta^{t-1}P(X<-a)$ is integrable for $t<(-1)/s^*$.\\
It follows from Theorem \ref{Thm:CharacterizingThm}(ii) that for any 
$a$ large enough and $x\in J(F)$,
\[
P(X<-a)
\leq
F(x)\left(1+\frac{s^*f(x)(-a-x)}{F(x)}\right)^{1/s^*}_+
=F(x)\left(1+\frac{-s^*f(x)(a+x)}{F(x)}\right)^{1/s^*}_+.
\]
Thus $ta^{t-1}P(X<-a)$ is integrable for $t<(-1)/s^*$, 
since
\begin{eqnarray*}
ta^{t-1}P(X<-a)
&\leq&
tF(x)a^{t-1}\left(1+\frac{-s^*f(x)(a+x)}{F(x)}\right)^{1/s^*}_+\\
&=&tF(x)\left(\frac{-s^*f(x)}{F(x)}\right)^{1/s^*}a^{t-1}\left(a+x+\frac{F(x)}{-s^*f(x)}\right)^{1/s^*}_+\\
&\leq&2tF(x)\left(\frac{-s^*f(x)}{F(x)}\right)^{1/s^*}a^{t+1/s^*-1}
\end{eqnarray*}
for $a$ large enough and $a^{t+1/s^*-1}$ is integrable for $t<(-1)/s^*$.\\
For other cases, the proof is similar and hence is omitted.
\hfill $\Box$
\medskip

\par\noindent
\textbf{Proof of Corollary \ref{Corollary:CorollaryInConfidenceBands}:}\\
Suppose that $x_0$ is a point in $J(F)$. 
Notice that for any $z\in\mathbb{R}$, 
\[
\phi(z)-\phi(x_0)=\int_{\mathbb{R}}\left(1_{[x_0\leq x< z]}-1_{[z\leq x< x_0]}\right)\phi^\prime(x)dx,
\]
and hence by Fubini's theorem, it follows that
\begin{eqnarray}
\label{Formula1:ProofOfCorollary3}
\int_{\mathbb{R}}\phi dG=\phi(x_0)+\int_{\mathbb{R}}\phi^\prime(x)\left(1_{[x\geq x_0]}-G(x)\right)dx,
\end{eqnarray}
provided that 
\[
\int_{\mathbb{R}}\left|\phi^\prime(x)\right|\left| 1_{[x\geq x_0]}-G(x)\right|dx<\infty.
\]
To prove the last display, note that for any $b_1\in (0,T_1(F))$ and $b_2\in(0,T_2(F))$, 
there exist points $x_1,x_2\in J(F)$ with $x_1\leq x_0\leq x_2$ and 
\[
\frac{f}{F^{1-s^*}}(x_1)>b_1,\ \frac{f}{(1-F)^{1-s^*}}(x_2)>b_2.
\]
Then it follows from Theorem \ref{Theorem:ConsistencyOfConfidenceBands}(ii) that with probability tending to one,
\[
U_n^o(x)\leq\left(U_n^{s^*}(x_1)+s^*b_1(x-x_1)\right)^{1/s^*}_+
\text{ for }x\leq x_1,
\]
and
\[
1-L_n^o(x)\leq \left((1-L_n(x_2))^{s^*}-s^*b_2(x-x_2)\right)^{1/s^*}_+
\text{ for }x\geq x_2.
\]
Hence for any $c>\max\{|x_1|,|x_2|\}$, it follows that 
\begin{eqnarray*}
\int_{-\infty}^{x_1-c}\left|\phi^\prime(x)\right|\left| 1_{[x\geq x_0]}-G(x)\right|dx
&=&
\int_{-\infty}^{x_1-c}\left|\phi^\prime(x)\right|G(x)dx\\
&\leq&
\int_{-\infty}^{x_1-c}\left|\phi^\prime(x)\right|U_n^o(x)dx\\
&\leq&\int_{-\infty}^{x_1-c}\left|\phi^\prime(x)\right|\left(U_n^{s^*}(x_1)+s^*b_1(x-x_1)\right)^{1/s^*}_+dx\\
&=&\int_{-\infty}^{x_1-c}\left|\phi^\prime(x)\right|\left(U_n^{s^*}(x_1)+s^*b_1(x-x_1)\right)^{1/s^*}dx.
\end{eqnarray*}
Since $|\phi^\prime(x)|\leq a|x|^{k-1}$, it follows that the last display is no larger than 
\[
\int_{-\infty}^{x_1-c}a|x|^{k-1}\left(U_n^{s^*}(x_1)+s^*b_1(x-x_1)\right)^{1/s^*}dx,
\]
which is finite by noting that $k-1+1/s^*<-1$.
Analogously, one can prove that for $c>\max\{|x_1|,|x_2|\}$,
\[
\int_{x_2+c}^{\infty}\left|\phi^\prime(x)\right|\left| 1_{[x\geq x_0]} 
  -G(x)\right|dx\leq \int_{x_2+c}^{\infty}\left|\phi^\prime(x)\right|\left| 1-L_n^o(x)\right|dx<\infty.
\] 
Since $\phi^\prime$ is continuous on $\mathbb{R}$, it follows that for any $c>\max\{|x_1|,|x_2|\}$,
\[
\int_{x_1-c}^{x_2+c}\left|\phi^\prime(x)\right|\left| 1_{[x\geq x_0]}-G(x)\right|dx<\infty
\]
and hence
\[
\int_{\mathbb{R}}\left|\phi^\prime(x)\right|\left| 1_{[x\geq x_0]}-G(x)\right|dx<\infty.
\]
By (\ref{Formula1:ProofOfCorollary3}), it follows that
\[
\sup_{G:L_n^o\leq G\leq U_n^o}\left|\int\phi dG-\int\phi dF\right|=
\sup_{G:L_n^o\leq G\leq U_n^o}\left|\int\phi^\prime(x)(F-G)(x)dx\right|,
\]
which is not larger than
\begin{eqnarray*}
& &\sup_{G:L_n^o\leq G\leq U_n^o}\|G-F\|_\infty\int_{x_1-c}^{x_2+c}|\phi^\prime(x)|dx\\
& &\ \ \ \ \ \ +
\int_{-\infty}^{x_1-c}|\phi^\prime(x)|(F+U_n^o)(x)dx
+
\int_{x_2+c}^{\infty}|\phi^\prime(x)|(1-F+1-L_n^o)(x)dx\\
&\leq&
o_p(1)+2\int_{-\infty}^{x_1-c}|\phi^\prime(x)|U_n^o(x)dx+2\int^{\infty}_{x_2+c}|\phi^\prime(x)|(1-L_n^o(x))dx.
\end{eqnarray*}
Note that the last two terms go to zero as $c$ goes to infinity by their integrability and hence 
\[
\sup_{G:L_n^o\leq G\leq U_n^o}\left|\int\phi dG-\int\phi dF\right|=o_p(1).
\]
\hfill $\Box$
\medskip

\par\noindent
\textbf{Proof of Theorem \ref{Theorem:RateOfConfidenceBands}:}
It follows from the proof of Corollary \ref{Corollary:CorollaryInConfidenceBands} that 
\[
\sup_{G:L_n^o\leq G\leq U_n^o}\left|\int \phi dG-\int \phi dF\right|
=
\sup_{G:L_n^o\leq G\leq U_n^o}\left|\int\phi^\prime(x)(F-G)(x)dx\right|
\]
and hence 
\[
\sup_{G:L_n^o\leq G\leq U_n^o}\left|\int \phi dG-\int \phi dF\right|
\leq
\sup_{G:L_n^o\leq G\leq U_n^o}\int\left|\phi^\prime(x)\right|\left|(G-F)(x)\right|dx.
\]
It suffices to bound $|G-F|$ on $\mathbb{R}$, where $G$ is between $L_n^o$ and $U_n^o$.\\
It follows from $G\leq U_n^o\leq U_n$ and Condition (*) that on the interval 
$
\{\lambda n^{-1/(2-2\gamma)}\leq \FF_n \leq 1-\lambda n^{-1/(2-2\gamma)}\},
$
\[
G-F
\leq U_n^o-F
\leq U_n-F
\leq U_n-\FF_n+\FF_n-F
\leq \kappa n^{-1/2}\left(\FF_n\left(1-\FF_n\right)\right)^\gamma
+\left|\FF_n-F\right|
\]
To bound $\left|\FF_n -F\right|$, it follows from Theorem 3.7.1, page 141, \cite{MR3396731} that
\[
\left\|\frac{\sqrt{n}\left(\FF_n -F\right)-\mathbb{U}\circ F}{\left(F(1-F)\right)^\gamma}\right\|\rightarrow_p 0
\]
by verifying that $q(t)\equiv (t(1-t))^\gamma$ with $0\leq\gamma<1/2$ is monotonically increasing on 
$[0,1/2]$, symmetric about $1/2$ and $\int_{0}^1q^{-2}(t)dt<\infty$, where $\mathbb{U}$ is Brownian bridge on $[0,1]$.\\
Hence for any fixed $\epsilon\in(0,1)$ there exists a constant $\kappa_\epsilon>0$ such that with probability at least $1-\epsilon$, 
\[
\left|\FF_n-F\right|\leq \kappa_\epsilon n^{-1/2}\left(F(1-F)\right)^\gamma
\]
on $\mathbb{R}$.
Thus, 
it follows that on the interval 
$
\{\lambda n^{-1/(2-2\gamma)}\leq \FF_n\leq 1-\lambda n^{-1/(2-2\gamma)}\},
$
\[
G-F
\leq 
\kappa n^{-1/2}\left(\FF_n\left(1-\FF_n \right)\right)^\gamma
+
\kappa_\epsilon n^{-1/2}\left(F(1-F)\right)^\gamma.
\]
To bound $\FF_n \left(1-\FF_n \right)$ by $F(1-F)$,
note that 
\begin{eqnarray*}
\FF_n\left(1-\FF_n \right)
&=&
\left(\FF_n-F+F\right)\left(1-F+F-\FF_n\right)\\
&=&\left(\FF_n-F\right)(1-F)+F(1-F)-\left(\FF_n -F\right)^2-F\left(\FF_n-F\right)\\
&=&F(1-F)+\left(\FF_n -F\right)(1-2F)-\left(\FF_n -F\right)^2\\
&\leq&F(1-F)+\left|\FF_n-F\right||1-2F|+\left|\FF_n-F\right|^2\\
&\leq&F(1-F)+\left|\FF_n-F\right|+\left|\FF_n-F\right|\\
&=&F(1-F)\cdot\left(1+\frac{2\left|\FF_n -F\right|}{F(1-F)}\right)\\
&\leq&F(1-F)\cdot\left(1+\frac{4\left|\FF_n-F\right|}{\min\{F,1-F\}}\right)\\
& &\ \ \ \ \ \ \ \ \text{ since }F(1-F)\geq \min\{F,1-F\}/2,\\
&\leq&F(1-F)\cdot\left(1+\frac{4\kappa_\epsilon n^{-1/2}\left(F(1-F)\right)^\gamma}{\min\{F,1-F\}}\right).
\end{eqnarray*}
For a constant $\lambda_\epsilon>0$ to be specified later, it follows from
$\lambda_\epsilon n^{-1/(2-2\gamma)}\leq F\leq 1- \lambda_\epsilon n^{-1/(2-2\gamma)}$ and $\gamma\in[0,1/2)$ that 
\[
\frac{\left(F(1-F)\right)^\gamma}{F}
=F^{\gamma-1}(1-F)^\gamma
\leq \lambda_\epsilon^{\gamma-1} n^{-(\gamma-1)/(2-2\gamma)}
=\lambda_\epsilon^{\gamma-1} n^{1/2}
\]
and
\[
\frac{\left(F(1-F)\right)^\gamma}{1-F}
=F^{\gamma}(1-F)^{\gamma-1}
\leq \lambda_\epsilon^{\gamma-1} n^{-(\gamma-1)/(2-2\gamma)}
=\lambda_\epsilon^{\gamma-1} n^{1/2}.
\]
Hence
\[
\FF_n\left(1-\FF_n \right)
\leq
F(1-F)\cdot\left(1+4\kappa_\epsilon n^{-1/2}\lambda_\epsilon^{\gamma-1} n^{1/2}\right)
=F(1-F)(1+4\kappa_\epsilon\lambda_\epsilon^{\gamma-1}).
\]
Thus, on the interval
\[
\{\lambda n^{-1/(2-2\gamma)}\leq \FF_n \leq 1-\lambda n^{-1/(2-2\gamma)}\}\cap
\{\lambda_\epsilon n^{-1/(2-2\gamma)}\leq F\leq 1-\lambda_\epsilon n^{-1/(2-2\gamma)}\},
\]
\begin{eqnarray*}
G-F
&\leq& 
\kappa n^{-1/2}\left(F(1-F)(1+4\kappa_\epsilon\lambda_\epsilon^{\gamma-1} )\right)^\gamma
+
\kappa_\epsilon n^{-1/2}\left(F(1-F)\right)^\gamma\\
&=&\nu_\epsilon n^{-1/2}\left(F(1-F)\right)^\gamma,
\end{eqnarray*}
where $\nu_\epsilon=\kappa(1+4\kappa_\epsilon\lambda_\epsilon^{\gamma-1} )^\gamma+\kappa_\epsilon$.\\
The following arguments show that for a large enough $\lambda_\epsilon$, the interval 
$\{\lambda_\epsilon n^{-1/(2-2\gamma)}\leq F\leq 1-\lambda_\epsilon n^{-1/(2-2\gamma)}\}$ is 
a subset of $\{\lambda n^{-1/(2-2\gamma)}\leq \FF_n \leq 1-\lambda n^{-1/(2-2\gamma)}\}$.
To see this, note that 
\begin{eqnarray*}
\FF_n 
&=&F+\FF_n -F\\
&\geq& \left(1-\frac{|\FF_n -F|}{F}\right)F\\
&\geq& \left(1-\kappa_\epsilon n^{-1/2}F^{\gamma-1}(1-F)^\gamma\right)F\\
&\geq& \left(1-\kappa_\epsilon n^{-1/2}\lambda_\epsilon^{\gamma-1}n^{1/2}\right)\lambda_\epsilon n^{-1/(2-2\gamma)}\\
&=& \left(\lambda_\epsilon-\kappa_\epsilon \lambda_\epsilon^{\gamma}\right) n^{-1/(2-2\gamma)}
\end{eqnarray*}
and analogously,
\begin{eqnarray*}
1-\FF_n
&\geq& \left(\lambda_\epsilon-\kappa_\epsilon \lambda_\epsilon^{\gamma}\right) n^{-1/(2-2\gamma)},
\end{eqnarray*}
it follows that by choosing a $\lambda_\epsilon$ large enough such that 
$\lambda_\epsilon-\kappa_\epsilon \lambda_\epsilon^{\gamma}>\lambda$,
the interval $\{\lambda_\epsilon n^{-1/(2-2\gamma)}\leq F\leq 1-\lambda_\epsilon n^{-1/(2-2\gamma)}\}$ 
is a subset of $\{\lambda n^{-1/(2-2\gamma)}\leq \FF_n \leq 1-\lambda n^{-1/(2-2\gamma)}\}$ and hence
on the interval
\[
\{\lambda_\epsilon n^{-1/(2-2\gamma)}\leq F\leq 1-\lambda_\epsilon n^{-1/(2-2\gamma)}\},
\]
\begin{eqnarray*}
G-F
&\leq&\nu_\epsilon n^{-1/2}\left(F(1-F)\right)^\gamma.
\end{eqnarray*}
Define $x_{n1}$ and $x_{n2}$, such that 
$F(x_{n1})=\lambda_\epsilon n^{-1/(2-2\gamma)}$ and  
$F(x_{n2})=1-\lambda_\epsilon n^{-1/(2-2\gamma)}$. 
Analogously, one can prove that 
$
F-G
\leq\nu_\epsilon n^{-1/2}\left(F(1-F)\right)^\gamma
$
on $[x_{n1},x_{n2}]$
and hence 
\begin{eqnarray}
|G-F|
\leq\nu_\epsilon n^{-1/2}\left(F(1-F)\right)^\gamma
\label{Formula3:ProofOfTheorem4}
\end{eqnarray}
on $[x_{n1},x_{n2}]$.
Thus for $G$ between $L_n^o$ and $U_n^o$,
\begin{eqnarray*}
\sup_{G:L_n^o\leq G\leq U_n^o}\left|\int \phi d(G-F)\right|
&=&\sup_{G:L_n^o\leq G\leq U_n^o}\left|\int \phi^\prime(x)\left(F(x)-G(x)\right) dx\right|\\
&\leq&
\nu_\epsilon n^{-1/2}\int_{x_{n1}}^{x_{n2}}\left|\phi^\prime(x)\right|F^\gamma(x)(1-F(x))^\gamma dx\\
& &\ \ \ \ +\int_{-\infty}^{x_{n1}}\left|\phi^\prime(x)\right|(F(x)+U_n^o(x))dx\\
& &\ \ \ \ +\int_{x_{n2}}^{\infty}\left|\phi^\prime(x)\right|(2-F(x)-L_n^o(x))dx.
\end{eqnarray*}
From here, we can see that if $F\in \mathcal{P}_{s^*}$ with $s^*>0$, it follows 
from Remark 1(i) that $J(F)$ is bounded and hence 
\[
\sup_{G:L_n^o\leq G\leq U_n^o}\left|\int \phi d(G-F)\right|=O_p(n^{-1/2})
\]
as long as $\phi^\prime$ is bounded on $J(F)$.\\
The similar argument works if $F\in \mathcal{P}_{s^*}$ with $s^*<0$ and $J(F)$ is bounded.
In the following proof, we get back to our case that $F\in \mathcal{P}_{s^*}$ with 
$s^*<0$ and without loss of generality, we assume $J(F)=(-\infty,\infty)$.\\
As in the proof of Corollary \ref{Corollary:CorollaryInConfidenceBands}, 
for $x_0\in J(F)$, $b_1\in(0,T_1(F))$ and $b_2\in(0,T_2(F))$, there exist points 
$x_1,x_2\in J(F)$ with $x_1<x_0<x_2$ such that $f(x_1)/F^{1-s^*}(x_1)>b_1$ and $f(x_2)/(1-F(x_2))^{1-s^*}>b_2$.
Then it follows from Theorem \ref{Theorem:ConsistencyOfConfidenceBands}(ii) that with asymptotic probability one,
\begin{eqnarray}
U_n^o(x)
&\leq&\left(U_n^{s^*}(x_1)+s^*b_1(x-x_1)\right)^{1/s^*}_+\notag\\
&=&U_n(x_1)\left(1+\frac{s^*b_1}{U_n^{s^*}(x_1)}(x-x_1)\right)^{1/s^*}
\text{ for }x\leq x_1,\label{Formula2:ProofOfTheorem4}
\end{eqnarray}
and
\begin{eqnarray*}
1-L_n^o(x)
&\leq&\left((1-L_n(x_2))^{s^*}-s^*b_2(x-x_2)\right)^{1/s^*}_+\\
&=&\left(1-L_n(x_2)\right)\left(1-\frac{s^*b_2}{(1-L_n(x_2))^{s^*}}(x-x_2)\right)^{1/s^*}
\text{ for }x\geq x_2.
\end{eqnarray*}
Similarly, it follows from Theorem \ref{Thm:CharacterizingThm}(ii) that
\begin{eqnarray}
F(x)
&\leq&F(x_1)\left(1+s^*\frac{f(x_1)}{F(x_1)}(x-x_1)\right)^{1/s^*}_+
\notag\\
&\leq&F(x_1)\left(1+\frac{s^*b_1}{F^{s^*}(x_1)}(x-x_1)\right)^{1/s^*}
\text{ for }x\leq x_1,\label{Formula1:ProofOfTheorem4}
\end{eqnarray}
and
\begin{eqnarray*}
1-F(x)
&\leq&(1-F(x_2))\left(1-s^*\frac{f(x_2)}{1-F(x_2)}(x-x_2)\right)^{1/s^*}_+\\
&\leq&(1-F(x_2))\left(1-\frac{s^*b_2}{\left(1-F(x_2)\right)^{s^*}}(x-x_2)\right)^{1/s^*}\text{ for }x\geq x_2.
\end{eqnarray*}
For large enough $n$, one can have $[x_1,x_2]\subset[x_{n1},x_{n2}]$ and hence 
\[
\sup_{G:L_n^o\leq G\leq U_n^o}\left|\int \phi d(G-F)\right|\leq I_{n0}+I_{n1}+I_{n1}^\prime+I_{n2}+I_{n2}^\prime,
\]
where
\[
I_{n0}\equiv\nu_\epsilon n^{-1/2}\int_{x_1}^{x_2}\left|\phi^\prime(x)\right|F^\gamma(x)(1-F(x))^\gamma dx,
\]
\[
I_{n1}\equiv\nu_\epsilon n^{-1/2}\int_{x_n1}^{x_1}\left|\phi^\prime(x)\right|F^\gamma(x)(1-F(x))^\gamma dx,
\]
\[
I_{n2}\equiv\nu_\epsilon n^{-1/2}\int_{x_2}^{x_{n2}}\left|\phi^\prime(x)\right|F^\gamma(x)(1-F(x))^\gamma dx,
\]
\[
I_{n1}^\prime\equiv\int_{-\infty}^{x_{n1}}\left|\phi^\prime(x)\right|(F(x)+U_n^o(x))dx,
\]
\[
I_{n2}^\prime\equiv\int_{x_{n2}}^{\infty}\left|\phi^\prime(x)\right|(2-F(x)-L_n^o(x))dx.
\]
Note that 
$I_{n0}\leq\nu_\epsilon n^{-1/2}\int_{x_1}^{x_2}\left|\phi^\prime(x)\right|dx=O(n^{-1/2})$.
For the other terms, first note that $F(x_{n1})=\lambda_\epsilon n^{-1/(2-2\gamma)}$ 
and hence it follows from (\ref{Formula1:ProofOfTheorem4}) that
\[
x_{n1}\geq x_1-\frac{F^{s^*}(x_1)}{s^*b_1}+\frac{\lambda_\epsilon^{s^*}}{s^*b_1}n^{-s^*/(2-2\gamma)}
=O(1)+\frac{\lambda_\epsilon^{s^*}}{s^*b_1}n^{-s^*/(2-2\gamma)}.
\]
Analogously, one can prove that 
\[
x_{n2}\leq x_2-\frac{(1-F(x_2))^{s^*}}{s^*b_2}
-\frac{\lambda_\epsilon^{s^*}}{s^*b_2}n^{-s^*/(2-2\gamma)}
=O(1)+\frac{\lambda_\epsilon^{s^*}}{s^*b_1}n^{-s^*/(2-2\gamma)}.
\]
Thus, it follows from (\ref{Formula1:ProofOfTheorem4}) and the upper bound of $|\phi^\prime|$ that 
\begin{eqnarray*}
I_{n1}
&\leq&\nu_\epsilon n^{-1/2}\int_{x_n1}^{x_1}\left|\phi^\prime(x)\right|F^\gamma(x)dx\\
&\leq&O\left(n^{-1/2}\int_{x_{n1}}^{x_1}\left|\phi^\prime(x) 
\right|\left(1+\frac{s^*b_1}{F^{s^*}(x_1)}(x-x_1)\right)^{\gamma/s^*}dx\right)\\
&\leq&O\left(n^{-1/2}\int_{0}^{O\left(n^{-s^*/(2-2\gamma)}\right)}
          \left|\phi^\prime(x)\right|\left(1+\frac{-s^*b_1}{F^{s^*}(x_1)}x\right)^{\gamma/s^*}dx\right)\\
&=&O\left(n^{-1/2}\int_{0}^{O\left(n^{-s^*/(2-2\gamma)}\right)}\left|\phi^\prime(x)\right|x^{\gamma/s^*}dx\right)\\
&\leq&O\left(n^{-1/2}\int_{0}^{O\left(n^{-s^*/(2-2\gamma)}\right)}x^{k-1}x^{\gamma/s^*}dx\right)\\
&=&O\left(n^{-1/2}n^{-(k+\gamma/s^*)s^*/(2-2\gamma)}\right)\\
&=&O\left(n^{-\frac{1}{2}\left(\frac{1+s^*k}{1-\gamma}\right)}\right).
\end{eqnarray*}
Analogously, one could show that 
\[
I_{n2}\leq O\left(n^{-\frac{1}{2}\left(\frac{1+s^*k}{1-\gamma}\right)}\right).
\]
To bound $I^\prime_{n1}$, note that for $x\leq x_{n1}$, it follows from  an 
analogous proof of (\ref{Formula1:ProofOfTheorem4}) that 
\[
F(x)
\leq\left(F^{s^*}(x_{n1})+s^*b_1(x-x_{n1})\right)^{1/s^*}
=\left(\lambda_\epsilon^{s^*}n^{-s^*/(2-2\gamma)}+s^*b_1(x-x_{n1})\right)^{1/s^*}.
\]
Analogously, it follows that for $x\leq x_{n1}$,
\[
U_n^o(x)
\leq
\left(U_n^{s^*}(x_{n1})+s^*b_1(x-x_{n1})\right)^{1/s^*}.
\]
Note that it follows from (\ref{Formula3:ProofOfTheorem4}) that 
\begin{eqnarray*}
U_n(x_{n1})
&=&U_n(x_{n1})-F(x_{n1})+F(x_{n1})\\
&\leq&\nu_\epsilon n^{-1/2}\left(F(x_{n1})(1-F(x_{n1}))\right)^\gamma+F(x_{n1})\\
&\leq&\nu_\epsilon n^{-1/2}F^{\gamma}(x_{n1})+F(x_{n1})\\
&=&(\nu_\epsilon\lambda_\epsilon^{\gamma}+\lambda_\epsilon)n^{-1/(2-2\gamma)}
\end{eqnarray*}
and hence for $x\leq x_{n1}$,
\begin{eqnarray*}
U_n^o(x)
&\leq&
\left((\nu_\epsilon\lambda_\epsilon^{\gamma}
    +\lambda_\epsilon)^{s^*}n^{-s^*/(2-2\gamma)}+s^*b_1(x-x_{n1})\right)^{1/s^*}.
\end{eqnarray*}
Thus, 
\begin{eqnarray*}
I^{\prime}_{n1}
&=&\int_{-\infty}^{x_{n1}}\left|\phi^\prime(x)\right|(F(x)+U_n^o(x))dx\\
&=&O\left(\int_{-\infty}^{x_{n1}}\left|\phi^\prime(x)\right|\left(n^{-s^*/(2-2\gamma)}+s^*b_1(x-x_{n1})\right)^{1/s^*}dx\right)\\
&=&O\left(\int_{-\infty}^{0}\left|\phi^\prime(x+x_{n1})\right|\left(n^{-s^*/(2-2\gamma)}+s^*b_1x\right)^{1/s^*}dx\right)\\
&=&O\left(\int_{-\infty}^{0}\left|x+x_{n1}\right|^{k-1}\left(n^{-s^*/(2-2\gamma)}+s^*b_1x\right)^{1/s^*}dx\right)\\
&=&O\left(n^{-1/(2-2\gamma)}\int_{-\infty}^{0}\left|x+x_{n1}\right|^{k-1}\left(1+s^*b_1x/n^{-s^*/(2-2\gamma)}\right)^{1/s^*}dx\right)\\
&=&O\left(n^{-1/(2-2\gamma)}n^{-s^*/(2-2\gamma)}\right.\\
& &\ \ \ \ \left. \cdot\int_{-\infty}^{0}\left|xn^{-s^*/(2-2\gamma)}+x_{n1}\right|^{k-1}\left(1+s^*b_1x\right)^{1/s^*}dx\right)\\
&=&O\left(n^{-1/(2-2\gamma)}n^{-s^*/(2-2\gamma)}\right.\\
& &\ \ \ \ \left. \cdot\int_{-\infty}^{0}\left|xn^{-s^*/(2-2\gamma)}+n^{-s^*/(2-2\gamma)}\right|^{k-1}\left(1+s^*b_1x\right)^{1/s^*}dx\right)\\
&=&O\left(n^{-1/(2-2\gamma)}n^{-ks^*/(2-2\gamma)}\int_{-\infty}^{0}|x|^{k-1}|x|^{1/s^*}dx\right)\\
&=&O\left(n^{-(ks^*+1)/(2-2\gamma)}\right).
\end{eqnarray*}
Analogously, one could show that 
\[
I^{\prime}_{n2}\leq O\left(n^{-(ks^*+1)/(2-2\gamma)}\right).
\]
Hence 
\begin{eqnarray*}
\sup_{G:L_n^o\leq G\leq U_n^o}\left|\int \phi d(G-F)\right|
&\leq& I_{n0}+I_{n1}+I_{n1}^\prime+I_{n2}+I_{n2}^\prime \\
&\leq& O(n^{-1/2})+O\left(n^{-(ks^*+1)/(2-2\gamma)}\right).
\end{eqnarray*}
\hfill $\Box$

\section{Appendix 1}
\label{Section:Appendix1} 
\textbf{Proof of the equivalence between Definition \ref{Definition1} and Definition \ref{Definition2}.}\\
Definition \ref{Definition1} implies Definition \ref{Definition2}:\\
For any $F\in\mathcal{P}_{s^*}$, Theorem \ref{Thm:CharacterizingThm} shows that 
$F$ is a continuous function on $\RR.$ By noticing that $J(F)\subset\RR$,  $J(F)\subset(\inf J(F),\infty)$ 
and $J(F)\subset(-\infty,\sup J(F))$, the convexity or concavity of $F^{s^*}$ or $(1-F)^{s^*}$ 
on $\RR$, $(\inf J(F),\infty)$ and $(-\infty,\sup J(F))$ imply the convexity or concavity of 
$F^{s^*}$ or $(1-F)^{s^*}$ on $J(F)$. 
Hence, Definition \ref{Definition1} implies Definition \ref{Definition2}.\\
Definition \ref{Definition2} implies Definition \ref{Definition1}:\\
Suppose $s^*<0$. By Definition \ref{Definition2}, for any $F\in \mathcal{P}_{s^*}$, $F^{s^*}$ 
and $(1-F)^{s^*}$ are convex on $J(F)$. 
Moreover, $F$ is continuous on $\RR$ and hence 
$J(F)=(a,b)$ where $a\equiv\inf J(F)$, $b\equiv\sup J(F)$.\\
To prove that $F^{s^*}$ is convex on $\RR$, by continuity of $F$ 
it suffices to prove that $F^{s^*}$ is mid-point convex:  that is, 
\begin{eqnarray}
F^{s^*}\left(\frac{x}{2}+\frac{y}{2}\right)\leq \frac{1}{2}F^{s^*}(x)+\frac{1}{2}F^{s^*}(y)
\label{Appendix:Formula1}
\end{eqnarray}
for any $x,y\in\RR$. 
Without loss of generality, we assume that $x<y$.\\
Note that if $a=-\infty$ and $b=\infty$, then there is nothing to prove. 
Without loss of generality, 
we assume that $a>-\infty$ and $b<\infty$.\\
Note that if $x\in(-\infty,a]$, then $F^{s^*}(x)=\infty$ and hence (\ref{Appendix:Formula1}) 
holds automatically. If $x\in(a,b)$ and $y\in(a,b)$, (\ref{Appendix:Formula1}) holds by the 
convexity of $F^{s^*}$ on $J(F)$. Moreover, by noticing the continuity of $F^{s^*}$ at $b$, 
(\ref{Appendix:Formula1}) holds for any $x\in(a,b)$ and $y\in(a,b]$. Since $F^{s^*}(y)=F^{s^*}(b)=1$ 
for $y\geq b$, (\ref{Appendix:Formula1}) holds for any $x\in(a,b)$ and $y\in[b,\infty)$. 
If $x,y\in[b,\infty)$, (\ref{Appendix:Formula1}) holds automatically since $F^{s^*}(x)=F^{s^*}(y)=1$. \\
The proof of the convexity of $(1-F)^{s^*}$ on $\RR$ is similar and hence is omitted.\\
For the cases that $s^*\geq 0$, the proof is similar and hence is omitted. 
\hfill $\Box$
\medskip

\par\noindent
\textbf{Proof of Theorem \ref{Thm:CharacterizingThm} ($0 \leq s^*\leq 1$):}\\
Recall that $a\equiv \inf J(F)$ and $b\equiv\sup J(F)$. Suppose $1\geq s^*>0$.\\
(i) implies (ii): \\
Suppose $F\in\mathcal{P}_{s^*}$. 
To prove that $F$ is continuous on $\RR$, we first note that $x\mapsto F^{s^*}(x)$ 
and $x\mapsto \left(1-F(x)\right)^{s^*}(x)$ are concave functions on $(a,\infty)$ and $(-\infty,b)$ respectively. 
By Theorem 10.1 (page 82) in \cite{10.2307/j.ctt14bs1ff}, $F^{s^*}$ and $\left(1-F(x)\right)^{s^*}$ 
are continuous on any open convex sets in their effective domains.
In particular, $F^{s^*}$ and $\left(1-F\right)^{s^*}$ are continuous on $(a,\infty)$ and $(-\infty,b)$ respectively. 
This yields that $F$ is continuous on $(a,\infty)$ and $(-\infty,b)$, or equivalently, on 
$(a,\infty)\cup(-\infty,b)=(-\infty,\infty)$ since $F$ is non-degenerate.\\
To prove that $F$ is differentiable on $J(F)$, note that $J(F)=(a,b)$ since $F$ is continuous on $\RR$.
By Theorem 23.1 (page 213) in \cite{10.2307/j.ctt14bs1ff}, for any $x\in J(F)$, the concavity of 
$F^{s^*}$ on $J(F)$ implies the existence of $\left(F^{s^*}\right)^{\prime}_+(x)$ and $\left(F^{s^*}\right)^{\prime}_-(x)$.
Moreover, $\left(F^{s^*}\right)^{\prime}_-(x)\geq \left(F^{s^*}\right)^{\prime}_+(x)$ by 
Theorem 24.1 (page 227) in \cite{10.2307/j.ctt14bs1ff}. Since $F=\left(F^{s^*}\right)^{1/s^*}$ 
on $J(F)$, the chain rule guarantees the existence of $F^\prime_\pm(x)$ and 
\[
F^\prime_\pm(x)=\frac{1}{s^*}\left(F^{s^*}\right)^{1/s^*-1}(x\pm)\left(F^{s^*}\right)^{\prime}_\pm(x).
\]
Since $F$ is continuous on $J(F)$, then
\[
F_\pm^\prime(x)=\frac{1}{s^*}\left(F^{s^*}\right)^{1/s^*-1}(x)\left(F^{s^*}\right)^{\prime}_\pm(x).
\]
Hence $F^\prime_-(x)\geq F^\prime_+(x)$ by $\left(F^{s^*}\right)^{\prime}_-(x)\geq \left(F^{s^*}\right)^{\prime}_+(x)$. \\
Similarly, one can prove $F^\prime_-(x)\leq F^\prime_+(x)$ by the concavity of $\left(1-F\right)^{s^*}$ on $J(F)$.\\
Thus $F^\prime_-(x)=F^\prime_+(x)=F^\prime(x)$ for any $x\in J(F)$, or equivalently, $F$ is differentiable on $J(F)$. 
The derivative of $F$ is denoted by $f$, i.e. $f\equiv F^\prime$.\\
To prove (\ref{Formula:TailBounds1}), note that the concavity of $x\mapsto F^{s^*}(x)$ on $J(F)$ implies that, 
for any $x,y\in J(F)$,
\[
F^{s^*}(y)-F^{s^*}(x)\leq (y-x) \left(F^{s^*}\right)^{\prime}(x)=(y-x)s^*F^{s^*-1}(x)f(x),
\]
or, with $x_+=\max\{x,0\}$,
\[
\frac{F^{s^*}(y)}{F^{s^*}(x)}\leq \left(1+s^*\frac{f(x)}{F(x)}(y-x)\right)_+.
\]
Hence
\[
\frac{F(y)}{F(x)}\leq \left(1+s^*\frac{f(x)}{F(x)}(y-x)\right)_+^{1/s^*},
\]
or,  equivalently,
\[
F(y)\leq F(x)\left(1+s^*\frac{f(x)}{F(x)}(y-x)\right)_+^{1/s^*}.
\]
Analogously, the convexity of $\left(1-F(x)\right)^{s^*}$ on $J(F)$ implies that 
for any $x,y\in J(F)$
\[
\left(1-F(y)\right)^{s^*}-\left(1-F(x)\right)^{s^*}\leq -(y-x)s^*\left(1-F(x)\right)^{s^*-1}f(x),
\]
or, equivalently,
\[
\left(\frac{1-F(y)}{1-F(x)}\right)^{s^*}\leq \left(1-s^*\frac{f(x)}{1-F(x)}(y-x)\right)_+,
\]
which yields
\[
F(y)\geq 1- \left(1-F(x)\right)\left(1-s^*\frac{f(x)}{1-F(x)}(y-x)\right)^{1/s^*}_+.
\]
The proof of (\ref{Formula:TailBounds1}) is complete.\\
(ii)  implies  (iii):\\
Applying (\ref{Formula:TailBounds1}) yields that for any $x,y\in J(F)$ with $x<y$, 
\[
\frac{F^{s^*}(x)}{F^{s^*}(y)}\leq 1+s^*\frac{f(y)}{F(y)}(x-y),
\]
and
\[
\frac{F^{s^*}(y)}{F^{s^*}(x)}\leq 1+s^*\frac{f(x)}{F(x)}(y-x),
\]
or, equivalently, 
\[
F^{s^*}(x)\leq F^{s^*}(y)+s^*\frac{f(y)}{F^{1-s^*}(y)}(x-y),
\]
and
\[
F^{s^*}(y)\leq F^{s^*}(x)+s^*\frac{f(x)}{F^{1-s^*}(x)}(y-x).
\]
By defining $h\equiv f/F^{1-s^*}$ on $J(F)$, it follows that 
\[
F^{s^*}(x)\leq F^{s^*}(y)+s^*h(y)(x-y),
\]
and
\[
F^{s^*}(y)\leq F^{s^*}(x)+s^*h(x)(y-x).
\]
After summing up the last two inequalities, it follows that 
\[
F^{s^*}(x)+F^{s^*}(y)\leq F^{s^*}(y)+s^*h(y)(x-y)+F^{s^*}(x)+s^*h(x)(y-x),
\]
or, equivalently,
\[
0\leq s^*\left(h(x)-h(y)\right)(y-x).
\]
Hence $h(x)\geq h(y)$, or equivalently, $h(\cdot)$ is a monotonically non-increasing function on $J(F)$.\\
The proof of the monotonicity of
$\tilde{h}\equiv f/(1-F)^{1-s^*}$ is similar and hence is omitted.\\
(iii)  implies  (iv):\\
If (iii) holds, it immediately follows that $f>0$ on $J(F)=(a,b)$. 
If not, suppose that $f(x_0)=0$ for some $x_0\in J(F)$.
It follows that $h(x_0)=f(x_0)/F^{1-s^*}(x_0)=0$.
Since $h$ is monotonically non-increasing on $J(F)$,
$h(x)=0$ for all $x\in [x_0,b)$,
or, equivalently,
$f=0$ on $[x_0,b)$.
Similarly, the non-decreasing monotonicity of $x\mapsto\tilde{h}(x)$ 
on $J(F)$ implies that $f=0$ on $(a,x_0]$.
Then $f=0$ on $J(F)$, which violates the continuity assumption in (iii) 
and hence $f>0$ on $J(F)$.\\
To prove $f$ is bounded on $J(F)$, note that the monotonicities of $h$ 
and $\tilde{h}$ imply that for any $x,x_0\in J(F)$,
\begin{eqnarray*}
f(x)=\left\{
\begin{tabular}{ll }
$F^{1-s^*}h(x)\leq h(x)\leq h(x_0)$, &if $x\geq x_0$,\\
$(1-F(x))^{1-s^*}\tilde{h}(x)\leq \tilde{h}(x)\leq \tilde{h}(x_0)$, & if $x\leq x_0$.
\end{tabular}\right.
\end{eqnarray*}
Hence $f(x)\leq \max\{h(x_0),\tilde{h}(x_0)\}$ for any $x,x_0\in J(F)$.\\
To prove that $f$ is differentiable on $J(F)$ almost every, we first prove that 
$f$ is Lipschitz continuous on $(c,d)$ for any $c,d\in J(F)$ with $c<d$.\\
By noticing the non-increasing monotonicity of $h$ on $J(F)$, the following 
arguments yield an upper bound of $\left(f(y)-f(x)\right)/(y-x)$ for $x,y\in(c,d)$:
\begin{eqnarray*}
\frac{f(y)-f(x)}{y-x}
&=&\frac{F^{1-s^*}(y)h(y)-F^{1-s^*}(x)h(x)}{y-x}\\
&=&h(y)\frac{F^{1-s^*}(y)-F^{1-s^*}(x)}{y-x}+F^{1-s^*}(x)\frac{h(y)-h(x)}{y-x}\\
&\leq&h(y)\frac{F^{1-s^*}(y)-F^{1-s^*}(x)}{y-x}\\
&=&h(y)(1-s^*)f(z)F^{-s^*}(z),
\end{eqnarray*}
where the last equality follows from the mean value theorem and $z$ is between $x$ and $y$.\\
Since $-s^*<0$, it follows that $F^{-s^*}(z)<F^{-s^*}(c)$ and hence 
\[
\frac{f(y)-f(x)}{y-x}\leq(1-s^*)f(z)h(y)F^{-s^*}(z)\leq (1-s^*)\max\{h(x_0),\tilde{h}(x_0)\}h(c)F^{-s^*}(c),
\]
for $x,y\in (c,d)$.\\
Similar arguments imply that 
\begin{eqnarray*}
\frac{f(y)-f(x)}{y-x}
&=&\frac{\bar{F}^{1-s^*}(y)\tilde{h}(y)-\bar{F}^{1-s^*}(x)\tilde{h}(x)}{y-x}\\
&=&\tilde{h}(y)\frac{\bar{F}^{1-s^*}(y)-\bar{F}^{1-s^*}(x)}{y-x}+\bar{F}^{1-s^*}(x)\frac{\tilde{h}(y)-\tilde{h}(x)}{y-x}\\
&\geq&\tilde{h}(y)\frac{\bar{F}^{1-s^*}(y)-\bar{F}^{1-s^*}(x)}{y-x}\\
&=&-\tilde{h}(y)(1-s^*)\bar{F}^{-s^*}(z)f(z)\\
&\geq&-(1-s^*)\max\{h(x_0),\tilde{h}(x_0)\}\tilde{h}(d)\bar{F}^{-s^*}(d).
\end{eqnarray*}
Hence 
\[
\left|\frac{f(y)-f(x)}{y-x}\right|\leq (1-s^*)\max\{h(x_0),\tilde{h}(x_0)\}\max\{h(c)F^{-s^*}(c),\tilde{h}(d)\bar{F}^{-s^*}(d)\}.
\]
The last display shows that $f$ is Lipschitz continuous on $(c,d)$.\\
By Proposition 4.1(iii) of \cite{Shorack2017}, page 82, $f$ is absolutely continuous on $(c,d)$, 
and hence $f$ is differentiable on $(c,d)$ almost everywhere.\\
Since $(c,d)$ is an arbitrary interval in $(a,b)$, the differentiability of $f$ on $(c,d)$ 
implies the differentiability of $f$ on $(a,b)$ and hence $f$ is differentiable on $(a,b)$ 
with $f^\prime=F^{\prime\prime}$ almost everywhere.\\
Since $f$ is differentiable almost everywhere, the non-increasing monotonicity of $h$ on $J(F)$ implies that
\[
h^\prime(x)\leq 0 \text{ almost everywhere on $J(F)$,}
\]
or, equivalently,
\[
\log(h)^\prime(x)\leq 0 \text{ almost everywhere on $J(F)$.}
\]
Straight-forward calculation yields that the last display is equivalent to 
\[
\frac{f^\prime}{f}-(1-s^*)\frac{f}{F}\leq 0 \text{ almost everywhere on $J(F)$,}
\]
or,
\[
f^{\prime} \leq(1-s^*)\frac{f^2}{F}\text{ almost everywhere on $J(F)$,}
\]
which is the right hand side of (\ref{Formula:fPrimeBounds}).\\
Similarly, the non-decreasing monotonicity of $\tilde{h}$ implies the left hand side of (\ref{Formula:fPrimeBounds}).\\
(iv)  implies  (i):\\
Since $F$ is continuous on $\RR$, it suffices to prove that $F^{s^*}$ is convex on $J(F)$ by Definition \ref{Definition2}.
Since we assume that $F$ is differentiable on $J(F)$ with derivative $f=F^\prime$, 
the concavity of $F^{s^*}$ on $J(F)$ can be proved by the non-increasing monotonicity 
of the first derivative of $F^{s^*}$ on $J(F)$. Since $f$ is differentiable almost everywhere 
on $J(F)$, the non-increasing monotonicity of $\left(F^{s^*}\right)^{\prime}$ on $J(F)$ 
can be proved by the non-positivity of $\left(F^{s^*}\right)^{\prime\prime}$ on $J(F)$ almost everywhere, which follows from
\begin{eqnarray*}
\left(F^{s^*}\right)^{\prime\prime}(x)=s^*F^{s^*-1}(x)\left(-(1-s^*)\frac{f^2(x)}{F(x)}+f^\prime(x)\right)\leq 0,
\end{eqnarray*}
where $f=F^\prime, f^\prime=F^{\prime\prime}$. 
The last inequality follows from the right hand side of (\ref{Formula:fPrimeBounds}).\\
Similarly, the concavity of $\left(1-F(x)\right)^{s^*}$, or $\bar{F}^{s^*}$, on $J(F)$ can be proved by the following arguments:
\[
\left(\bar{F}^{s^*}\right)^{\prime\prime}(x)=s^*\bar{F}^{s^*-1}(x)\left(-(1-s^*)\frac{f^2(x)}{\bar{F}(x)}-f^\prime(x)\right)\leq 0,
\]
where the last inequality follows from the left part of (\ref{Formula:fPrimeBounds}).
\hfill $\Box$

\section{Appendix 2:   The Figures again, step by step}
\label{Section:Appendix2} 
\subsection{Bi-$s^*$-concave bands based on KS bands:  $t_1$ sample of size 100 }
\begin{figure}[!htb]
\center{\includegraphics[width=\textwidth]{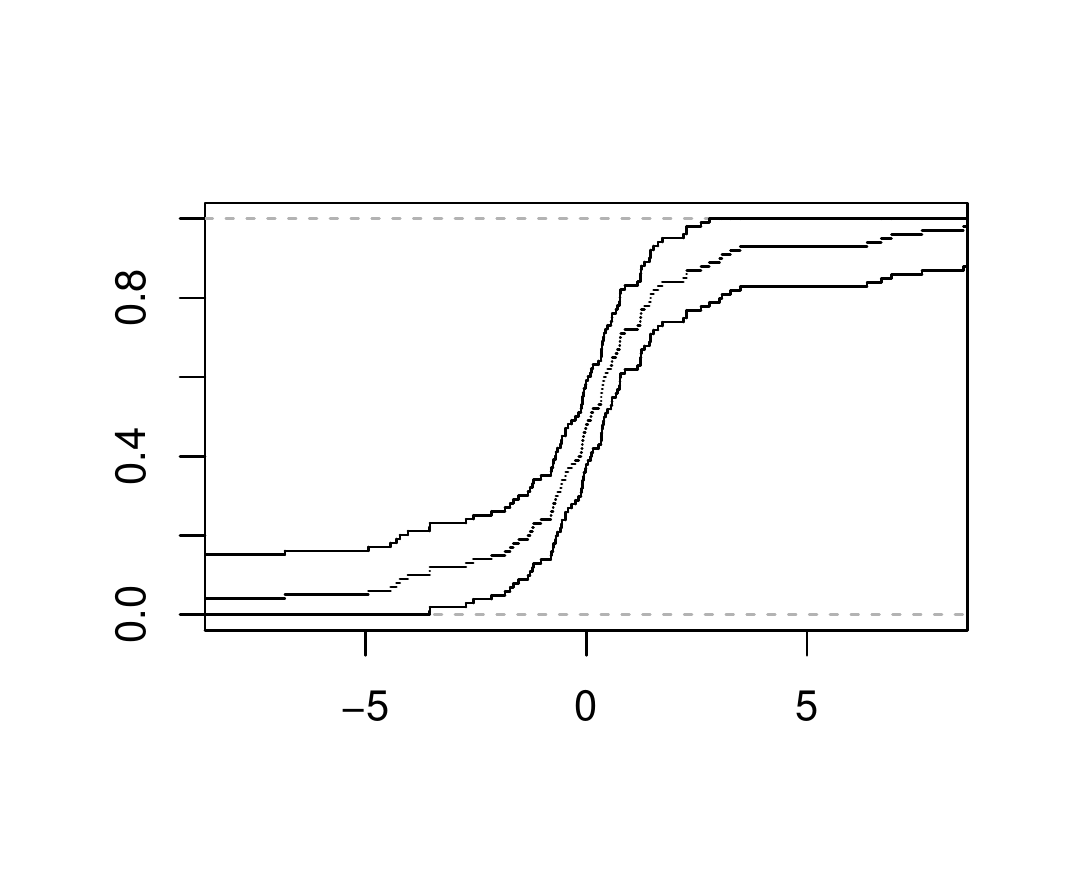}}
\caption{
KS confidence bands based on a sample of size $100$ generated from Student-$t$ distribution with one degree of freedom.
The step function (black) in the middle is the empirical distribution function. 
The upper and lower black lines correspond to the confidence bands.
}
\label{plot:KS1a-apdx}
\end{figure}

\begin{figure}[!htb]
\center{\includegraphics[width=\textwidth]{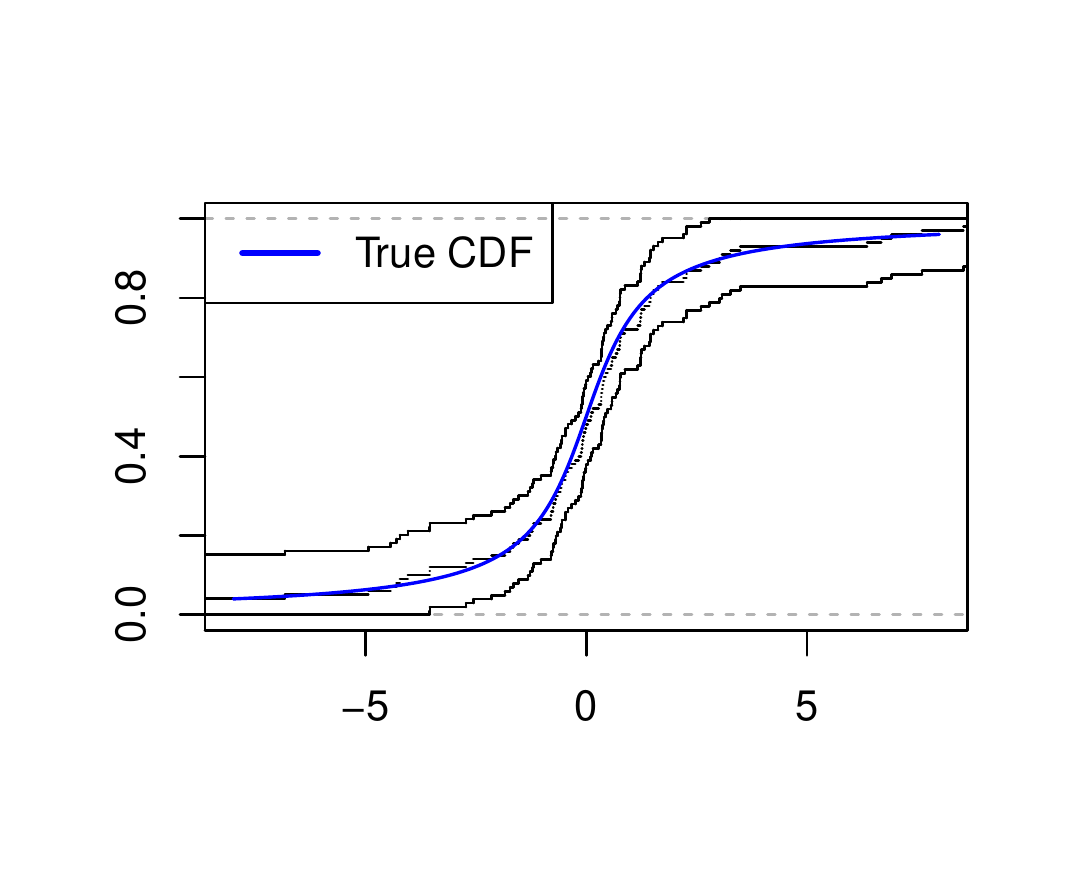}}
\caption{
KS confidence bands  based on a sample of size $100$ generated from Student-$t$ distribution  with one degree of freedom.
The step function (black) in the middle is the empirical distribution function. 
The upper and lower black lines correspond to the confidence bands. 
The blue  line corresponds to the true distribution function. 
}
\label{plot:KS1b-apdx}
\end{figure}

\begin{figure}[!htb]
\center{\includegraphics[width=\textwidth]{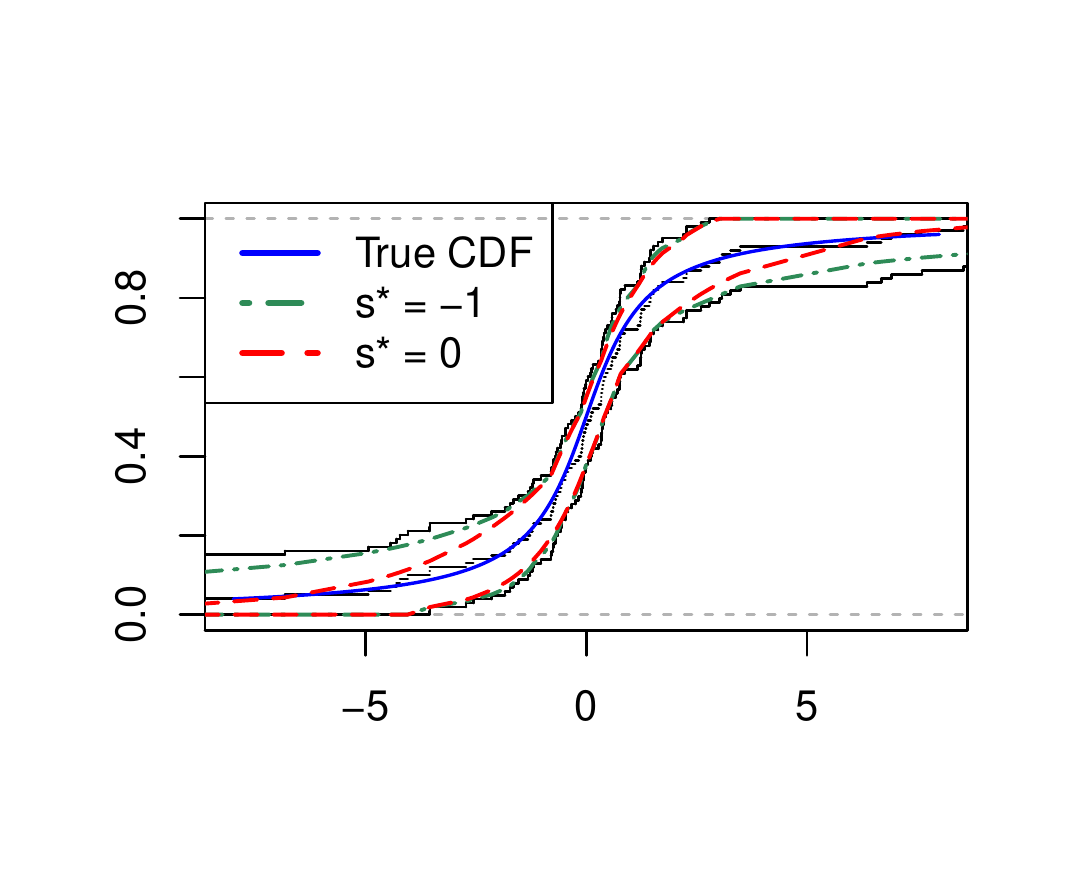}}
\caption{
Confidence bands for bi-$s^*$-concave distribution functions from KS bands 
based on  a sample of size $100$ generated from Student-$t$ distribution with $1$ degree of freedom. 
The blue dashed line corresponds to the true distribution function. 
The two black lines give the KS band and lines in other colors are 
refined confidence bands under the bi-$s^*$-concave assumption. In this case, $s^*_0=-1$.
The step function (black) in the middle is the empirical distribution function.
}
\label{plot:KS1c-apdx} 
\end{figure}
\FloatBarrier
\subsection{Bi-$s^*$-concave bands based on WKS bands:  $t_1$ sample of size 100 }
\begin{figure}[!htb]
\center{\includegraphics[width=\textwidth]{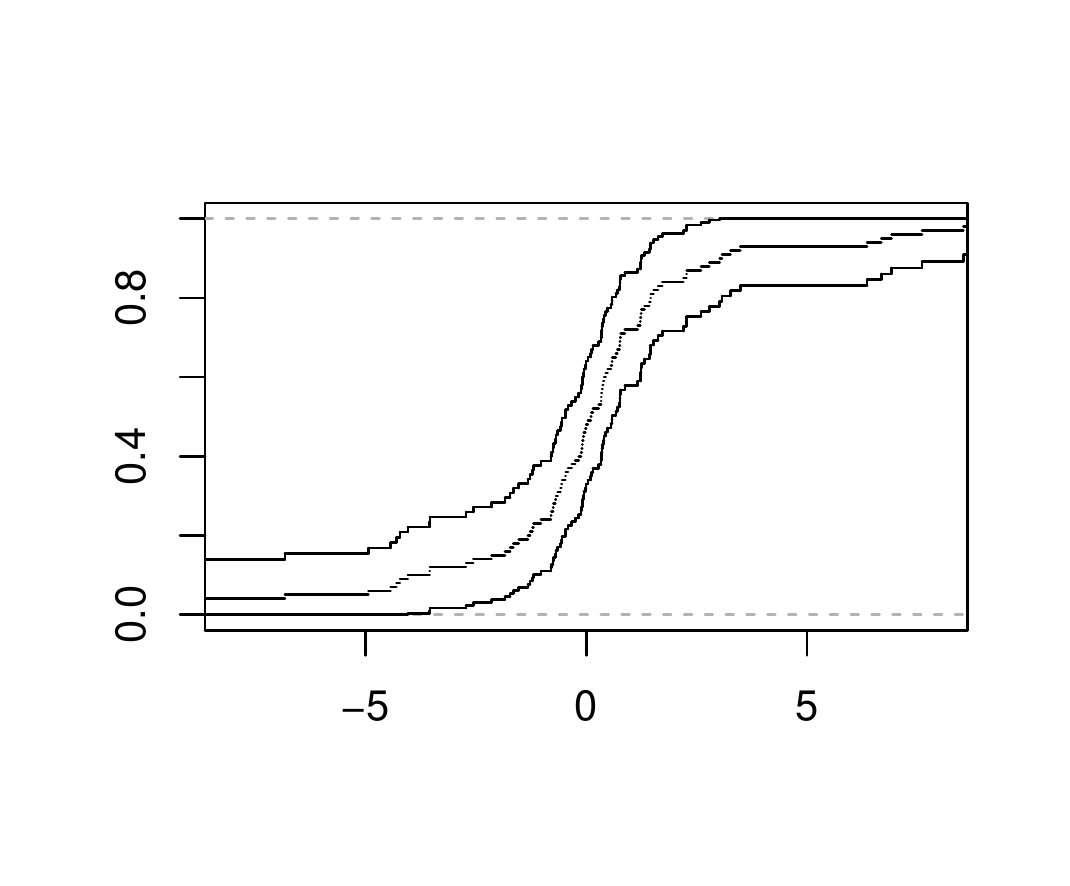}}
\caption{
WKS confidence bands based on a sample of size $100$ generated from Student-$t$ distribution with  $1$ degree of freedom.
The step function (black) in the middle is the empirical distribution function. 
The upper and lower black lines correspond to the confidence bands.
}
\label{plot:WKS2a-apdx} 
\end{figure}

\begin{figure}[!htb]
\center{\includegraphics[width=\textwidth]{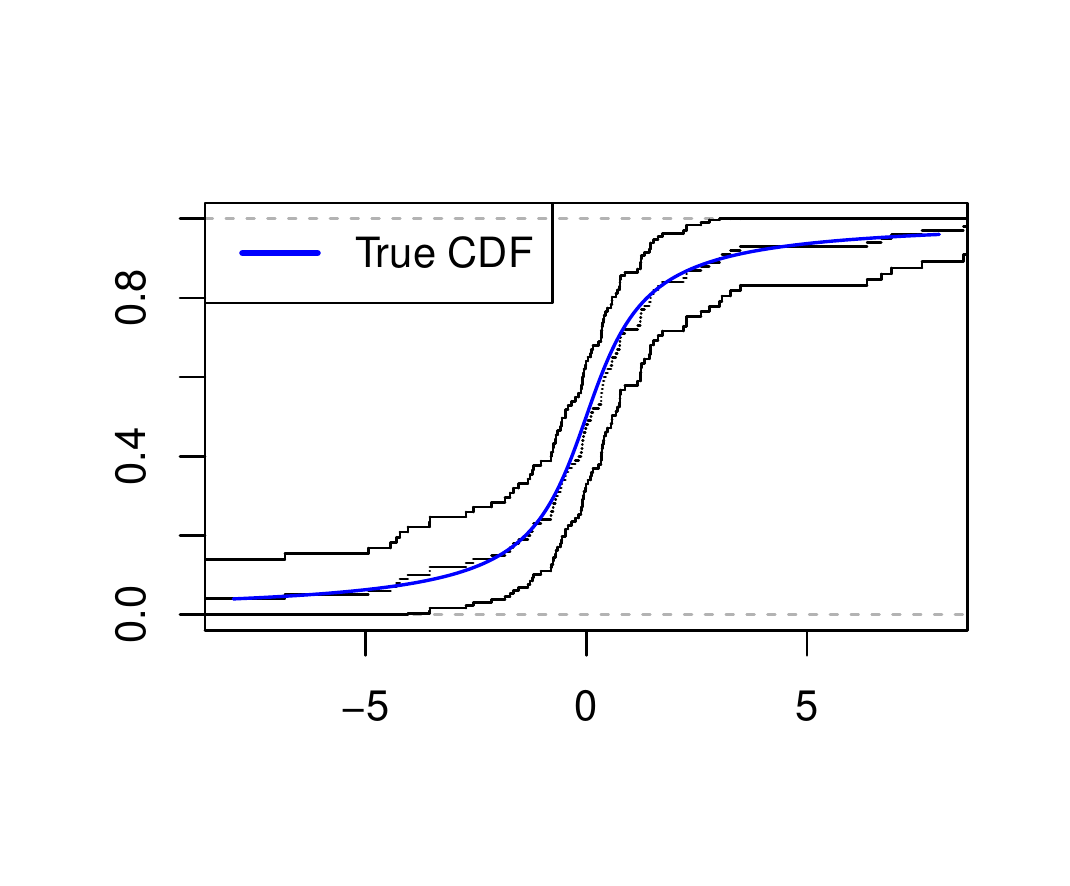}}
\caption{
WKS confidence bands  based on a sample of size $100$ generated from Student-$t$ distribution  with one degree of freedom.
The step function (black) in the middle is the empirical distribution function. 
The upper and lower black lines correspond to the confidence bands. The blue  line corresponds to the true distribution function. 
}
\label{plot:WKS2b-apdx} 
\end{figure}

\begin{figure}[!htb]
\center{\includegraphics[width=\textwidth]{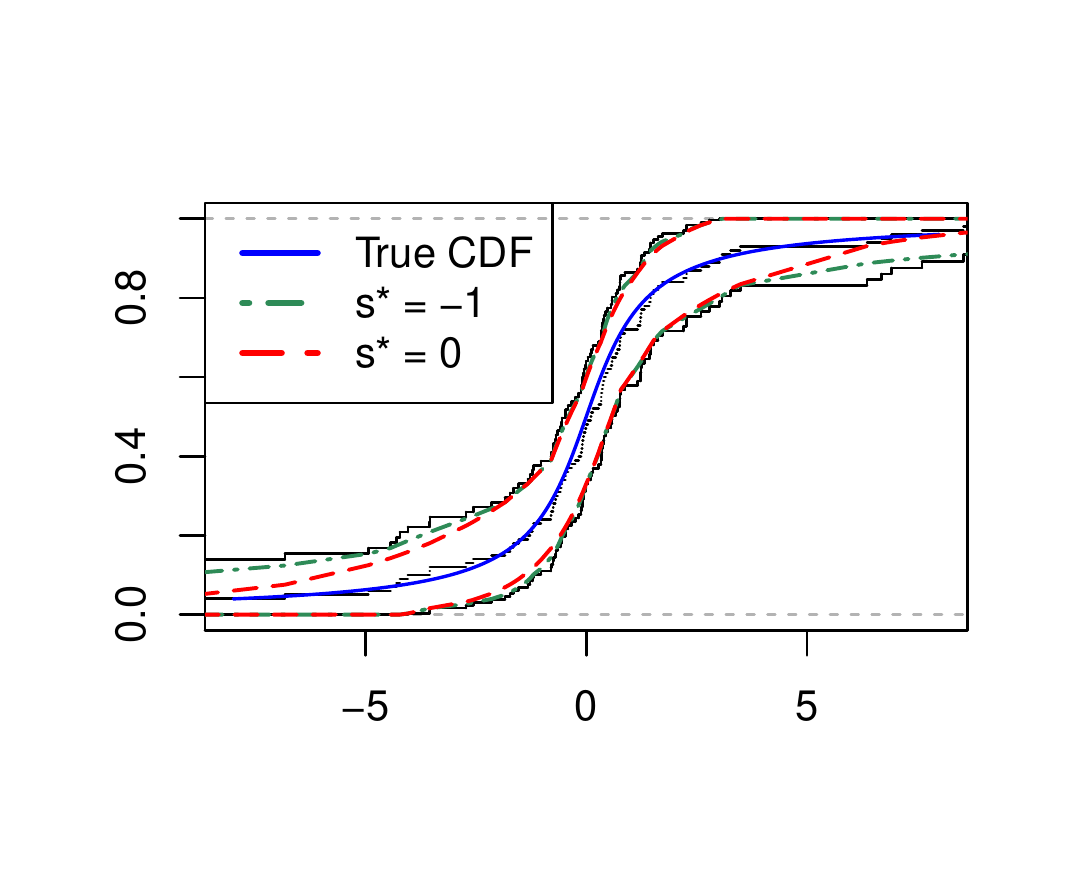}}
\caption{
Confidence bands for bi-$s^*$-concave distribution functions from WKS bands based on  a sample of 
size $100$ generated from Student-$t$ distribution with $1$ degree of freedom. 
The blue dashed line corresponds to the true distribution function. 
The two black lines give the WKS band and lines in other colors are 
refined confidence bands under the bi-$s^*$-concave assumption. In this case, $s^*_0=-1$.
The step function (black) in the middle is the empirical distribution function.
}
\label{plot:WKS2c-apdx} 
\end{figure}
\FloatBarrier
\subsection{Bi-$s^*$-concave bands based on KS bands:  $t_1$ sample of size 1000 }

\begin{figure}[!htb]
\center{\includegraphics[width=\textwidth]{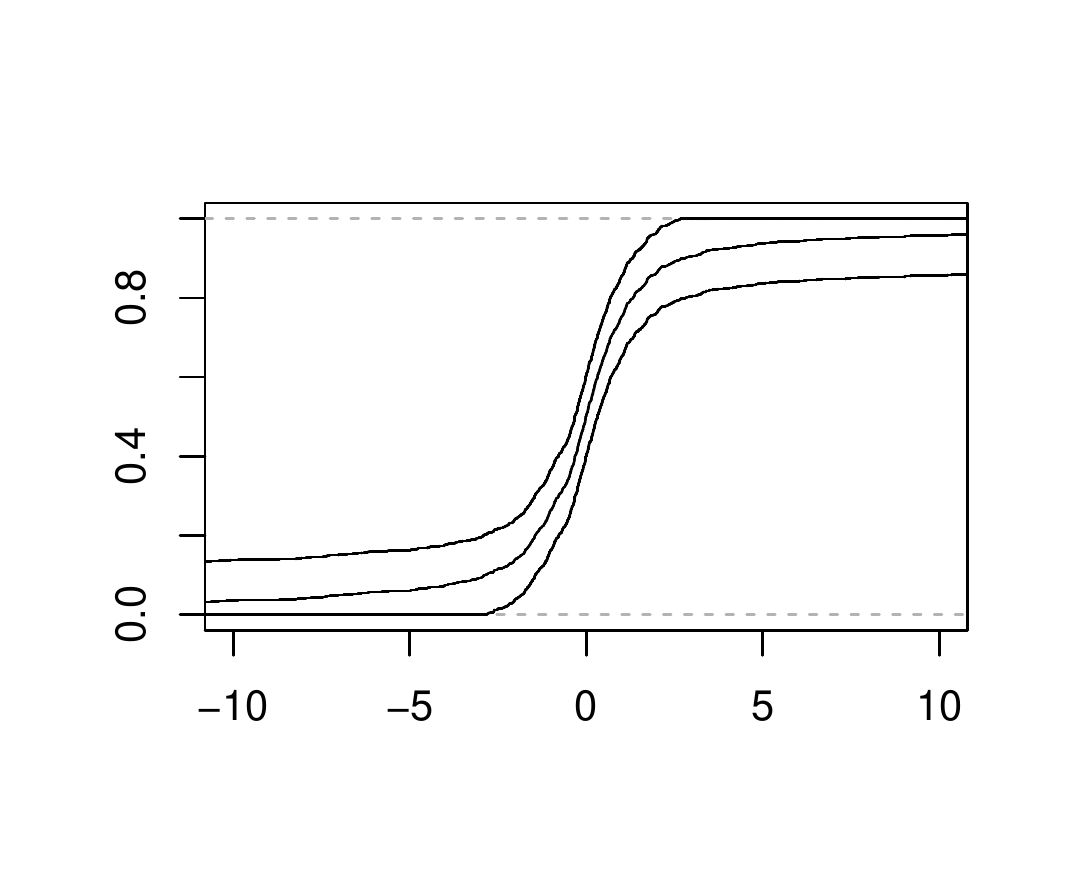}}
\caption{
KS confidence bands based on a sample of size $1000$ generated from Student-$t$ distribution with degree of freedom one.
The step function (black) in the middle is the empirical distribution function. 
The upper and lower black lines correspond to the confidence bands.
}
\label{plot:KS3a-apdx} 
\end{figure}

\begin{figure}[!htb]
\center{\includegraphics[width=\textwidth]{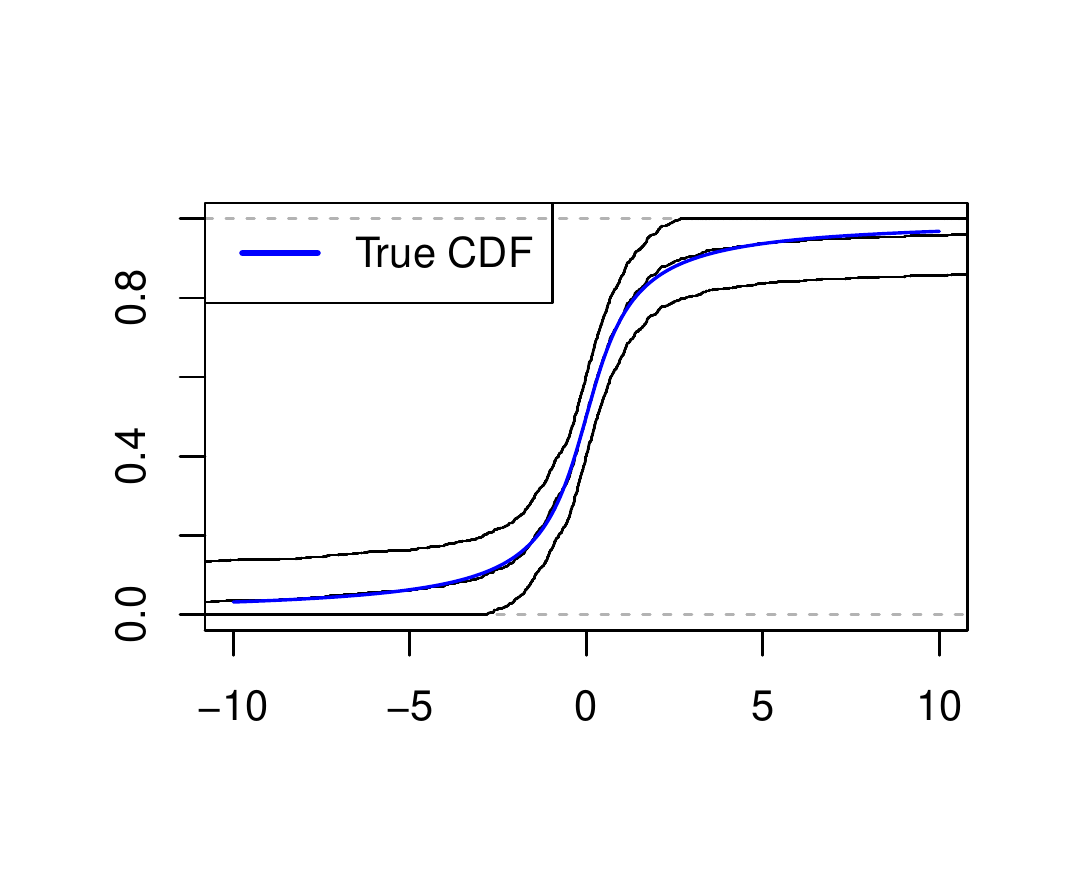}}
\caption{ 
KS confidence bands  based on a sample of size $1000$ generated from 
Student-$t$ distribution  with degree of freedom one.
The step function (black) in the middle is the empirical distribution function. 
The upper and lower black lines correspond to the confidence bands. 
The blue  line corresponds to the true distribution function. 
}
\label{plot:KS3b-apdx}
\end{figure}

\begin{figure}[!htb]
\center{\includegraphics[width=\textwidth]{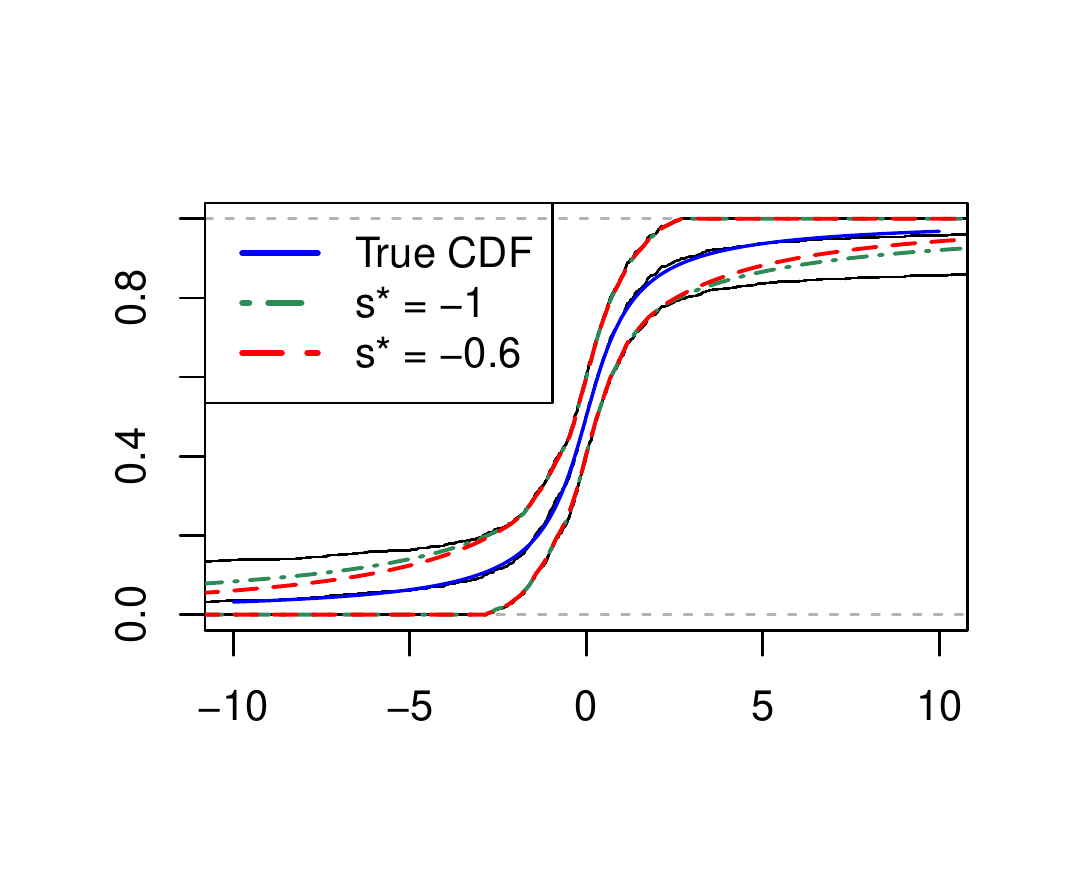}}
\caption{
Confidence Bands for bi-$s^*$-concave distribution functions from KS bands based on  a sample of size 
$1000$ generated from Student-$t$ distribution with $1$ degree of freedom. 
The blue dashed line corresponds to the true distribution function. 
The two black lines give the KS band and lines in other colors are 
refined confidence bands under the bi-$s^*$-concave assumption. In this case, $s^*_0=-1$.
The step function (black) in the middle is the empirical distribution function.
}
\label{plot:KS3c-apdx} 
\end{figure}

\FloatBarrier

\subsection{Bi-$s^*$-concave bands based on WKS bands:  $t_1$ sample of size 1000}

\begin{figure}[!htb]
\center{\includegraphics[width=\textwidth]{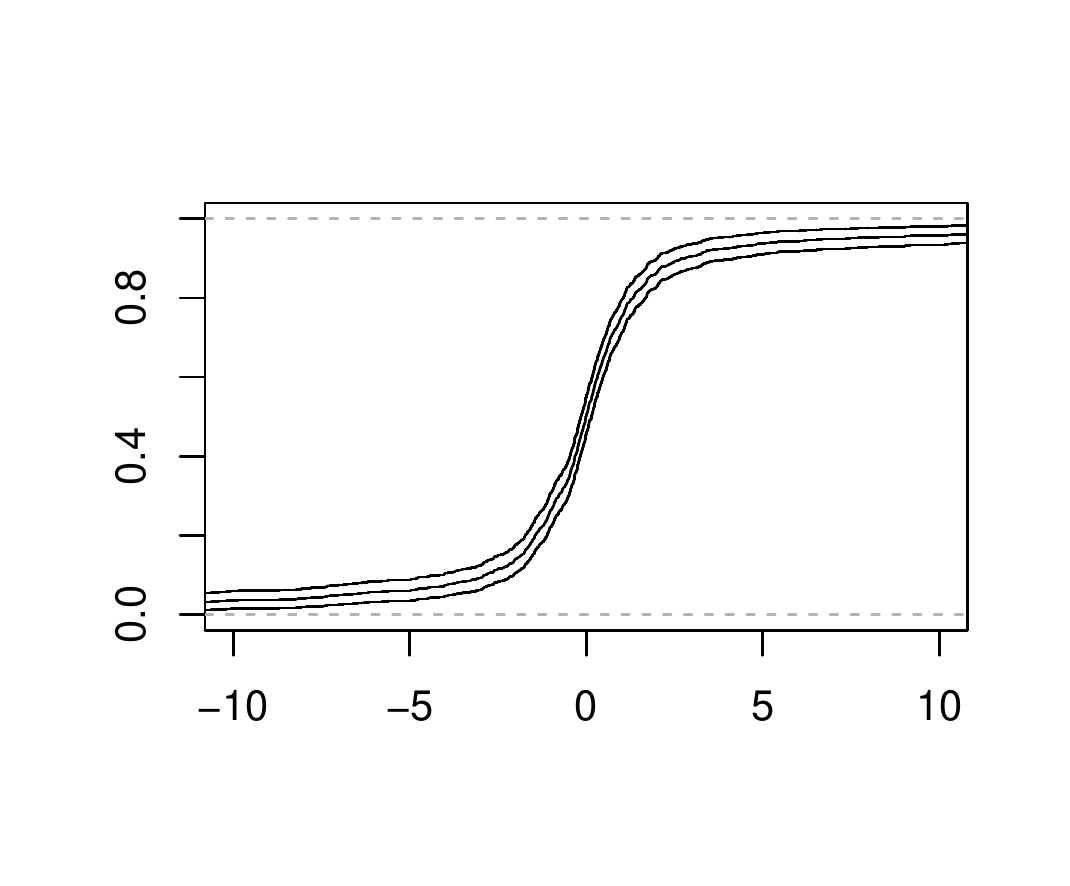}}
\caption{
WKS confidence bands based on a sample of size $1000$ generated from 
Student-$t$ distribution with one degree of freedom.
The step function (black) in the middle is the empirical distribution function. 
The upper and lower black lines correspond to the confidence bands.
}
\label{plot:WKS4a-apdx} 
\end{figure}

\begin{figure}[!htb]
\center{\includegraphics[width=\textwidth]{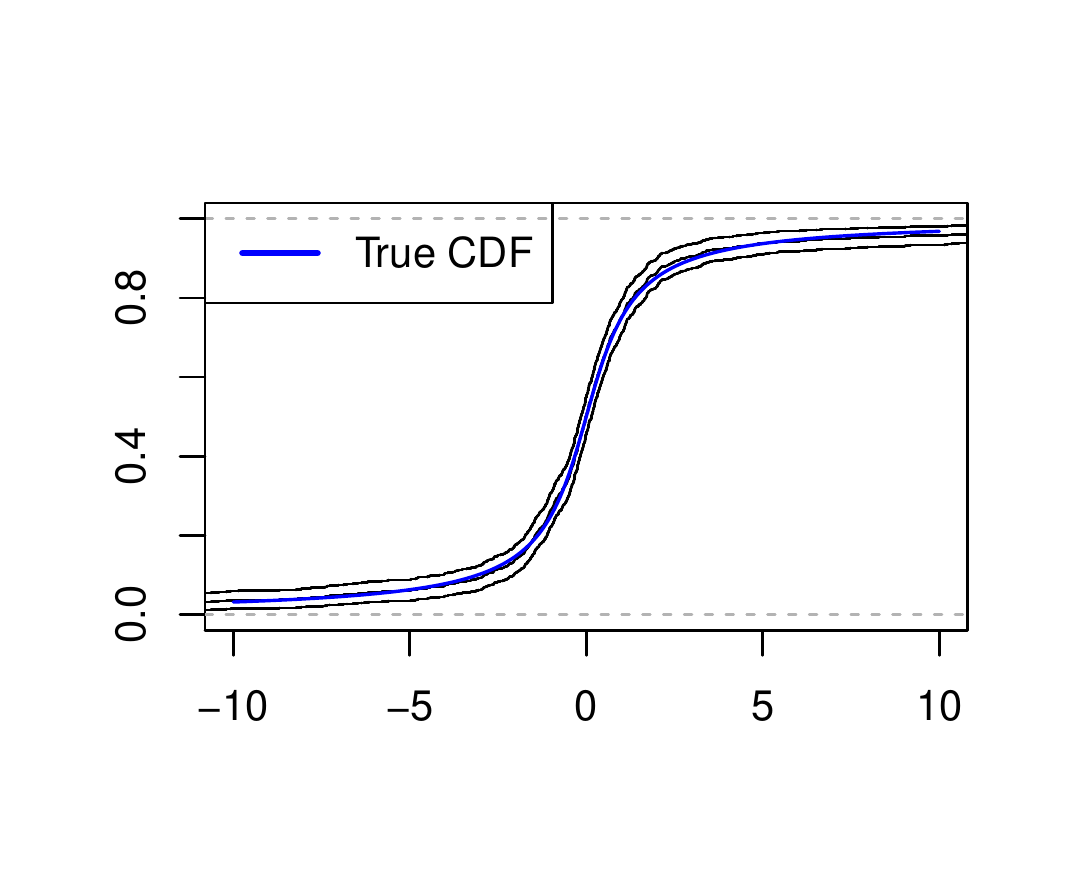}}
\caption{
KS confidence bands  based on a sample of size $1000$ generated from 
Student-$t$ distribution  with one degree of freedom.
The step function (black) in the middle is the empirical distribution function. 
The upper and lower black lines correspond to the confidence bands. 
The blue  line corresponds to the true distribution function.
}
\label{plot:WKS4b-apdx} 
\end{figure}

\begin{figure}[!htb]
\center{\includegraphics[width=\textwidth]{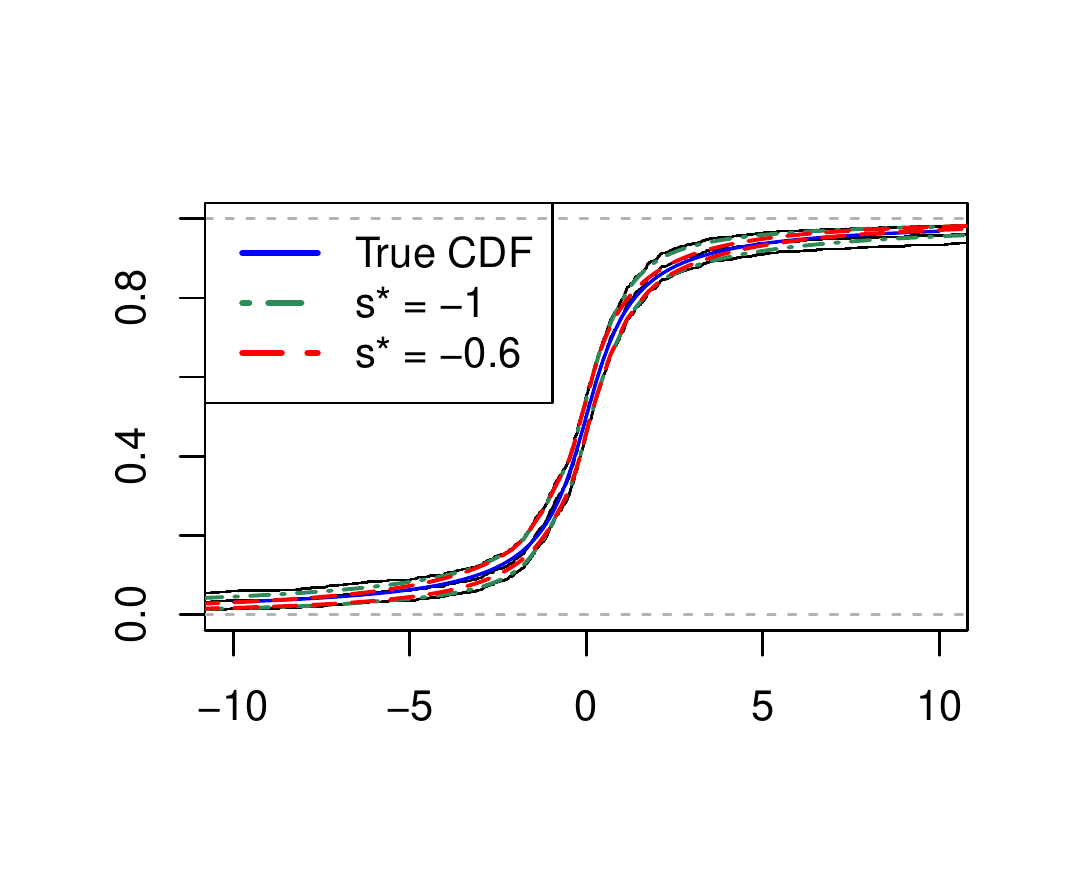}}
\caption{
Confidence bands for bi-$s^*$-concave distribution functions from WKS bands based on  a sample of size 
$1000$ generated from Student-$t$ distribution with $1$ degree of freedom $1$. 
The blue dashed line corresponds to the true distribution function. 
The two black lines give the KS band and lines in other colors are 
refined confidence bands under the bi-$s^*$-concave assumption. In this case, $s^*_0=-1$.
The step function (black) in the middle is the empirical distribution function.
}
\label{plot:WKS4c-apdx} 
\end{figure}
\FloatBarrier
\subsection{Bi-$s^*$-concave bands based on KS bands:  application to the CEO salary data}

\begin{figure}[!htb]
\center{\includegraphics[width=\textwidth]{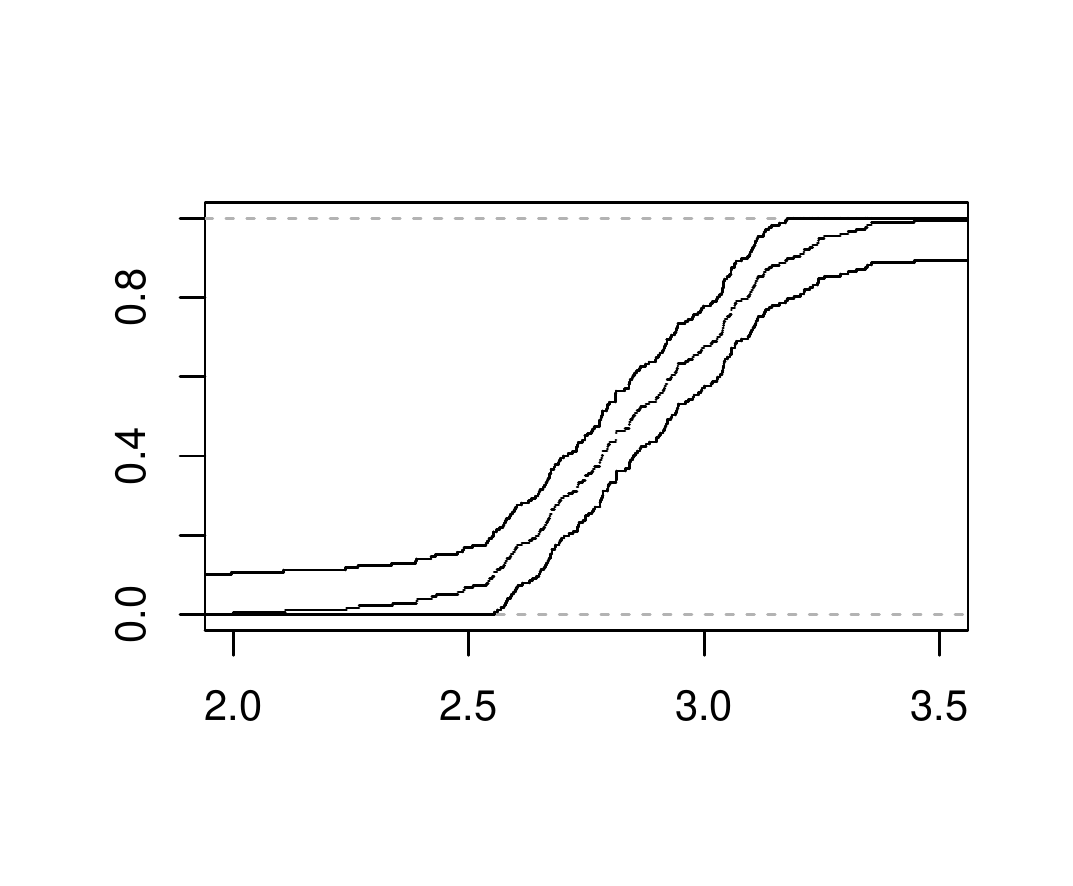}}
\caption{
KS confidence bands for the CEO salary data.
The step function in the middle is the empirical distribution function.
The upper and lower lines give the KS bands 
}
\label{plot:KS5a-apdx} 
\end{figure}

\begin{figure}[!htb]
\center{\includegraphics[width=\textwidth]{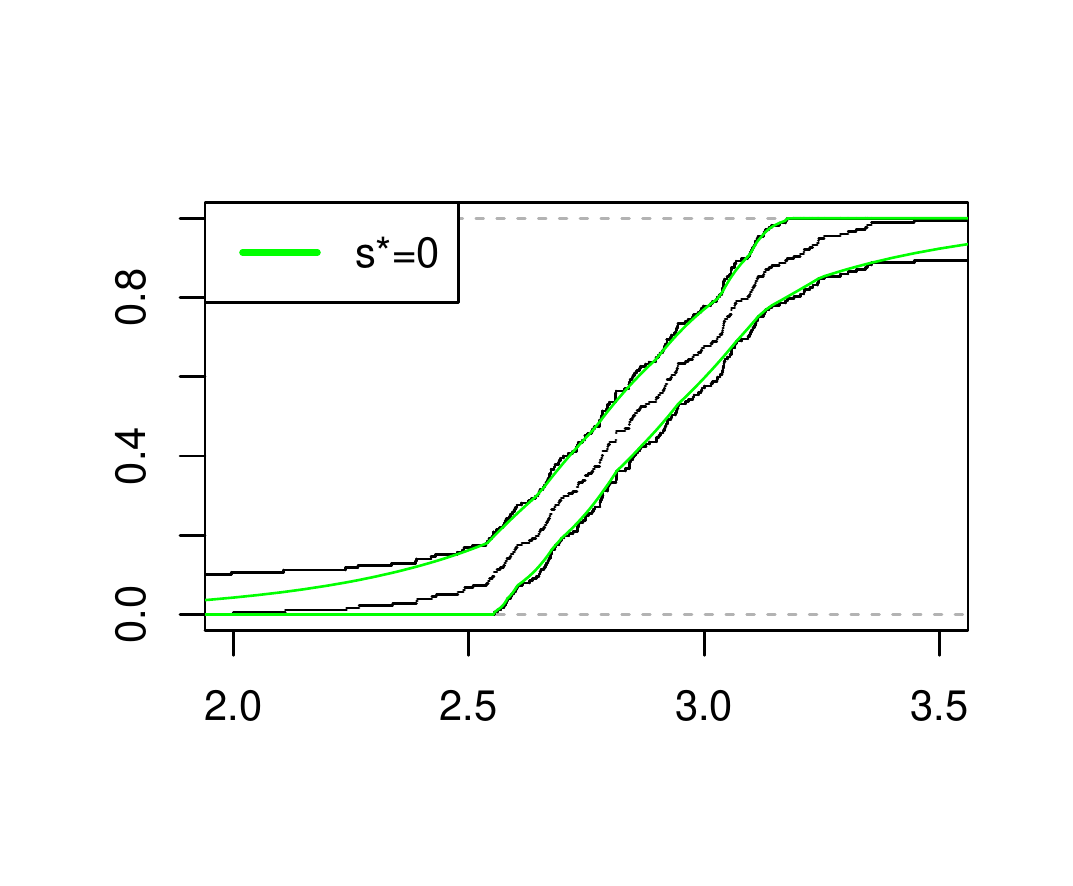}} 
\caption{
KS confidence bands for the CEO salary data.
The step function in the middle is the empirical distribution function.
The upper and lower lines give the KS bands and the green lines  are  the bi-$s^*$-concave bands for $s^*=0$.
}
\label{plot:KS5b-apdx}
\end{figure}

\begin{figure}[!htb]
\center{\includegraphics[width=\textwidth]{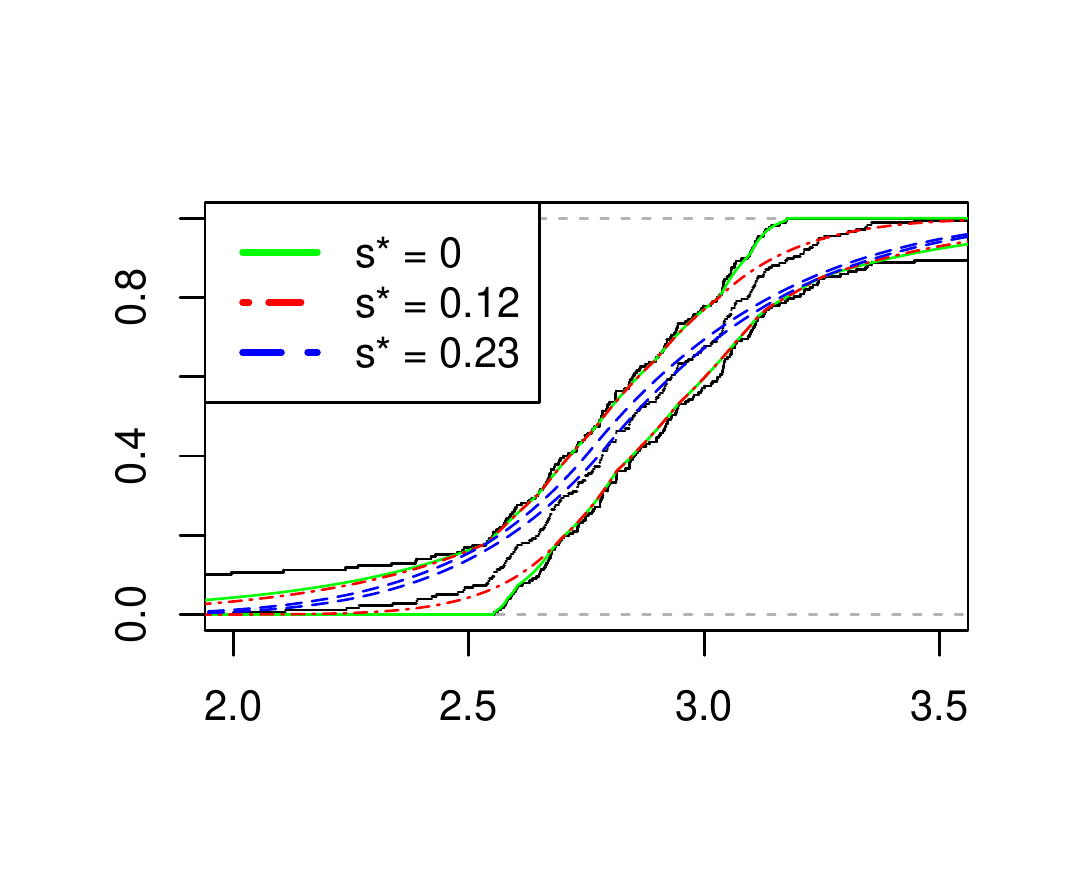}}
\caption{
Confidence bands from an initial KS band for the CEO salary data.
The step function in the middle is the empirical distribution function.
The two gray-black lines give the KS band and lines in other colors are refined 
confidence bands under the bi-$s^*$-concave assumption.
}
\label{plot:KS5c-apdx} 
\end{figure}

\subsection{Bi-$s^*$-concave bands based on WKS bands: application to the CEO salary data}

\begin{figure}[!htb]
\center{\includegraphics[width=\textwidth]{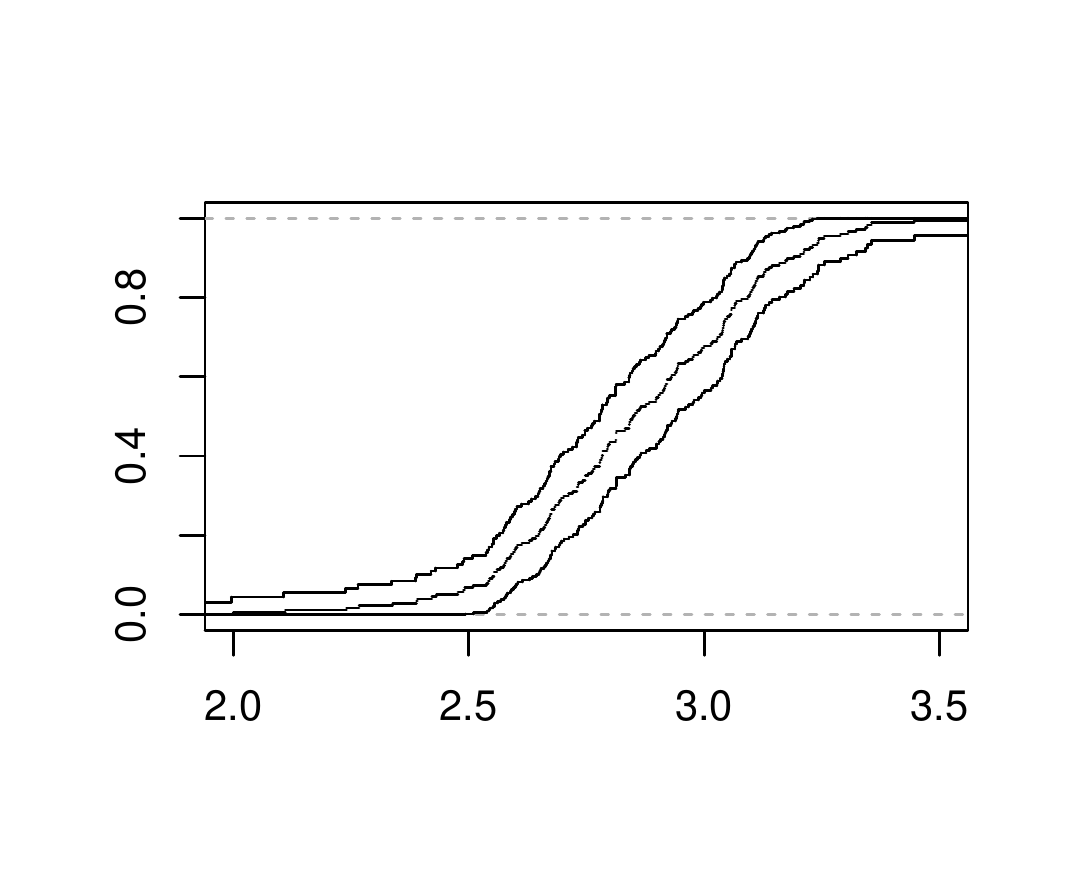}}
\caption{
WKS confidence bands for the CEO salary data.
The step function in the middle is the empirical distribution function.
The upper and lower lines give the WKS bands 
}
\label{plot:WKS6a-apdx} 
\end{figure}

\begin{figure}[!htb]
\center{\includegraphics[width=\textwidth]{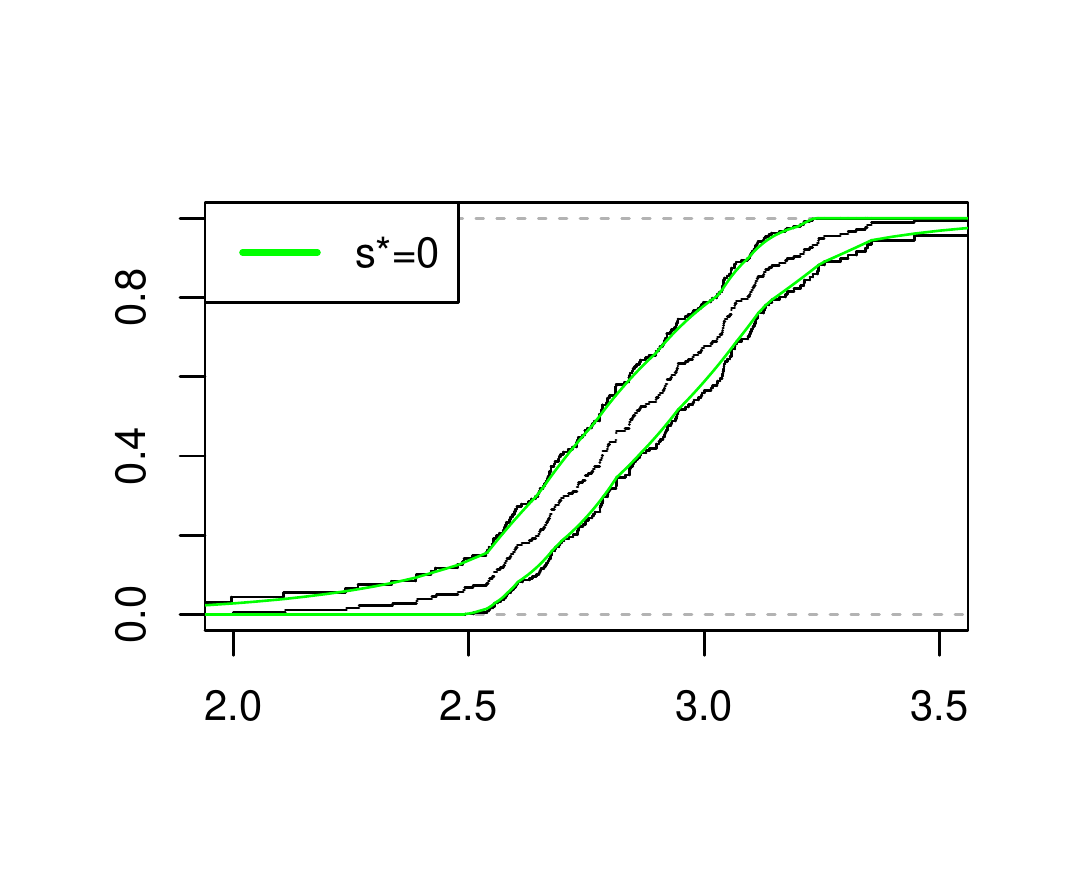}}
\label{plot:WKS-6b-apdx} 
\caption{
WKS confidence bands for the CEO salary data.
The step function in the middle is the empirical distribution function.
The upper and lower lines give the WKS bands and the green lines  are  the bi-$s^*$-concave bands for $s^*=0$.
}
\label{plot:WKS6b-apdx} 
\end{figure}

\begin{figure}[!htb]
\center{\includegraphics[width=\textwidth]{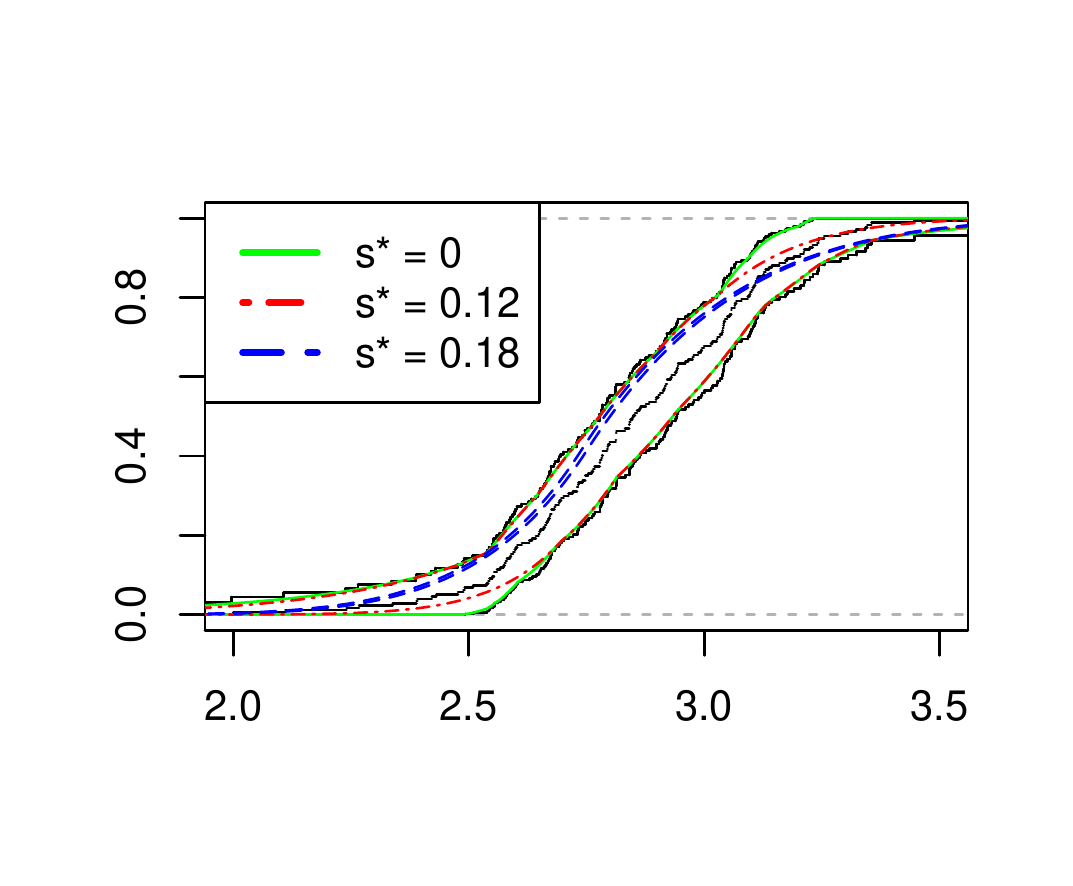}} 
\caption{
Confidence bands from an initial WKS band for the CEO salary data.
The step function in the middle is the empirical distribution function.
The two gray-black lines give the WKS band and lines in other colors are refined confidence bands under the bi-$s^*$-concave assumption.
}
\label{plot:WKS6c-apdx}
\end{figure}
\FloatBarrier

\section*{Acknowledgements:}   We owe thanks to Lutz D\"umbgen for several helpful discussions. 
We also thank two referees for their positive comments and suggestions.


\end{document}